\documentclass[10pt,a4paper]{article}

\usepackage[hidelinks]{hyperref}
\usepackage{lineno}

\usepackage{amssymb, bm}
\usepackage{amsmath}
\usepackage{amsthm}
\usepackage{mathrsfs}
\usepackage{siunitx}
\usepackage[margin=2cm]{geometry}
\usepackage{subfig}
\usepackage{float}
\usepackage{xcolor}

\usepackage{mathtools}
\usepackage{bbm}
\usepackage{multirow}
\usepackage{stmaryrd}

\usepackage{algorithm}     %serve per scrivere lo pseudo codice
\usepackage{algpseudocode} %serve per scrivere lo pseudo codice

\usepackage[utf8]{inputenc} % usually not needed (loaded by default)
\usepackage[T1]{fontenc}

\usepackage{tikz}
\usetikzlibrary{decorations.markings}
\usetikzlibrary{shapes.geometric}
\usetikzlibrary{shapes.arrows}
\usetikzlibrary{fit,calc}
\usetikzlibrary{quotes,angles}
\usetikzlibrary{arrows.meta}
\usepackage{animate}  
\usepackage{graphicx}
\usepackage{pgf}

\usepackage{pgfplots}
\usepackage[normalem]{ulem}

\newcommand{\vertiii}[1]{{\left\vert\kern-0.25ex\left\vert\kern-0.25ex\left\vert #1 
   \right\vert\kern-0.25ex\right\vert\kern-0.25ex\right\vert}}

\usepackage{geometry}
\makeatletter
\def\blfootnote{\xdef\@thefnmark{$\star$}\@footnotetext}
\makeatother
\newenvironment{Authors}%
  {\begin{center}\begin{bfseries}}%
  {\end{bfseries}\end{center}}
\newenvironment{Addresses}%
  {\begin{flushleft}\begin{itshape}}%
  {\end{itshape}\end{flushleft}}
  \newcommand{\email}[1]{\hspace*{\stretch{1}}\emph{\texttt{#1}}}

 \usepackage{fancyhdr}  
  \fancypagestyle{plain}{
\fancyhead{}
 }

\newtheorem{theorem}{Theorem}[section]

\newtheorem{remark}[theorem]{Remark}
\newtheorem*{remark*}{Remark}  % Define an unnumbered remark environment

\begin{document}

\title{Interface-aware dynamic mode decomposition with control for partitioned multiphysics systems}

\date{}

 \maketitle
\vspace{-50pt} 

\begin{Authors}
Jiarui Tian$^{1}$,
Dunhui Xiao$^{1}$,
Lei Zhang$^1$.
\end{Authors}

\begin{Addresses}
$^1$
School of Mathematical Sciences, Key Laboratory of Intelligent Computing and Applications (Ministry of Education), Tongji University, Shanghai 200092, China,
\\ \email{2433970@tongji.edu.cn, xiaodunhui@tongji.edu.cn, 22210@tongji.edu.cn} \\ \email{} \\[3mm]

\end{Addresses}

\begin{abstract}
We develop an interface-aware dynamic mode decomposition with control (DMDc) formulation for multiphysics systems, within an optimization-based partitioned coupling framework. The finite-element interface-input operator is projected onto the reduced space and kept fixed during identification, while the unknown state-dependent local dynamics are learned from trajectory data. The common interface control is recomputed online by minimizing the mismatch between the reconstructed physical interface quantities. The formulation accommodates local subproblems governed by ordinary differential equations (ODEs) and differential--algebraic equations (DAEs). For a DAE subproblem, tangent projection identifies the state dynamics, whereas a residualized normal regression identifies the pressure balance and provides an algebraic pressure reconstruction. We further propose regularized quadratic correction DMDc to enrich the tangent dynamics and, for DAEs, the normal pressure balance through jointly fitted linear and compact quadratic operators. Feature normalization and block-dependent regularization control the enlarged regression while retaining the prescribed interface operators and identified DAE geometry. The proposed framework is evaluated on problems involving both ODE--ODE and DAE--ODE couplings. The numerical results demonstrate accurate predictions and stable long-time behavior. 
\end{abstract}

\noindent
\emph{Keywords:} 
partitioned multiphysics; reduced-order modeling; dynamic mode decomposition with control; operator inference; fluid--structure interaction.

\section{Introduction}
\label{sec:introduction}

%\subsection{Reduced modeling of partitioned multiphysics problems}

Partitioned simulations of multiphysics systems require both the evolution of local physical fields and the enforcement of transmission conditions at their interfaces. In fluid–structure interaction, for example, the fluid and structural subproblems exchange velocities and tractions while satisfying their respective governing equations. 
Such simulations arise in cardiovascular mechanics, aeroelasticity, and many related applications \cite{hou2012numerical,bazilevs2013computational}, but their repeated use in many-query or real-time settings is often prohibitively expensive. Model order reduction (MOR) seeks to alleviate this cost by replacing high-dimensional
discretized problems with low-dimensional surrogates
\cite{benner2015survey}.  
Projection-based reduced models have been developed for both monolithic and partitioned formulations of coupled physical systems \cite{kalashnikova2013stable,nonino2023projection}, whereas non-intrusive approaches infer reduced dynamics from data without requiring access to full-order operators \cite{xiao2016nonintrusive}. Recent data-driven have further considered kinematics-consistent reduced dynamics \cite{dupont2023kinematics}, partitioned deep-learning surrogates \cite{gupta2022hybridfsi}, and non-intrusive partitioned ROM--FOM coupling
\cite{tiba2024partitioned}.
For partitioned multiphysics systems, however, reducing the local equations is only part of the task. 
The reduced models must also preserve the mechanism through which the subproblems interact. 
In particular, an accurate approximation of each isolated trajectory does not guarantee an accurate coupled trajectory if the interface quantities are represented incorrectly. 
This issue is especially important when an exchanged flux or traction has a prescribed physical action on each subsystem and must simultaneously satisfy transmission conditions. If the local dynamics and interface forcing are identified together without further structure, their respective effects can become confounded in the learned model. We therefore seek a reduced formulation in which the unknown local dynamics are identified from data while the interface control and the prescribed operator through which it acts on each subproblem remain explicit.

Optimization-based partitioned coupling provides a suitable architecture for this purpose. In this formulation, a common interface flux or traction is introduced as a control variable, the local physical problems are imposed as equality constraints, and continuity of the primary interface quantity is enforced by minimizing a mismatch objective \cite{kuberry2013decoupling,kuberry2015analysis}. This viewpoint has subsequently been combined with component-based MOR \cite{taddei2024nonoverlap} and projection-based FSI reduction \cite{taddei2025optimization}. The present work retains this coupling architecture but replaces the expensive state-dependent local operators by data-driven surrogates. The common interface control remains an online unknown determined by the coupled optimization problem; only the local state-dependent dynamics are identified from trajectory data. 
In this respect, the formulation also differs from learning a monolithic coupled evolution operator, while sharing the broader objective of exploiting known block structure in data-driven multiphysics reduction \cite{zastrow2026block}.

Dynamic mode decomposition (DMD) provides a data-driven representation of dynamical systems from snapshot data \cite{schmid2010dynamic,tu2014dynamic}. Its connection with Koopman operator theory gives a spectral interpretation of the identified dynamics \cite{rowley2009spectral}, and DMD and its variants are now widely used both for modal analysis and for reduced-order prediction \cite{schmid2022variants}. 
 For systems subject to external inputs, dynamic mode decomposition with control (DMDc) separates the state-dependent evolution from the action of an input \cite{proctor2016dmdc}; related input--output and controlled-Koopman formulations have been developed in \cite{benner2018iodmd,proctor2018koopman,korda2018linear}. This separation is particularly relevant to partitioned multiphysics systems, but with an important distinction: the input in a partitioned multiphysics problem is not an arbitrary exogenous signal. It is the interface flux or traction determined jointly by the coupled problem. 
 Moreover, the operator that maps this interface control into each local semi-discrete equation is known from the spatial discretization. 
 Consequently, re-identifying this input operator together with the local dynamics discards available physical information. Our starting point is therefore a known-input DMDc formulation in which the projected interface-input operator is prescribed and identification is restricted to the unknown local dynamics. This use of available structure is consistent with the broader use of known physical structure to restrict data-driven operators \cite{baddoo2023pidmd,morandin2023phdmd}. 
 
 A further difficulty arises because the local semi-discrete models need not all have the same mathematical structure. For example, diffusion, and structural subproblems lead to ordinary differential equations (ODEs); incompressible flow, insteads, gives a differential--algebraic equation (DAE), in which the velocity evolves dynamically while pressure acts as an algebraic multiplier. 
More generally, MOR for descriptor systems requires explicit treatment of the differential and algebraic components \cite{benner2017daesurvey}. 
 Treating pressure as an additional unconstrained evolving state would destroy this distinction. 
 Data-driven reduction of related structures has been considered through operator inference for incompressible flows, where velocity dynamics are separated from pressure reconstruction \cite{benner2022operator}, and through non-intrusive reduction of DAEs arising from lifting transformations \cite{khodabakhshi2022dae}.
In the present partitioned setting, the pressure-related reduced geometry therefore must be identified from trajectory data while remaining compatible with a prescribed interface-input operator and the optimization-based partitioned coupling. 

To address these requirements, we propose {interface-aware DMD with control} (IA-DMDc). The finite-element interface operator is projected onto the reduced spaces and kept fixed during identification. For an ODE subproblem, IA-DMDc therefore learns only the reduced state-dependent dynamics. For a DAE subproblem, admissible velocity snapshots define the differential, or tangent, reduced coordinates. 
Projecting the momentum balance onto this tangent space removes the pressure contribution and consequently provides no equation from which pressure can be recovered. 
 We therefore identify a complementary pressure-active normal subspace from the residualized momentum data and use the resulting normal balance to reconstruct the identifiable pressure algebraically.
 The reduced DAE consequently retains a descriptor interpretation: its finite dynamics are carried by the tangent state, whereas the pressure remains associated with the algebraic constraint. The same common interface control is then recomputed at every time step through the  optimization-based coupling objective.

A linear state operator can capture the dominant response while leaving systematic nonlinear state dependence unresolved. We therefore further propose regularized quadratic correction DMDc (RQC-DMDc) to address this systematic model defect. The construction of this quadratic correction builds on operator inference, which provides a framework for learning reduced vector fields with prescribed polynomial structure \cite{peherstorfer2016operator,kramer2024learning,qian2022pdes}. The number of coefficients, however, grows rapidly when quadratic terms are introduced. 
Regularization is therefore particularly important when the available training trajectories provide only limited coverage of the reduced state space  \cite{mcquarrie2021regularised,sawant2023physics}. 
 RQC-DMDc follows the regularized operator-inference paradigm, and enriches both the tangent dynamics and, for a DAE subproblem, the normal pressure balance with compact quadratic state features while retaining the prescribed interface operators and the identified DAE geometry. Feature normalization and block-dependent regularization are used to favor a linear-dominant model while allowing nonlinear defects to be represented.

%The method is designed to improve predictive accuracy; unlike structure-preserving quadratic formulations \cite{koike2024energy,goyal2025stablequadratic}, it does not impose a global nonlinear stability or energy condition.

The main contributions are threefold. First, we formulate an interface-aware DMDc framework for optimization-based partitioned coupling, retaining the prescribed interface-input operators during identification and recomputing the common interface control online. Second, we develop a common treatment of ODE and DAE subproblems in which residualized normal regression provides algebraic reconstruction without assigning independent dynamics to algebraic variable. Third, we introduce a normalized, block-regularized quadratic enrichment on the same interface and DAE architecture and incorporate it into the implicit online coupling formulation.

The remainder of the paper is organized as follows. Section~\ref{sec:optimization_framework}  introduces the optimization-based full-order coupling formulation and its spatial and temporal discretization.  
 Section~\ref{sec:interface_aware_dmdc} develops IA-DMDc for ODE and DAE subproblems and the corresponding online coupling formulation.  
Section~\ref{sec:regularized_quadratic_dmdc} presents RQC-DMDc and its online solution.
 Section~\ref{sec:numerical_experiments} reports the numerical experiments. Section~\ref{sec:conclusions} summarizes the main conclusions.

\section{Optimization-based coupling and finite-element setting}
\label{sec:optimization_framework}
This section establishes the full-order formulation underlying the reduced
models. The local subproblems are coupled through an optimization-based
interface control, discretized in space by finite elements, and advanced in time using a
\(\theta\)-scheme.

\paragraph{Strong formulation}
We consider a multiphysics problem posed on a domain \(\Omega\) partitioned
into two non-overlapping subdomains, \(\Omega_1\) and \(\Omega_2\). The
interface between the two subdomains is
\[
\Gamma:=\partial\Omega_1\cap\partial\Omega_2.
\]
For \(i=1,2\), let \(u_i\) denote the time-dependent state field on
\(\Omega_i\) and let \(p_i\), when present, denote an algebraic multiplier.
We represent each controlled local problem abstractly as
\begin{equation}
\label{eq:subproblem}
\begin{cases}
\partial_t u_i=\mathcal F_i(u_i,p_i;x,t),
&x\in\Omega_i,\quad t\in(0,T),\\
\mathcal C_i(u_i)=0,
&x\in\Omega_i,\quad t\in(0,T),\\
\mathcal B_i(u_i,p_i;x,t)=0,
&x\in\partial\Omega_i\setminus\Gamma,\quad t\in(0,T),\\
u_i(x,0)=u_{i,0}(x),
&x\in\Omega_i.
\end{cases}
\end{equation}
In the applications considered in this work, the diffusion and Allen--Cahn
subproblems have no algebraic multiplier; equivalently,
\(\mathcal C_i\equiv0\), i.e., the second relation in
\eqref{eq:subproblem} is absent. In the incompressible Stokes and
Navier--Stokes subproblems, \(p_i\) is the pressure and the constraint is
\(
\mathcal C_i(u_i)=\nabla\!\cdot u_i=0
\).

At the interface \(\Gamma\), the two subproblems are coupled through two 
 transmission conditions: continuity
of a primary interface variable,
\begin{equation}
\label{eq:interface_continuity}
\mathcal R_{\Gamma,1}(u_1,x,t)
=
\mathcal R_{\Gamma,2}(u_2,x,t);
\end{equation}
and balance of an interface flux,
\begin{equation}
\label{eq:interface_balance}
\mathcal T_1(u_1,p_1;x,t)+\mathcal T_2(u_2,p_2;x,t)=0.
\end{equation}

We solve the coupled problem using the optimization-based coupling method
described in
\cite{kuberry2013decoupling,kuberry2015analysis,taddei2025optimization}. The central idea is to treat the balanced flux  as a common interface
control while minimizing the mismatch in the primary interface variable.
Accordingly, we define
\[
g(x,t)=\mathcal T_1(u_1,p_1;x,t)=-\mathcal T_2(u_2,p_2;x,t),
\]
and set \(\sigma_1=1\) and \(\sigma_2=-1\). The two local Neumann-type
conditions can then be written as
\begin{equation}
\label{eq:subproblems}
\mathcal T_i(u_i,p_i;x,t)=\sigma_i g(x,t),
\qquad x\in\Gamma,\quad i=1,2.
\end{equation}
Thus, for any prescribed interface control \(g\), the flux or traction
balance is satisfied by construction, since
\[
\mathcal T_1+\mathcal T_2
=(\sigma_1+\sigma_2)g=0.
\]
Moreover, the conditions in \eqref{eq:subproblems} provide the
Neumann-type boundary data required to close the two local subproblems. 
An arbitrary choice of \(g\), however, does not generally ensure continuity
of the primary interface quantity. We therefore determine \(g\) by minimizing
the interface mismatch, together with a regularization term:
\begin{equation}
\label{eq:optimization}
J(u_1,u_2,g;t)
:=
\frac12\int_\Gamma
\left\|
\mathcal R_{\Gamma,1}(u_1)
-
\mathcal R_{\Gamma,2}(u_2)
\right\|^2\,\mathrm ds
+
\frac{\lambda_{\mathrm g}}2
\int_\Gamma\|g\|^2\,\mathrm ds,
\end{equation}
where \(\lambda_{\mathrm g}\geq0\) is the control-regularization parameter.

The coupled problem can thus be formulated as the constrained optimization
problem
\begin{equation}
\label{eq:continuous_coupled_system}
\begin{aligned}
\min_{u_1,u_2,g}\quad
&J(u_1,u_2,g;t), \quad
\text{subject to the two local controlled subproblems \eqref{eq:subproblem}.}
\end{aligned}
\end{equation}

\paragraph{Finite-element formulation}
We use finite elements to construct the high-fidelity spatial discretization.
Let \(g(t)\in\mathbb R^m\) collect the discrete interface-control
coefficients, let \(\widetilde U_{i,h}(t)\in\mathbb R^{N_i}\) be the
 finite-element coefficient vector of the dynamic state, and let
\(P_i(t)\in\mathbb R^{q_i}\), when present, be the coefficient vector of
the algebraic multiplier. We denote the  finite-element mass,
constraint, and interface-load matrices by \(\widetilde M_i\),
\(\widetilde C_{i,h}\), and \(\widetilde B_{i,h}\), respectively, and the
assembled state-dependent term by \(\widetilde F_{i,h}\). Thus, 
finite-element discretization of the local problem
\eqref{eq:subproblem}, subject to the controlled interface condition
\eqref{eq:subproblems}, yields the saddle-point system
\[
\widetilde M_i\dot{\widetilde U}_{i,h}
=\widetilde F_{i,h}(\widetilde U_{i,h},t)
+\widetilde C_{i,h}^\top P_i
+\sigma_i\widetilde B_{i,h}g,
\qquad
\widetilde C_{i,h}\widetilde U_{i,h}=0.
\]
Here we use $\widetilde{(\cdot)}$ to denote physical variables. To formulate 
POD construction and the tangent--normal
decomposition later using Euclidean orthogonality, we introduce mass-orthonormal coordinates. Specifically, we factor the mass matrix \(\widetilde M_i=R_i^\top R_i\), and use untilded symbol as variable in mass-orthonormal coordinates, for example $U_i=R_i\widetilde U_{i,h}$, we obtain 
the differential-algebraic equation (DAE)
%define
%\[
%U_i=R_i\widetilde U_{i,h},
%\qquad
%C_i=\widetilde C_{i,h}R_i^{-1},
%\qquad
%B_i=R_i^{-\top}\widetilde B_{i,h},
%\qquad
%F_i(U_i,t)
%=R_i^{-\top}\widetilde F_{i,h}(R_i^{-1}U_i,t).
%\]
%Henceforth, a tilde identifies a vector or operator expressed in the original
%finite-element coordinates, whereas the corresponding untilded symbol is
%expressed in mass-orthonormal coordinates. The multiplier coefficients
%\(P_i\) and the physical interface-control coefficients \(g\) do not belong
%to the transformed dynamic-state space and therefore remain untilded.
%Multiplying the semi-discrete evolution equation by \(R_i^{-\top}\) gives
\begin{equation}
\label{eq:semi_discrete_dae_balanced}
\begin{aligned}
\dot U_i
&=F_i(U_i,t)+C_i^\top P_i+\sigma_iB_i g,\\
0&=C_iU_i.
\end{aligned}
\end{equation}

If the spatially discretized PDE has no algebraic constraint, the pressure
and constraint blocks are removed and
\eqref{eq:semi_discrete_dae_balanced} reduces instead to the ordinary
differential equation (ODE)
\begin{equation}
\label{eq:semi_discrete_ode_balanced}
\dot U_i=F_i(U_i,t)+\sigma_iB_i g.
\end{equation}
The structural subproblems considered in this work are governed by
second-order equations for the displacement \(\widetilde d_{s,h}\).
Introducing the velocity
\(\widetilde v_{s,h}=\dot{\widetilde d}_{s,h}\) and the augmented state
\(\widetilde U_{s,h}
=[\widetilde d_{s,h}^\top,\widetilde v_{s,h}^\top]^\top\)
converts each structural problem into a first-order system. Applying the same
mass-orthonormal transformation to this augmented system yields an ODE of the
form \eqref{eq:semi_discrete_ode_balanced}.

\paragraph{Fully discrete formulation}
To generate the full-order trajectories used for training, we discretize
\eqref{eq:semi_discrete_dae_balanced} and
\eqref{eq:semi_discrete_ode_balanced} using the one-step
\(\theta\)-scheme. Let
\[
t^n=n\Delta t,
\qquad n=0,1,\ldots,N,
\qquad \Delta t=T/N,
\qquad t^{n+\theta}=(1-\theta)t^n+\theta t^{n+1},
\]
and fix \(\theta\in[1/2,1]\). For the dynamic state and interface control,
define
\[
U_i^{n+\theta}=(1-\theta)U_i^n+\theta U_i^{n+1},
\qquad
g^{n+\theta}=(1-\theta)g^n+\theta g^{n+1}.
\]

For a DAE subproblem, the one-step \(\theta\)-scheme gives
\begin{equation}
\label{eq:fully_discrete_dae_i}
\begin{aligned}
\frac{U_i^{n+1}-U_i^n}{\Delta t}
&=
F_i(U_i^{n+\theta},t^{n+\theta})
+C_i^\top P_i^{n+\theta}
+\sigma_iB_i g^{n+\theta},\\
0&=C_iU_i^{n+1}.
\end{aligned}
\end{equation}
Here \(P_i^{n+\theta}\) is an algebraic unknown evaluated at
\(t^{n+\theta}\); it is not advanced by a separate evolution
equation; and the initial state is assumed to be consistent:
\(C_iU_i^0=0\). 
For an ODE subproblem, the corresponding time-discrete equation is
\begin{equation}
\label{eq:fully_discrete_ode_i}
\frac{U_i^{n+1}-U_i^n}{\Delta t}
=
F_i(U_i^{n+\theta},t^{n+\theta})
+\sigma_iB_i g^{n+\theta}.
\end{equation}

At time \(t^{n+1}\), the coupled problem is
\begin{equation}
\label{eq:discrete_optimization}
\begin{aligned}
\min_{\substack{
U_1^{n+1},\,U_2^{n+1},\,g^{n+1},\\
P_i^{n+\theta}\ \text{for DAE subproblems}}}
\quad
J^{n+1}
:={}&
\frac12
\left\|
{\mathcal R}_{\Gamma,1,h}(U_1^{n+1})
-
{\mathcal R}_{\Gamma,2,h}(U_2^{n+1})
\right\|_{W_\Gamma}^2
+
\frac{\lambda_{\mathrm g}}2\|g^{n+1}\|_{W_{\mathrm g}}^2,
\\
\text{subject to}\quad
&\eqref{eq:fully_discrete_dae_i}
\quad\text{for each DAE subproblem},\\
&\eqref{eq:fully_discrete_ode_i}
\quad\text{for each ODE subproblem}.
\end{aligned}
\end{equation}
Here \(\|v\|_W^2=v^\top Wv\), with symmetric positive-definite weight
matrices \(W_\Gamma\) and \(W_{\mathrm g}\). The algebraic multiplier is
internal to a local subproblem and does not enter the interface objective.
Equations
\eqref{eq:fully_discrete_dae_i}--\eqref{eq:discrete_optimization} define the
mass-orthonormal full-order trajectories from which the reduced models below
are learned.

\section{Interface-aware DMDc for ODE and DAE subproblems}
\label{sec:interface_aware_dmdc}

We now construct reduced models for the local systems introduced in Section~\ref{sec:optimization_framework}. The objective is to learn the state-dependent dynamics while preserving two features of the full-order coupling formulation: the known interface-input operator and, for a constrained subproblem, the algebraic pressure balance. 
Preserving the latter is essential to the DAE formulation: whereas an ODE requires only an evolution equation for the reduced state, a DAE must additionally enforce the constraint and determine the associated multiplier algebraically, rather than evolve it as an independent state variable.

The construction consists of an offline identification phase and an online prediction phase. During the offline phase, synchronized state, pressure, and interface-control snapshots are used to identify the local reduced operators. During the online phase, the identified local models are coupled by recomputing a common interface control at every time step. We first present the construction for a DAE subproblem; the ODE formulation then follows by omitting the normal and pressure blocks. Throughout this section, all state snapshots and interface operators are expressed in the mass-orthonormal coordinates of Section~\ref{sec:optimization_framework}. 

Consider a single local subproblem and temporarily omit its index. In the momentum equation of \eqref{eq:semi_discrete_dae_balanced}, linear IA-DMDc approximates only the unknown state-dependent term by \(F(U,t)\approx AU\), while retaining the prescribed interface action \(Bg\). We denote the pressure-action operator by \(G_{\mathrm p}\). It equals \(C^\top\) when the full constraint operator is available, whereas the non-intrusive construction below identifies only its trajectory-active component. These choices give the momentum ansatz
\begin{equation}
\label{eq:data_level_dae}
\dot U\approx AU+G_{\mathrm p}P+\sigma Bg.
\end{equation}
Equation~\eqref{eq:data_level_dae} represents the momentum balance only; the
algebraic constraint is recovered below through a normal-force balance. 

\subsection{Snapshot coordinates and constraint geometry}
\label{subsec:dmdc_projected_interface}
\label{subsec:common_dae_geometry_identification}

\paragraph{POD coordinates and the prescribed interface action}
Let \(N_{\mathrm{tr}}+1\) synchronized state and interface-control snapshots
be available at \(t^n=n\Delta t\), \(n=0,\ldots,N_{\mathrm{tr}}\). For a
DAE, we also collect the \(N_{\mathrm{tr}}\) pressure snapshots returned by
the time integrator at the evaluation points \(t^{n+\theta}\):
\[
S_U=
\begin{bmatrix}U^0&U^1&\cdots&U^{N_{\mathrm{tr}}}\end{bmatrix},
\qquad
S_{\mathrm g}=
\begin{bmatrix}g^0&g^1&\cdots&g^{N_{\mathrm{tr}}}\end{bmatrix},
\]
\[
S_P^\theta=
\begin{bmatrix}
P^\theta&P^{1+\theta}&\cdots&P^{N_{\mathrm{tr}}-1+\theta}
\end{bmatrix}.
\]
Let \(\Phi_{\mathrm u}\), \(\Phi_{\mathrm p}\), and
\(\Phi_{\mathrm g}\) be POD bases~\cite{sirovich1987turbulence} for the
state, pressure, and interface control, respectively, and 
\(\alpha_\bullet\) the corresponding coordinate vector. We write
\begin{equation}
\label{eq:dae_pod_coordinates}
U^n\approx\Phi_{\mathrm u}\alpha_{\mathrm u}^n,
\qquad
P^{n+\theta}\approx\Phi_{\mathrm p}\alpha_{\mathrm p}^{n+\theta},
\qquad
g^n\approx\Phi_{\mathrm g}\alpha_{\mathrm g}^n,
\end{equation}
The bases have orthonormal columns, with
\(\Phi_{\mathrm u}\in\mathbb R^{N_{\mathrm u}\times r_{\mathrm u}}\),
\(\Phi_{\mathrm p}\in\mathbb R^{N_{\mathrm p}\times r_{\mathrm p}}\), and
\(\Phi_{\mathrm g}\in\mathbb R^{N_{\mathrm g}\times r_{\mathrm g}}\),
where \(N_{\mathrm u}\), \(N_{\mathrm p}\), and \(N_{\mathrm g}\) are the
corresponding full-order dimensions.
The corresponding reduced coordinates are
\begin{equation}
\label{eq:dae_reduced_coordinates}
\alpha_{\mathrm u}^n=\Phi_{\mathrm u}^\top U^n,
\qquad
\alpha_{\mathrm p}^{n+\theta}
=\Phi_{\mathrm p}^\top P^{n+\theta},
\qquad
\alpha_{\mathrm g}^n=\Phi_{\mathrm g}^\top g^n.
\end{equation}
The pressure coordinate in \eqref{eq:dae_reduced_coordinates} is obtained from
the value returned by the full-order solver at \(t^{n+\theta}\); it is not
formed by interpolating endpoint pressures. The state and control coordinates
at the same evaluation point are
\[
\alpha_{\mathrm u}^{n+\theta}
=(1-\theta)\alpha_{\mathrm u}^n+\theta\alpha_{\mathrm u}^{n+1},
\qquad
\alpha_{\mathrm g}^{n+\theta}
=(1-\theta)\alpha_{\mathrm g}^n+\theta\alpha_{\mathrm g}^{n+1}.
\]

For later regression, define the difference-quotient data and the reduced
coordinate matrices
\begin{equation}
\label{eq:dae_stage_data}
\begin{aligned}
\dot U_\Delta
&:=
\frac1{\Delta t}
\begin{bmatrix}
U^1-U^0&\cdots&U^{N_{\mathrm{tr}}}-U^{N_{\mathrm{tr}}-1}
\end{bmatrix},\\
\mathcal S_{\mathrm u}
&:=
\begin{bmatrix}
\alpha_{\mathrm u}^\theta&
\alpha_{\mathrm u}^{1+\theta}&\cdots&
\alpha_{\mathrm u}^{N_{\mathrm{tr}}-1+\theta}
\end{bmatrix},\\
\mathcal S_{\mathrm p}
&:=
\begin{bmatrix}
\alpha_{\mathrm p}^\theta&
\alpha_{\mathrm p}^{1+\theta}&\cdots&
\alpha_{\mathrm p}^{N_{\mathrm{tr}}-1+\theta}
\end{bmatrix},\\
\mathcal S_{\mathrm g}
&:=
\begin{bmatrix}
\alpha_{\mathrm g}^\theta&
\alpha_{\mathrm g}^{1+\theta}&\cdots&
\alpha_{\mathrm g}^{N_{\mathrm{tr}}-1+\theta}
\end{bmatrix}.
\end{aligned}
\end{equation}
The known interface action in the reduced
control coordinates is
\begin{equation}
\label{eq:reduced_interface_action}
B_{\mathrm g}:=B\Phi_{\mathrm g}.
\end{equation}

\paragraph{Why a normal space is required}
The columns of \(\Phi_{\mathrm u}\) are constructed from admissible state
snapshots and therefore approximately satisfy
\(C\Phi_{\mathrm u}=0\). Since the pressure action belongs to
\(\operatorname{range}(C^\top)\), it is orthogonal to the reduced state
space:
\[
\Phi_{\mathrm u}^\top G_{\mathrm p}\Phi_{\mathrm p}\approx0.
\]
Consequently, substituting
\(U\approx\Phi_{\mathrm u}\alpha_{\mathrm u}\),
\(P\approx\Phi_{\mathrm p}\alpha_{\mathrm p}\), and
\(g\approx\Phi_{\mathrm g}\alpha_{\mathrm g}\) into
\eqref{eq:data_level_dae} and projecting onto \(\Phi_{\mathrm u}\)
eliminates the pressure term.
This tangent projection therefore determines the reduced-state evolution but
provides no equation for the pressure. To recover the pressure, we therefore introduce
a reduced normal basis \(\Phi_{\mathrm n}\), with
\(\Phi_{\mathrm u}^\top\Phi_{\mathrm n}=0\), and use the normal component of the momentum balance
to construct an auxiliary algebraic equation for the pressure. 

The pressure-active normal basis \(\Phi_{\mathrm n}\) is constructed from
the snapshot data in the remainder of this subsection. 
Subsection~\ref{subsec:dmdc_prescribed_input} then uses the resulting
tangent--normal decomposition to identify the operators of the local
IA-DMDc model.

\paragraph{Identification of the normal pressure geometry}
We first remove
the known interface contribution from the momentum equation and project the
remainder onto \(\operatorname{range}(\Phi_{\mathrm u})^\perp\):
\begin{equation}
\label{eq:dae_input_corrected_data}
\mathcal Y:=\dot U_\Delta-\sigma B_{\mathrm g}\mathcal S_{\mathrm g},
\qquad
\mathcal Y_\perp
:=(I-\Phi_{\mathrm u}\Phi_{\mathrm u}^\top)\mathcal Y.
\end{equation}

Projection onto \(\operatorname{range}(\Phi_{\mathrm u})^\perp\) removes
only the component contained in the retained reduced state space. Because
\(\Phi_{\mathrm u}\) is truncated, its orthogonal complement contains not
only pressure-active normal directions but also tangent directions omitted
by the POD basis. Consequently, \(\mathcal Y_\perp\) contains both the
pressure action and a residual state contribution.

Along the training trajectory, we represent this residual state contribution
by
\[
(I-\Phi_{\mathrm u}\Phi_{\mathrm u}^{\top})
A\Phi_{\mathrm u}\mathcal S_{\mathrm u}
\approx
K_{\perp}\mathcal S_{\mathrm u},
\]
and denote the projected pressure operator by
\[
G_{\perp}
:=
(I-\Phi_{\mathrm u}\Phi_{\mathrm u}^{\top})
G_{\mathrm p}\Phi_{\mathrm p}.
\]
The projected snapshot data (applying \(I-\Phi_{\mathrm u}\Phi_{\mathrm u}^{\top}\) to \eqref{eq:data_level_dae}) therefore satisfy
\begin{equation}
\label{eq:linear_geometry_surrogate}
\mathcal Y_\perp
\approx
K_{\perp}\mathcal S_{\mathrm u}
+
G_{\perp}\mathcal S_{\mathrm p}.
\end{equation}
To isolate the pressure contribution, we eliminate the state-dependent term
\(K_{\perp}\mathcal S_{\mathrm u}\) in the temporal snapshot space. Let
\[
\Pi_{\mathrm u}^{\perp}
:=I_{N_{\mathrm{tr}}}
-\mathcal S_{\mathrm u}^\dagger\mathcal S_{\mathrm u}
\]
be the orthogonal projector onto the complement of the temporal row space of
\(\mathcal S_{\mathrm u}\). Since
\(\mathcal S_{\mathrm u}\Pi_{\mathrm u}^{\perp}=0\), right multiplication of
\eqref{eq:linear_geometry_surrogate} by \(\Pi_{\mathrm u}^{\perp}\) eliminates the
unknown term \(K_{\perp}\mathcal S_{\mathrm u}\) and gives
\begin{equation}
\label{eq:geometry_residualized_relation}
\mathcal Y_\perp \Pi_{\mathrm u}^{\perp}
\approx G_{\perp}\mathcal S_{\mathrm p}\Pi_{\mathrm u}^{\perp}.
\end{equation}
We therefore obtain the candidate pressure
action
\begin{equation}
\label{eq:geometry_G_candidate}
G_\star
:=\mathcal Y_\perp \Pi_{\mathrm u}^{\perp}
\bigl(\mathcal S_{\mathrm p}\Pi_{\mathrm u}^{\perp}\bigr)^\dagger.
\end{equation}

\begin{remark}[Practical computation of \(G_\star\)]
The dense temporal projector
\[
\Pi_{\mathrm u}^{\perp}
=
I_{N_{\mathrm{tr}}}
-
\mathcal S_{\mathrm u}^{\dagger}\mathcal S_{\mathrm u}
\]
need not be formed explicitly. Let
\begin{equation}
\label{eq:geometry_compact_svd}
\mathcal S_{\mathrm u}
=
\mathcal U_{\mathrm u}
\Sigma_{\mathrm u}
\mathcal W_{\mathrm u}^{\top}
\end{equation}
be a compact SVD of numerical rank \(\rho_{\mathrm u}\). Because the columns
of
\(\mathcal W_{\mathrm u}
\in\mathbb R^{N_{\mathrm{tr}}\times\rho_{\mathrm u}}\)
span the temporal row space of \(\mathcal S_{\mathrm u}\), the residualized
data can be computed as
\begin{equation}
\label{eq:geometry_G_candidate_implementation}
\begin{aligned}
\mathcal Y_{\mathrm{res}}
&=
\mathcal Y_\perp
-
(\mathcal Y_\perp\mathcal W_{\mathrm u})
\mathcal W_{\mathrm u}^{\top},
\\
\mathcal S_{\mathrm p,\mathrm{res}}
&=
\mathcal S_{\mathrm p}
-
(\mathcal S_{\mathrm p}\mathcal W_{\mathrm u})
\mathcal W_{\mathrm u}^{\top}.
\end{aligned}
\end{equation}
The candidate pressure-action operator is then evaluated as
\[
G_\star
=
\mathcal Y_{\mathrm{res}}
\mathcal S_{\mathrm p,\mathrm{res}}^{\dagger}.
\]
\end{remark}

Before truncation, we replace \(G_\star\) by
\((I-\Phi_{\mathrm u}\Phi_{\mathrm u}^\top)G_\star\) to remove roundoff
leakage into the tangent space. We then compute a rank-\(r_{\mathrm n}\)
truncated SVD,
\begin{equation}
\label{eq:identifiable_pressure_geometry}
G_\star
\approx\Phi_{\mathrm n}\Sigma_{\mathrm n}T_{\mathrm p}^\top,
\end{equation}
and define
\begin{equation}
\label{eq:learned_CG_pair}
\widehat G_{\mathrm p}:=\Phi_{\mathrm n}\Sigma_{\mathrm n},
\qquad
\widehat C:=\widehat G_{\mathrm p}^\top,
\qquad
\Phi_{\mathrm p,\mathrm{id}}:=\Phi_{\mathrm p}T_{\mathrm p},
\qquad
\alpha_{\mathrm p,\mathrm{id}}
:=T_{\mathrm p}^\top\alpha_{\mathrm p},
\qquad
\mathcal S_{\mathrm p,\mathrm{id}}
:=T_{\mathrm p}^\top\mathcal S_{\mathrm p}.
\end{equation}
The retained singular-value matrix
\(\Sigma_{\mathrm n}\in
\mathbb R^{r_{\mathrm n}\times r_{\mathrm n}}\) is nonsingular. Moreover,
\begin{equation}
\label{eq:learned_CG_identities}
\Phi_{\mathrm u}^\top\Phi_{\mathrm n}=0,
\qquad
\Phi_{\mathrm u}^\top\widehat G_{\mathrm p}=0,
\qquad
\widehat C\Phi_{\mathrm u}=0,
\qquad
\Phi_{\mathrm n}^\top\widehat G_{\mathrm p}=\Sigma_{\mathrm n}.
\end{equation}
Thus \(\Phi_{\mathrm n}\) supplies the learned normal directions, while
\(\Phi_{\mathrm p,\mathrm{id}}\) spans the associated identifiable pressure
space. Hence \(\operatorname{range}(\Phi_{\mathrm n})\) approximates only
the trajectory-active portion of \(\operatorname{range}(C^\top)\), rather
than the complete normal space. If
\(\mathcal S_{\mathrm p}\Pi_{\mathrm u}^{\perp}\) has zero numerical rank, the training
data contain no pressure component that can be distinguished from the chosen
linear state library.

\subsection{Offline identification of the local IA-DMDc model}
\label{subsec:dmdc_prescribed_input}

\paragraph{Tangent and normal regressions}
Once the normal geometry has been identified, the prescribed interface action
is projected onto the tangent and normal spaces:
\begin{equation}
\label{eq:dae_tangent_normal_input}
B_{\mathrm t}:=\Phi_{\mathrm u}^\top B_{\mathrm g},
\qquad
B_{\mathrm n}:=\Phi_{\mathrm n}^\top B_{\mathrm g}.
\end{equation}
Keeping these input operators fixed separates the known coupling physics from
the unknown local dynamics. The tangent generator is determined by the
known-input DMDc regression~\cite{proctor2016dmdc}
\begin{equation}
\label{eq:dmdc_tangent_regression}
\min_L
\left\|
\dot{\mathcal S}_{\mathrm u}
-L\mathcal S_{\mathrm u}
-\sigma B_{\mathrm t}\mathcal S_{\mathrm g}
\right\|_{\mathrm F}^2,
\qquad
\dot{\mathcal S}_{\mathrm u}
:=\Phi_{\mathrm u}^\top\dot U_\Delta,
\end{equation}
whose minimum-norm solution is
\begin{equation}
\label{eq:reduced_A_solution}
L=
\left(
\dot{\mathcal S}_{\mathrm u}
-\sigma B_{\mathrm t}\mathcal S_{\mathrm g}
\right)\mathcal S_{\mathrm u}^\dagger.
\end{equation}
To obtain the remaining normal operator, project the input-corrected momentum
data $\mathcal Y$ onto \(\Phi_{\mathrm n}\). Using
\(\Phi_{\mathrm n}^{\top}\Phi_{\mathrm u}=0\), \(\Phi_{\mathrm n}^{\top}\mathcal Y
=
\Phi_{\mathrm n}^{\top}\mathcal Y_\perp
\), 
\(G_\star\approx
\Phi_{\mathrm n}\Sigma_{\mathrm n}T_{\mathrm p}^\top\), and
\(\mathcal S_{\mathrm p,\mathrm{id}}
=T_{\mathrm p}^\top\mathcal S_{\mathrm p}\) gives
\[
\Phi_{\mathrm n}^\top\mathcal Y
\approx K_{\mathrm n}\mathcal S_{\mathrm u}
+\Sigma_{\mathrm n}\mathcal S_{\mathrm p,\mathrm{id}}.
\]
We therefore solve
\[
\min_{K_{\mathrm n}}
\left\|
\Phi_{\mathrm n}^\top\mathcal Y
-\Sigma_{\mathrm n}\mathcal S_{\mathrm p,\mathrm{id}}
-K_{\mathrm n}\mathcal S_{\mathrm u}
\right\|_{\mathrm F}^2.
\]
The minimum-norm solution is
\begin{equation}
\label{eq:dmdc_normal_regression}
K_{\mathrm n}
=
\left(
\Phi_{\mathrm n}^\top\mathcal Y
-\Sigma_{\mathrm n}\mathcal S_{\mathrm p,\mathrm{id}}
\right)\mathcal S_{\mathrm u}^\dagger.
\end{equation}
The resulting local DAE model is
\begin{subequations}
\label{eq:dmdc_tangent_normal_model}
\begin{align}
\dot{\alpha}_{\mathrm u}
&=L\alpha_{\mathrm u}+\sigma B_{\mathrm t}\alpha_{\mathrm g},
\label{eq:dmdc_tangent_model}\\
0&=K_{\mathrm n}\alpha_{\mathrm u}
+\Sigma_{\mathrm n}\alpha_{\mathrm p,\mathrm{id}}
+\sigma B_{\mathrm n}\alpha_{\mathrm g}.
\label{eq:dmdc_normal_model}
\end{align}
\end{subequations}
Only the tangent state coordinate \(\alpha_{\mathrm u}\) is advanced in
time. The second equation is an algebraic normal-force balance that
determines the identifiable pressure coordinate
\(\alpha_{\mathrm p,\mathrm{id}}\).

\paragraph{Descriptor interpretation and pressure recovery}
We now recast the reduced DAE
\eqref{eq:dmdc_tangent_normal_model} in descriptor form. 
This representation makes the differential--algebraic structure explicit and
will later allow stabilization to act only on the finite dynamics. Introduce
an auxiliary normal coordinate \(\alpha_{\mathrm n}\) through
\[
U\approx
\Phi_{\mathrm u}\alpha_{\mathrm u}
+\Phi_{\mathrm n}\alpha_{\mathrm n}.
\]
Since
\[
\widehat C
\left(
\Phi_{\mathrm u}\alpha_{\mathrm u}
+\Phi_{\mathrm n}\alpha_{\mathrm n}
\right)
=\Sigma_{\mathrm n}^\top\alpha_{\mathrm n},
\]
the learned constraint imposes
\(\Sigma_{\mathrm n}^\top\alpha_{\mathrm n}=0\). Equations
\eqref{eq:dmdc_tangent_normal_model} are equivalently embedded in
\begin{equation}
\label{eq:dmdc_reduced_dae_descriptor}
E_{\mathrm r}\frac{\mathrm d}{\mathrm dt}
\begin{bmatrix}
\alpha_{\mathrm u}\\
\alpha_{\mathrm n}\\
\alpha_{\mathrm p,\mathrm{id}}
\end{bmatrix}
=
A_{\mathrm r}^{\mathrm{IA}}
\begin{bmatrix}
\alpha_{\mathrm u}\\
\alpha_{\mathrm n}\\
\alpha_{\mathrm p,\mathrm{id}}
\end{bmatrix}
+\sigma B_{\mathrm r}^{\mathrm{DAE}}\alpha_{\mathrm g},
\end{equation}
where
\begin{equation}
\label{eq:dmdc_reduced_dae_matrices}
E_{\mathrm r}=
\begin{bmatrix}
I_{r_{\mathrm u}}&0&0\\
0&I_{r_{\mathrm n}}&0\\
0&0&0
\end{bmatrix},
\qquad
A_{\mathrm r}^{\mathrm{IA}}=
\begin{bmatrix}
L&0&0\\
K_{\mathrm n}&0&\Sigma_{\mathrm n}\\
0&\Sigma_{\mathrm n}^\top&0
\end{bmatrix},
\qquad
B_{\mathrm r}^{\mathrm{DAE}}=
\begin{bmatrix}B_{\mathrm t}\\B_{\mathrm n}\\0\end{bmatrix}.
\end{equation}
Because \(\Sigma_{\mathrm n}\) is nonsingular, the constraint gives
\(\alpha_{\mathrm n}=0\); this coordinate is introduced only to display the
descriptor structure and is eliminated in the online solve. The pressure coordinate and the corresponding
full-dimensional pressure approximation are
\begin{equation}
\label{eq:dmdc_pressure_recovery}
\alpha_{\mathrm p,\mathrm{id}}
=-\Sigma_{\mathrm n}^{-1}
\left(
K_{\mathrm n}\alpha_{\mathrm u}
+\sigma B_{\mathrm n}\alpha_{\mathrm g}
\right),
\qquad
P\approx
\Phi_{\mathrm p,\mathrm{id}}\alpha_{\mathrm p,\mathrm{id}}.
\end{equation}

\paragraph{ODE formulation}
For an unconstrained subproblem, pressure snapshots and the normal geometry
are unnecessary. We use the same symbol \(\Phi_{\mathrm u}\) for the POD
basis of its state snapshots. The known reduced interface operator is
\[
B_{\mathrm r}^{\mathrm{ODE}}
=\Phi_{\mathrm u}^\top B\Phi_{\mathrm g},
\]
and the local model reduces to
\begin{equation}
\label{eq:dmdc_ode_model}
\dot{\alpha}_{\mathrm u}
=L\alpha_{\mathrm u}
+\sigma B_{\mathrm r}^{\mathrm{ODE}}\alpha_{\mathrm g},
\end{equation}
where \(L\) is obtained from the same state-coordinate regression. 

\paragraph{Linear comparison model}
To isolate the effect of prescribing the interface action, the numerical
experiments also use a standard DMDc baseline that identifies the state and
input operators jointly~\cite{proctor2016dmdc,benner2018iodmd}:
\begin{equation}
\label{eq:dmdc_baseline_regression}
\min_{L_{\mathrm b},B_{\mathrm t,\mathrm b}}
\left\|
\dot{\mathcal S}_{\mathrm u}-
\begin{bmatrix}L_{\mathrm b}&B_{\mathrm t,\mathrm b}\end{bmatrix}
\begin{bmatrix}
\mathcal S_{\mathrm u}\\
\sigma\mathcal S_{\mathrm g}
\end{bmatrix}
\right\|_{\mathrm F}^2.
\end{equation}
In IA-DMDc, by contrast,  \(B_{\mathrm t}\) and \(B_{\mathrm n}\) are determined from the
known interface operator before \(L\) and
\(K_{\mathrm n}\) are identified.  Preserving this available information is an
essential feature of the proposed approach: it enforces the prescribed action
of the interface control and restricts data-driven identification to the
unknown local dynamics. By contrast, learning the state and input operators
simultaneously may confound the internal dynamics with the interface forcing.

\subsection{Online optimization-based coupling}
\label{subsec:ia_dmdc_online}

The offline procedure identifies one local IA-DMDc model for each subproblem.
During the online phase, the identified operators remain fixed, and only the
reduced states and the common interface-control coordinate are determined.
At each time step, the local reduced models are coupled through this common
control, which is selected by minimizing the mismatch between the reconstructed
physical interface quantities. For a DAE subproblem, the pressure coordinate
is recovered algebraically after the state and control have been determined.

We first consider one local subproblem and omit its index for simplicity.
Applying the same one-step \(\theta\)-scheme as in the full-order model to the
tangent evolution equation \eqref{eq:dmdc_tangent_model} gives
\begin{equation}
\label{eq:dmdc_rom_constraint}
\frac{
\alpha_{\mathrm u}^{n+1}-\alpha_{\mathrm u}^{n}
}{\Delta t}
=
L\alpha_{\mathrm u}^{n+\theta}
+\sigma B_{\mathrm t}\alpha_{\mathrm g}^{n+\theta}.
\end{equation}
Once \(\alpha_{\mathrm u}^{n+1}\) and
\(\alpha_{\mathrm g}^{n+1}\) are known, the algebraic normal balance
\eqref{eq:dmdc_normal_model} can be solved for the identifiable pressure
coordinate. Since \(\Sigma_{\mathrm n}\) is nonsingular,
\begin{equation}
\label{eq:dmdc_stage_pressure}
\alpha_{\mathrm p,\mathrm{id}}^{n+\theta}
=
-\Sigma_{\mathrm n}^{-1}
\left(
K_{\mathrm n}\alpha_{\mathrm u}^{n+\theta}
+\sigma B_{\mathrm n}\alpha_{\mathrm g}^{n+\theta}
\right).
\end{equation}
Thus, the pressure is recovered from the state and interface control rather
than advanced by an independent evolution equation or introduced as an
additional optimization variable.
For an ODE subproblem, the algebraic normal balance is absent. Its reduced
state is advanced using \eqref{eq:dmdc_rom_constraint}, with
\(B_{\mathrm t}\) replaced by \(B_{\mathrm r}^{\mathrm{ODE}}\), and no
pressure recovery is required.

For either an ODE or a DAE subproblem,
\(\Phi_{\mathrm u,i}\) denotes the corresponding dynamic-state basis. Since
the physical interface control is approximated by
\(g\approx\Phi_{\mathrm g}\alpha_{\mathrm g}\), its regularization norm is
represented in reduced coordinates by
\[
W_{\mathrm g,\mathrm r}
:=
\Phi_{\mathrm g}^{\top}W_{\mathrm g}\Phi_{\mathrm g}.
\]
At time \(t^{n+1}\), the two reduced states and the common interface control
are determined from
\begin{equation}
\label{eq:dmdc_rom_optimization}
\begin{aligned}
\min_{\substack{
\alpha_{\mathrm u,1}^{n+1},\alpha_{\mathrm u,2}^{n+1},\\
\alpha_{\mathrm g}^{n+1}}}
\quad
J_{\mathrm r}^{n+1}
:={}&
\frac12
\left\|
\widehat{\mathcal R}_{\Gamma,1,h}
\left(
\Phi_{\mathrm u,1}\alpha_{\mathrm u,1}^{n+1},
t^{n+1}
\right)
-
\widehat{\mathcal R}_{\Gamma,2,h}
\left(
\Phi_{\mathrm u,2}\alpha_{\mathrm u,2}^{n+1},
t^{n+1}
\right)
\right\|_{W_\Gamma}^{2}
+
\frac{\lambda_{\mathrm g}}2
\left\|
\alpha_{\mathrm g}^{n+1}
\right\|_{W_{\mathrm g,\mathrm r}}^{2},
\\
\text{subject to}\quad
&\text{\eqref{eq:dmdc_rom_constraint}
for each DAE subproblem},\\
&\text{the corresponding ODE update for each ODE subproblem}.
\end{aligned}
\end{equation}

The online formulation therefore preserves the optimization-based coupling principle introduced in Section~\ref{sec:optimization_framework} while replacing each full-order local solve with a fixed IA-DMDc reduced model. It also provides the interface and DAE structure retained by RQC-DMDc in Section~\ref{sec:regularized_quadratic_dmdc}.

\section{Regularized quadratic correction DMDc}
\label{sec:regularized_quadratic_dmdc}

The linear IA-DMDc model preserves the prescribed interface-input operators and the identified DAE constraint geometry, but a constant state operator can leave nonlinear modal interactions unresolved. 
 We therefore propose regularized quadratic correction DMDc (RQC-DMDc) to enrich the state dependence with quadratic terms. In RQC-DMDc, the linear and quadratic operators are identified jointly, while the prescribed interface operators and DAE geometry remain fixed. The construction can be interpreted as a regularized quadratic operator-inference method \cite{peherstorfer2016operator,qian2022pdes} adapted to the tangent evolution and normal pressure balance of IA-DMDc.  Throughout this section, the local subproblem index is omitted and the notation introduced in Sections~\ref{sec:optimization_framework} and~\ref{sec:interface_aware_dmdc} is retained.

A quadratic model offers greater approximation power, but it also introduces \(\mathcal O(r_{\mathrm u}^3)\) coefficients. Moreover, linear and quadratic features can be strongly correlated when the training trajectory is short or explores a limited region of state space. An unrestricted least-squares fit may then reduce the training residual by assigning excessive weight to the quadratic block, leading to poor recursive prediction. Motivated by block-regularized operator inference \cite{mcquarrie2021regularised,sawant2023physics}, RQC-DMDc mitigates this issue through
root-mean-square feature normalization and block-dependent regularization.

The remainder of this section develops the main components of the RQC-DMDc formulation. We first construct the tangent and normal regression targets by removing the prescribed interface and pressure actions from the snapshot equations. Feature normalization and block-dependent regularization are then introduced to identify the linear and quadratic operators on comparable scales. For comparison, we also introduce an unregularized counterpart to isolate the effect of regularization. Finally, the identified operators are embedded into the IA-DMDc DAE structure and discretized for online coupling. 

\paragraph{Quadratic library and regression targets}
\label{subsec:rqc_fixed_geometry}

The objective is to enrich only the unknown state dependence without
allowing the additional features to redefine the interface action or the
pressure constraint. We therefore reuse the normal geometry
\(
(\widehat C,\widehat G_{\mathrm p},\Phi_{\mathrm n},
\Sigma_{\mathrm n},T_{\mathrm p})
\)
identified from the linear state library as described in Section~3. The quadratic features enrich only the tangent and normal state
actions; the pressure remains an algebraic multiplier and is not included in
the feature library.

For \(\alpha_{\mathrm u}\in\mathbb R^{r_{\mathrm u}}\), define the
\(n_{\mathrm q}=r_{\mathrm u}(r_{\mathrm u}+1)/2\) distinct quadratic
features by
\begin{equation}
\label{eq:qcdmdc_feature_map}
q(\alpha_{\mathrm u})
:=
\left[
\alpha_{\mathrm u,j}\alpha_{\mathrm u,k}
\right]_{1\leq j\leq k\leq r_{\mathrm u}}
\in\mathbb R^{n_{\mathrm q}},
\end{equation}
where any fixed ordering of the monomials may be used. This compact map
removes duplicate off-diagonal products from
\(\alpha_{\mathrm u}\otimes\alpha_{\mathrm u}\) without restricting the
class of homogeneous quadratic vector fields. 

Evaluating these features at the training points \(t^{n+\theta}\) gives the
quadratic snapshot matrix
\begin{equation}
\label{eq:qcdmdc_quadratic_snapshot_matrix}
\mathcal S_{\mathrm q}
:=
\begin{bmatrix}
q(\alpha_{\mathrm u}^{\theta})&
q(\alpha_{\mathrm u}^{1+\theta})&\cdots&
q(\alpha_{\mathrm u}^{N_{\mathrm{tr}}-1+\theta})
\end{bmatrix}
\in\mathbb R^{n_{\mathrm q}\times N_{\mathrm{tr}}}.
\end{equation}
The combined state library is
\begin{equation}
\label{eq:qcdmdc_state_library}
\Xi
:=
\begin{bmatrix}
\mathcal S_{\mathrm u}\\
\mathcal S_{\mathrm q}
\end{bmatrix}
\in
\mathbb R^{(r_{\mathrm u}+n_{\mathrm q})\times N_{\mathrm{tr}}}.
\end{equation}

Because the interface action is known, it is removed before fitting the state
operators. From
Section~\ref{sec:interface_aware_dmdc},
the tangent target is
\begin{equation}
\label{eq:rqc_tangent_target}
\mathcal Y_{\mathrm t}
:=
\dot{\mathcal S}_{\mathrm u}
-\sigma B_{\mathrm t}\mathcal S_{\mathrm g},
\end{equation}
whereas the normal target associated with the fixed pressure geometry is
\begin{equation}
\label{eq:rqc_normal_target}
\mathcal Y_{\mathrm n}
:=
\Phi_{\mathrm n}^{\top}\mathcal Y
-\Sigma_{\mathrm n}\mathcal S_{\mathrm p,\mathrm{id}}.
\end{equation}
The corresponding tangent and normal relations are
\begin{subequations}
\label{eq:rqc_tangent_normal_relations}
\begin{align}
\mathcal Y_{\mathrm t}
&\approx
A_{\mathrm t}\mathcal S_{\mathrm u}
+H_{\mathrm t}\mathcal S_{\mathrm q},
\label{eq:rqc_tangent_relation}\\
\mathcal Y_{\mathrm n}
&\approx
A_{\mathrm n}\mathcal S_{\mathrm u}
+H_{\mathrm n}\mathcal S_{\mathrm q}.
\label{eq:rqc_normal_relation}
\end{align}
\end{subequations}
Thus the prescribed matrices \(B_{\mathrm t}\), \(B_{\mathrm n}\), and
\(\Sigma_{\mathrm n}\), together with the normal geometry, remain unchanged;
only the linear and quadratic state operators are learned. 

\paragraph{Feature normalization and block-selective regularization}

Both regressions
in \eqref{eq:rqc_tangent_normal_relations} have the form
\begin{equation}
\label{eq:rqc_generic_regression}
\mathcal Y_{\ell}
\approx
\Theta_{\ell}\Xi,
\qquad
\Theta_{\ell}
:=
\begin{bmatrix}A_{\ell}&H_{\ell}\end{bmatrix},
\qquad
\ell\in\{\mathrm t,\mathrm n\}.
\end{equation}
They are solved independently and may use different regularization
parameters. Before regularization, the feature rows are normalized so that
the penalty is not determined by their physical scales. For
\(j=1,\ldots,r_{\mathrm u}+n_{\mathrm q}\), define
\begin{equation}
\label{eq:qcdmdc_rms_scale}
s_j
:=
\max\left\{
\left(
\frac{1}{N_{\mathrm{tr}}}
\sum_{m=1}^{N_{\mathrm{tr}}}\Xi_{jm}^2
\right)^{1/2},
s_{\min}
\right\},
\end{equation}
where \(s_{\min}>0\) prevents division by a vanishing feature scale. Let
\begin{equation}
\label{eq:qcdmdc_scaling_matrix}
D_{\Xi}:=\operatorname{diag}
\left(s_1,\ldots,s_{r_{\mathrm u}+n_{\mathrm q}}\right)
\end{equation}
and define
\begin{equation}
\label{eq:qcdmdc_scaled_variables}
\Xi_{\mathrm s}:=D_{\Xi}^{-1}\Xi,
\qquad
\Theta_{\ell,\mathrm s}:=\Theta_{\ell}D_{\Xi}.
\end{equation}
This change of variables leaves the model response invariant because
\(\Theta_{\ell,\mathrm s}\Xi_{\mathrm s}=\Theta_{\ell}\Xi\), while
placing the feature rows on comparable root-mean-square scales. 

To retain the linear component as the leading contribution and treat the
quadratic block as a correction, set
\begin{equation}
\label{eq:qcdmdc_regularization_matrix}
\Omega_{\ell}
:=
\begin{bmatrix}
\sqrt{\omega_{\mathrm A,\ell}}I_{r_{\mathrm u}}&0\\
0&\sqrt{\omega_{\mathrm H,\ell}}I_{n_{\mathrm q}}
\end{bmatrix},
\qquad
0<\omega_{\mathrm A,\ell}\ll\omega_{\mathrm H,\ell}.
\end{equation}
The normalized block ridge problem is
\begin{equation}
\label{eq:qcdmdc_regularized_ls}
\min_{\Theta_{\ell,\mathrm s}}
\quad
\frac{1}{N_{\mathrm{tr}}}
\left\|
\mathcal Y_{\ell}
-\Theta_{\ell,\mathrm s}\Xi_{\mathrm s}
\right\|_{\mathrm F}^2
+
\left\|
\Omega_{\ell}\Theta_{\ell,\mathrm s}^{\top}
\right\|_{\mathrm F}^2.
\end{equation}
The larger value of \(\omega_{\mathrm H,\ell}\) suppresses unsupported
quadratic coefficients while allowing them to represent repeatable nonlinear
defects. The solution is
\begin{equation}
\label{eq:qcdmdc_feature_space_solution}
\Theta_{\ell,\mathrm s}
=
\mathcal Y_{\ell}\Xi_{\mathrm s}^{\top}
\left(
\Xi_{\mathrm s}\Xi_{\mathrm s}^{\top}
+N_{\mathrm{tr}}\Omega_{\ell}^2
\right)^{-1},
\end{equation}
and the operators in the original feature coordinates are recovered as
\begin{equation}
\label{eq:qcdmdc_operator_recovery}
\Theta_{\ell}
=\Theta_{\ell,\mathrm s}D_{\Xi}^{-1}.
\end{equation}

\paragraph{Selection of the regularization weights}

We fix the linear regularization weight at \(\omega_{\mathrm A,\ell}=10^{-8}\) on the normalized features. This weak positive penalty keeps the regularized normal matrix nonsingular when the feature library is rank deficient. For a DAE subproblem, the same quadratic weight is used in the
tangent and normal regressions: 
$
\omega_{\mathrm H,\mathrm t}
=
\omega_{\mathrm H,\mathrm n}
=
\omega_{\mathrm H}.
$
To assess the direct contribution of the quadratic terms, we define
\begin{equation}
\label{eq:relative_quadratic_contribution}
\rho_{Q,\ell}^{\mathrm{train}}
:=
\frac{
\|H_\ell\mathcal S_{\mathrm q}\|_{\mathrm F}
}{
\|A_\ell\mathcal S_{\mathrm u}\|_{\mathrm F}
},
\qquad
\ell\in\{\mathrm t,\mathrm n\},
\end{equation}
The ratio is evaluated using the recovered operators and the unscaled projected snapshots and quadratic features. For a DAE, we monitor
$
\rho_Q^{\mathrm{train}}
:=
\max\left\{
\rho_{Q,\mathrm t}^{\mathrm{train}},
\rho_{Q,\mathrm n}^{\mathrm{train}}
\right\}.
$
For an ODE subproblem, only the tangential differential ratio is used.
We select the quadratic weight empirically so that
$
\rho_Q^{\mathrm{train}}\lesssim10^{-3}.
$
A very small ratio indicates an almost linear fit on the training data,
whereas excessive quadratic action can prevent convergence of the online
nonlinear iterations. The weight is increased to reduce an overly large
ratio and may be decreased when the quadratic correction is overly
suppressed.

\paragraph{IA-OpInf baseline}

To isolate the effect of block-selective regularization, we introduce an IA-OpInf baseline method, which uses the same reduced coordinates, normalized feature library, prescribed
interface operators, and DAE normal geometry, but applies a common weak regularization weight to both operator blocks. It solves
\begin{equation}
\label{eq:ia_opinf_baseline_ls}
\min_{\Theta_{\ell,\mathrm s}^{\mathrm{IA-OpInf}}}
\quad
\frac{1}{N_{\mathrm{tr}}}
\left\|
\mathcal Y_{\ell}
-\Theta_{\ell,\mathrm s}^{\mathrm{IA-OpInf}}\Xi_{\mathrm s}
\right\|_{\mathrm F}^2 
+ 
\omega_{\mathrm A,\ell}
\left\|
\Theta_{\ell,\mathrm s}^{\mathrm{IA-OpInf}}
\right\|_{\mathrm F}^2,
\qquad
\ell\in\{\mathrm t,\mathrm n\}.
\end{equation}
Its minimum-norm coefficients are
\begin{equation}
\label{eq:ia_opinf_baseline_solution}
\Theta_{\ell,\mathrm s}^{\mathrm{IA-OpInf}}
=\mathcal Y_{\ell}\Xi_{\mathrm s}^{\top}
\left(\Xi_{\mathrm s}\Xi_{\mathrm s}^{\top}
+N_{\mathrm{tr}}\,\omega_{\mathrm A,\ell}I\right)^{-1}, 
\qquad
\Theta_{\ell}^{\mathrm{IA-OpInf}}
=\Theta_{\ell,\mathrm s}^{\mathrm{IA-OpInf}}D_{\Xi}^{-1}.
\end{equation}
This baseline therefore differs from RQC-DMDc only through the
block-dependent penalties. It is distinct from the standard DMDc baseline of
Section~3, which additionally learns the interface-input matrix. As
\(\omega_{\mathrm H,\ell}\to\infty\) with
\(\omega_{\mathrm A,\ell}\) fixed, the quadratic block vanishes and
RQC-DMDc approaches a ridge-regularized linear fit.

\paragraph{Reduced model and online solution}
\label{subsec:rqc_reduced_dae}
Restoring the
prescribed interface terms, the identified tangent evolution and normal
balance are
\begin{subequations}
\label{eq:qcdmdc_surrogate}
\begin{align}
\dot{\alpha}_{\mathrm u}
&=
A_{\mathrm t}\alpha_{\mathrm u}
+H_{\mathrm t}q(\alpha_{\mathrm u})
+\sigma B_{\mathrm t}\alpha_{\mathrm g},
\label{eq:rqc_tangent_model}\\
0
&=
A_{\mathrm n}\alpha_{\mathrm u}
+H_{\mathrm n}q(\alpha_{\mathrm u})
+\Sigma_{\mathrm n}\alpha_{\mathrm p,\mathrm{id}}
+\sigma B_{\mathrm n}\alpha_{\mathrm g}.
\label{eq:rqc_normal_model}
\end{align}
\end{subequations}
Using the auxiliary normal coordinate \(\alpha_{\mathrm n}\), the
corresponding descriptor representation is
\begin{equation}
\label{eq:rqc_reduced_dae_descriptor}
E_{\mathrm r}\frac{\mathrm d}{\mathrm dt}
\begin{bmatrix}
\alpha_{\mathrm u}\\
\alpha_{\mathrm n}\\
\alpha_{\mathrm p,\mathrm{id}}
\end{bmatrix}
=
A_{\mathrm r}^{\mathrm{RQC}}
\begin{bmatrix}
\alpha_{\mathrm u}\\
\alpha_{\mathrm n}\\
\alpha_{\mathrm p,\mathrm{id}}
\end{bmatrix}
+
H_{\mathrm r}^{\mathrm{RQC}}q(\alpha_{\mathrm u})
+\sigma B_{\mathrm r}^{\mathrm{DAE}}\alpha_{\mathrm g},
\end{equation}
where \(E_{\mathrm r}\) and \(B_{\mathrm r}^{\mathrm{DAE}}\) are the
same as in the linear model, and
\begin{equation}
\label{eq:rqc_reduced_dae_matrices}
A_{\mathrm r}^{\mathrm{RQC}}
=
\begin{bmatrix}
A_{\mathrm t}&0&0\\
A_{\mathrm n}&0&\Sigma_{\mathrm n}\\
0&\Sigma_{\mathrm n}^{\top}&0
\end{bmatrix},
\qquad
H_{\mathrm r}^{\mathrm{RQC}}
=
\begin{bmatrix}
H_{\mathrm t}\\
H_{\mathrm n}\\
0
\end{bmatrix}.
\end{equation}
The three block rows of \eqref{eq:rqc_reduced_dae_descriptor} are the tangent
evolution, an auxiliary normal evolution, and the constraint
\(0=\Sigma_{\mathrm n}^{\top}\alpha_{\mathrm n}\), respectively. Since
\(\Sigma_{\mathrm n}\) is nonsingular, a consistent initial condition
\(\alpha_{\mathrm n}^0=0\) gives
\(\alpha_{\mathrm n}=\dot\alpha_{\mathrm n}=0\). The second block row then
reduces exactly to the algebraic normal balance in
\eqref{eq:rqc_normal_model}. Thus the quadratic terms change the reduced
vector field without changing the differential--algebraic structure.

Solving the normal balance for the identifiable pressure coordinate gives
\begin{equation}
\label{eq:rqc_pressure_recovery}
\alpha_{\mathrm p,\mathrm{id}}
(\alpha_{\mathrm u},\alpha_{\mathrm g})
=
-\Sigma_{\mathrm n}^{-1}
\left(
A_{\mathrm n}\alpha_{\mathrm u}
+H_{\mathrm n}q(\alpha_{\mathrm u})
+\sigma B_{\mathrm n}\alpha_{\mathrm g}
\right).
\end{equation}
For an ODE subproblem, the normal equation is omitted, and the model reduces
to
\[
\dot{\alpha}_{\mathrm u}
=A_{\mathrm t}\alpha_{\mathrm u}
+H_{\mathrm t}q(\alpha_{\mathrm u})
+\sigma B_{\mathrm r}^{\mathrm{ODE}}\alpha_{\mathrm g}.
\]

Applying the same one-step \(\theta\)-scheme as in the full-order model to
the tangent equation gives
\begin{equation}
\label{eq:qcdmdc_rom_constraint}
\frac{
\alpha_{\mathrm u}^{n+1}-\alpha_{\mathrm u}^{n}
}{\Delta t}
=
A_{\mathrm t}\alpha_{\mathrm u}^{n+\theta}
+H_{\mathrm t}q(\alpha_{\mathrm u}^{n+\theta})
+\sigma B_{\mathrm t}\alpha_{\mathrm g}^{n+\theta}.
\end{equation}
The pressure coordinate at the same evaluation point is
\begin{equation}
\label{eq:rqc_stage_pressure}
\alpha_{\mathrm p,\mathrm{id}}^{n+\theta}
=
-\Sigma_{\mathrm n}^{-1}
\left(
A_{\mathrm n}\alpha_{\mathrm u}^{n+\theta}
+H_{\mathrm n}q(\alpha_{\mathrm u}^{n+\theta})
+\sigma B_{\mathrm n}\alpha_{\mathrm g}^{n+\theta}
\right).
\end{equation}
We note that direct discretization of the descriptor system is algebraically equivalent
to advancing the tangent equation and then recovering the pressure from
\eqref{eq:rqc_stage_pressure}. Here, the eliminated form is used within the online
interface optimization.

Unlike linear IA-DMDc, the new time-step equation is nonlinear because
\(q(\alpha_{\mathrm u}^{n+\theta})\) depends quadratically on the unknown
state \(\alpha_{\mathrm u}^{n+1}\). For a fixed \(\alpha_{\mathrm g}^{n+1}\), define the nonlinear tangent
residual
\begin{equation}
\label{eq:qcdmdc_nonlinear_residual}
\begin{aligned}
\mathcal R^{n+1}
(\alpha_{\mathrm u}^{n+1},\alpha_{\mathrm g}^{n+1})
:={}&
\alpha_{\mathrm u}^{n+1}-\alpha_{\mathrm u}^{n}
-\Delta t
\left[
A_{\mathrm t}\alpha_{\mathrm u}^{n+\theta}
+H_{\mathrm t}q(\alpha_{\mathrm u}^{n+\theta})
+\sigma B_{\mathrm t}\alpha_{\mathrm g}^{n+\theta}
\right].
\end{aligned}
\end{equation}
If
\(Dq(\alpha_{\mathrm u})
\in\mathbb R^{n_{\mathrm q}\times r_{\mathrm u}}\)
is the Jacobian of the compact feature map, then
\begin{equation}
\label{eq:qcdmdc_state_jacobian}
\frac{\partial\mathcal R^{n+1}}
{\partial\alpha_{\mathrm u}^{n+1}}
=
I_{r_{\mathrm u}}
-\theta\Delta t
\left[
A_{\mathrm t}
+H_{\mathrm t}Dq(\alpha_{\mathrm u}^{n+\theta})
\right].
\end{equation}
The small nonlinear system can be solved by Newton iterations
inside the interface optimization.

\paragraph{Interpretation of quadratic regularization}
\label{subsec:qcdmdc_stability_scope}
%The stronger penalty on the quadratic block is intended to reduce overfitting
%and excessive nonlinear sensitivity, but it is not a nonlinear stability
%constraint. 
%
%Similarly, choosing \(\theta\in[1/2,1]\) provides \(A\)-stability for the
%linear test equation but does not guarantee global solvability, boundedness,
%or unconditional stability of the quadratic model. Even if
%\(A_{\mathrm t}\) is Hurwitz, stability of the complete quadratic vector
%field requires additional dissipative or Lyapunov structure
%\cite{sawant2023physics,goyal2025stablequadratic}.
%
%RQC-DMDc therefore has a specific, accuracy-oriented role. It enriches the
%tangent dynamics and normal pressure balance while retaining the prescribed
%interface-input operators and constraint geometry. Relative to linear
%IA-DMDc, its online cost is limited to evaluating the compact quadratic
%features and solving a small nonlinear tangent equation. In the
%defect--amplification terminology of
%Subsection~\ref{subsec:defect_amplification}, Hurwitz-DMDc targets error
%amplification, whereas RQC-DMDc targets systematic nonlinear defects. Neither
%correction is expected to dominate in every regime.
RQC-DMDc is used here as an accuracy-oriented enrichment. It adds regularized quadratic state dependence to the tangent dynamics and, for a DAE subproblem, to the normal pressure balance, while retaining the prescribed interface-input operators and the identified DAE geometry. 
The pressure remains algebraic and is recovered at each time step through \eqref{eq:rqc_stage_pressure}, and the additional online work relative to linear IA-DMDc consists of evaluating the quadratic terms and solving the nonlinear reduced tangent equation \eqref{eq:qcdmdc_nonlinear_residual}.

\section{Numerical experiments}
\label{sec:numerical_experiments}

We assess IA-DMDc and RQC-DMDc on four model problems. The heterogeneous Fisher--KPP problem couples two reaction--diffusion ODE subproblems. The Astorino--Grandmont, freely oscillating-cylinder, and Turek problems couple an incompressible-flow DAE with a structural ODE. The Fisher--KPP interval is short relative to the reaction time scales, and convection vanishes along the manufactured Astorino--Grandmont trajectory. The two wake-driven FSI problems assess prediction beyond training under sustained vortex shedding and structural motion. All models use full-order trajectories generated by the optimization-based coupling formulation. IA-DMDc and RQC-DMDc use the same snapshots, reduced spaces, time discretization, and projected interface-input operators. The DMDc baseline in \eqref{eq:dmdc_baseline_regression} assesses joint state--input identification. The IA-OpInf baseline in \eqref{eq:ia_opinf_baseline_ls} uses the same quadratic library and interface-aware structure as RQC-DMDc, with the common weak penalty on both operator blocks. For two coupled subsystems, \(\mathcal M_1\mid\mathcal M_2\) denotes the local models on the two sides of the interface.

At time \(t^k\), the relative error of a reported field \(Y_i\) is measured in the Frobenius norm as
\begin{equation}
\label{eq:re}
\operatorname{RE}_{Y_i}(t^k) =
\frac{
\left\| Y_{i,\mathrm{ROM}}^k - Y_{i,\mathrm{ref}}^k \right\|_{\mathrm F}
}{
\left\| Y_{i,\mathrm{ref}}^k \right\|_{\mathrm F}
},
\end{equation}
where \(Y_{i,\mathrm{ROM}}^k\) is reconstructed in the full finite-element coordinate space before the norm is evaluated. Kinetic energies are evaluated from the physical velocity fields as
\begin{equation}
\label{eq:kinetic_energy_definition}
E_{\mathrm f}^{\mathrm{kin}}(t)
=\frac{\rho_{\mathrm f}}{2}\int_{\Omega_{\mathrm f}}|u_{\mathrm f}(x,t)|^2\,\mathrm{d}x,
\qquad
E_{\mathrm s}^{\mathrm{kin}}(t)
=\frac{\rho_{\mathrm s}}{2}\int_{\Omega_{\mathrm s}}|\partial_t d_{\mathrm s}(x,t)|^2\,\mathrm{d}x.
\end{equation}

\subsection{Heterogeneous Fisher--KPP problem}
\label{subsec:fisher_kpp}

We consider a heterogeneous Fisher--KPP problem~\cite{Fisher1937}
on the unit disk $
\Omega=\{\boldsymbol{x}\in\mathbb{R}^2:
\|\boldsymbol{x}\|_2<1\}.
$
The circular geometry and zero-flux exterior boundary are motivated
by the diffusion test in~\cite{HerreraHernandez2017}.
The domain is partitioned into
\(\Omega_1=\Omega\cap\{x_1<0\}\) and
\(\Omega_2=\Omega\cap\{x_1>0\}\), with interface
\(\Gamma=\{0\}\times(-1,1)\).
In each subdomain, the normalized density satisfies
\begin{equation*}
\frac{\partial u_i}{\partial t}
=
D_i\Delta u_i+\rho_i u_i(1-u_i),
\qquad
(\boldsymbol{x},t)\in\Omega_i\times(0,T],
\quad i=1,2.
\end{equation*}
We choose
$
(D_1,\rho_1)=(0.60,1.00),
(D_2,\rho_2)=(0.25,1.40),
$
so that the diffusion coefficient is larger in \(\Omega_1\) while the reaction rate is larger in \(\Omega_2\).
Zero diffusive flux is prescribed on the exterior boundary.
Across \(\Gamma\), we impose continuity of the density and
balance of the normal diffusive flux:
$
u_1=u_2,
g=D_1\nabla u_1\cdot n_1
=-D_2\nabla u_2\cdot n_2,
$
where \(n_i\) denotes the outward unit normal to \(\Omega_i\).
With
\(r(\boldsymbol{x})=\|\boldsymbol{x}-(0.08,0)\|_2\),
the initial condition is
\begin{equation*}
u_i(\boldsymbol{x},0)
=
\begin{cases}
0, & r(\boldsymbol{x})<R_0,\\
\frac12, & r(\boldsymbol{x})=R_0,\\
1, & r(\boldsymbol{x})>R_0,
\end{cases}
\end{equation*}
with $R_0=0.30.$ The initially empty region remains part of the computational
domain and introduces no internal boundary.

The discretization contains 3121 density degrees of freedom in each subdomain and 97 interface-control degrees of freedom. The full-order problem is advanced by backward Euler with
\(\Delta t=3\times10^{-5}\) to \(T=0.03\). Snapshots are saved at every time step.
The first 250 snapshot columns, ending at \(t=7.5\times10^{-3}\), are used for training, and the reduced models predict the remaining evolution. The reduced dimensions of both density fields and the interface control are all 10.

We compare IA-DMDc and RQC-DMDc on both sides of the interface,
using the same snapshots, reduced spaces, time discretization,
and prescribed interface-input operators.
The DMDc and IA-OpInf baselines are included for comparison.
For RQC-DMDc, we use \(\omega_{\mathrm H}=10^{-3}\) in both
subdomains. The corresponding quadratic-to-linear training
action ratios are \(\rho_Q^{\mathrm{train}}=6.786\times10^{-4}\) in \(\Omega_1\) and
\(\rho_Q^{\mathrm{train}}=2.436\times10^{-4}\) in \(\Omega_2\).
Although the reaction term is quadratic, the simulation
interval is short relative to the reaction time scales,
with \(\max_i\rho_iT=0.042\).
This example therefore examines whether the regularized
quadratic enrichment retains the accuracy of the linear
interface-aware model over a short prediction horizon.

\begin{figure}[H]
    \centering
    \captionsetup[subfloat]{font=footnotesize,captionskip=2pt}

    {\subfloat[{$t=0.0006$}]{%
        \begin{minipage}[t]{0.31\textwidth}
        \centering
        \includegraphics[width=\linewidth]{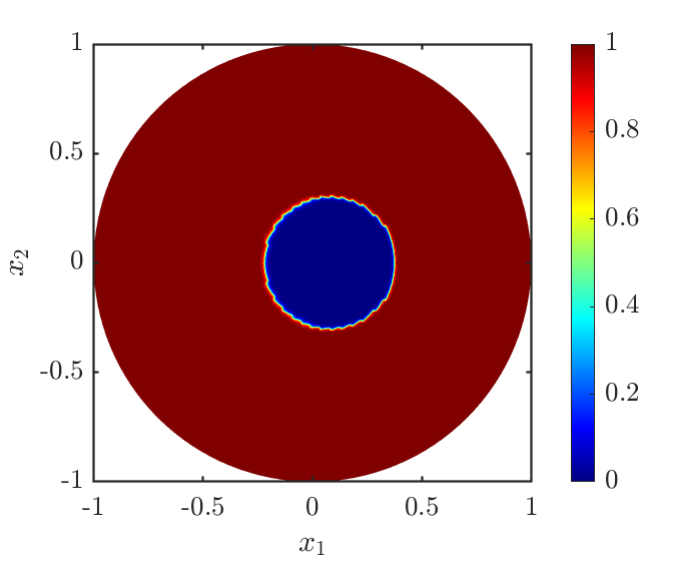}
        \end{minipage}%
    }}%
    \hspace{0.02\textwidth}%
    {\subfloat[{$t=0.01002$}]{%
        \begin{minipage}[t]{0.31\textwidth}
        \centering
        \includegraphics[width=\linewidth]{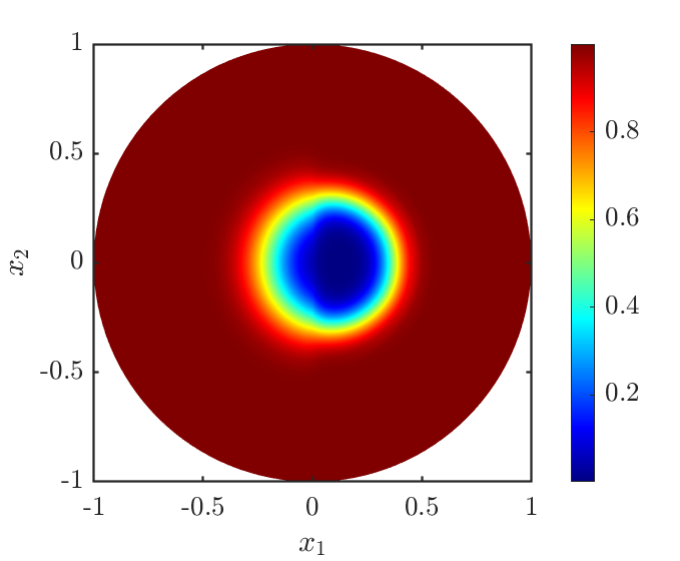}
        \end{minipage}%
    }}%
    \hspace{0.02\textwidth}%
    {\subfloat[{$t=0.03$}]{%
        \begin{minipage}[t]{0.31\textwidth}
        \centering
        \includegraphics[width=\linewidth]{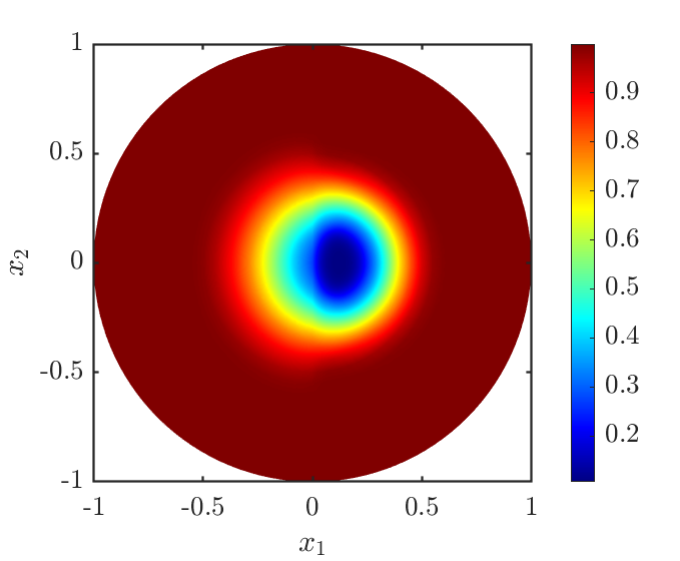}
        \end{minipage}%
    }}

    \vspace{0.25em}

    {\subfloat[{$t=0.0006$}]{%
        \begin{minipage}[t]{0.31\textwidth}
        \centering
        \includegraphics[width=\linewidth]{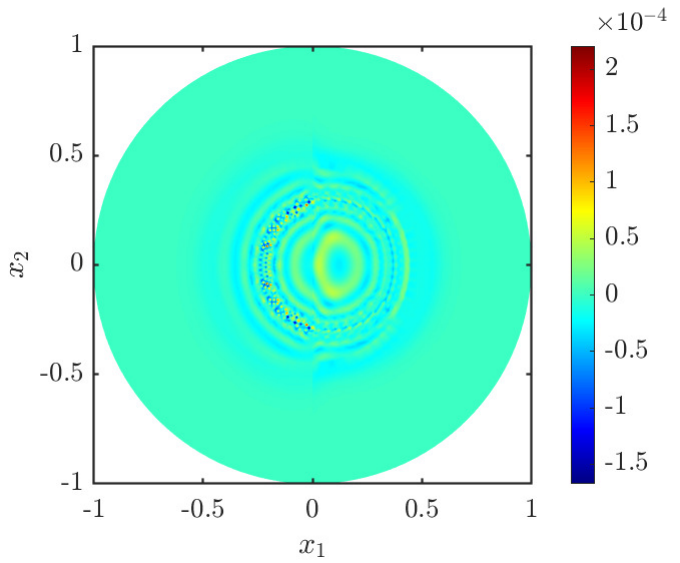}
        \end{minipage}%
    }}%
    \hspace{0.02\textwidth}%
    {\subfloat[{$t=0.01002$}]{%
        \begin{minipage}[t]{0.31\textwidth}
        \centering
        \includegraphics[width=\linewidth]{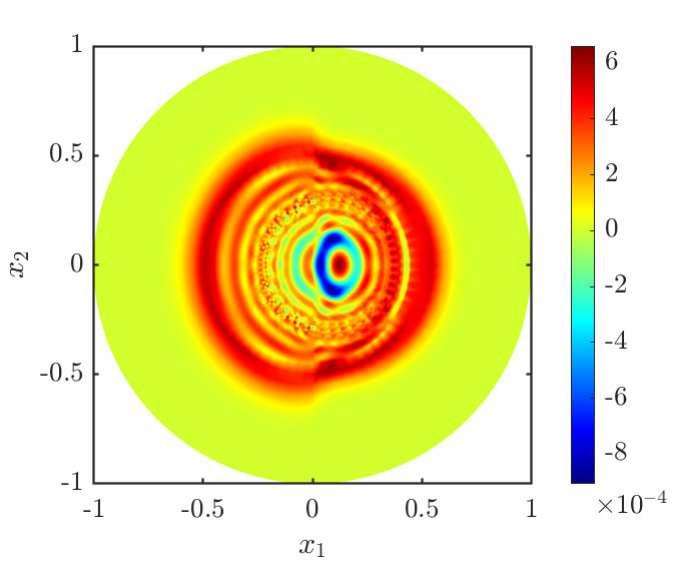}
        \end{minipage}%
    }}%
    \hspace{0.02\textwidth}%
    {\subfloat[{$t=0.03$}]{%
        \begin{minipage}[t]{0.31\textwidth}
        \centering
        \includegraphics[width=\linewidth]{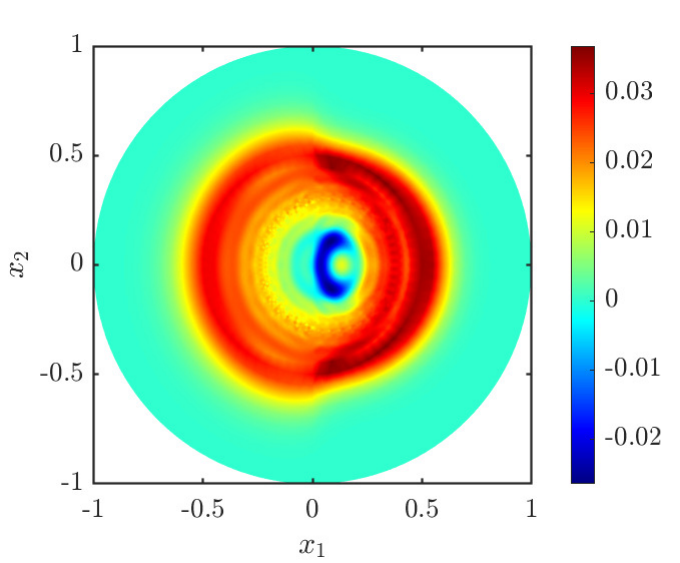}
        \end{minipage}%
    }}

    \caption{Heterogeneous Fisher--KPP problem; IA-DMDc prediction of the density
    (top) and pointwise error relative to the full-order solution (bottom).}
    \label{fig:fisher_dmdc_fields}
\end{figure}

\begin{figure}[H]
    \centering
    \captionsetup[subfloat]{font=footnotesize,captionskip=2pt}

    {\subfloat[{$t=0.0006$}]{%
        \begin{minipage}[t]{0.31\textwidth}
        \centering
        \includegraphics[width=\linewidth]{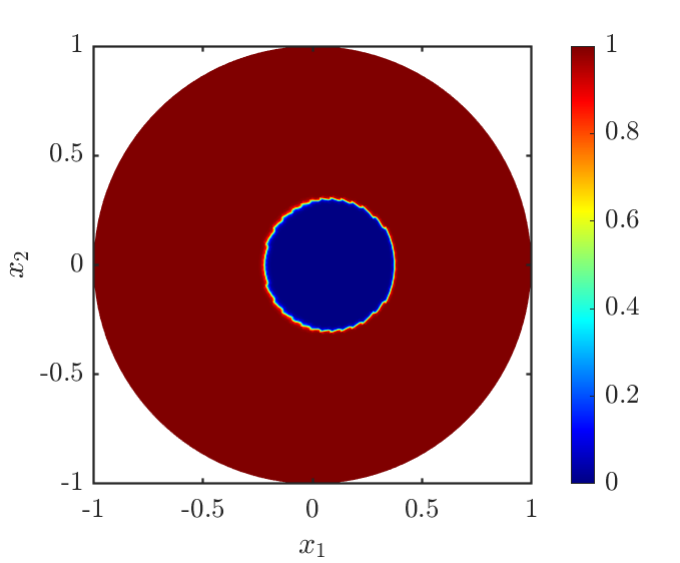}
        \end{minipage}%
    }}%
    \hspace{0.02\textwidth}%
    {\subfloat[{$t=0.01002$}]{%
        \begin{minipage}[t]{0.31\textwidth}
        \centering
        \includegraphics[width=\linewidth]{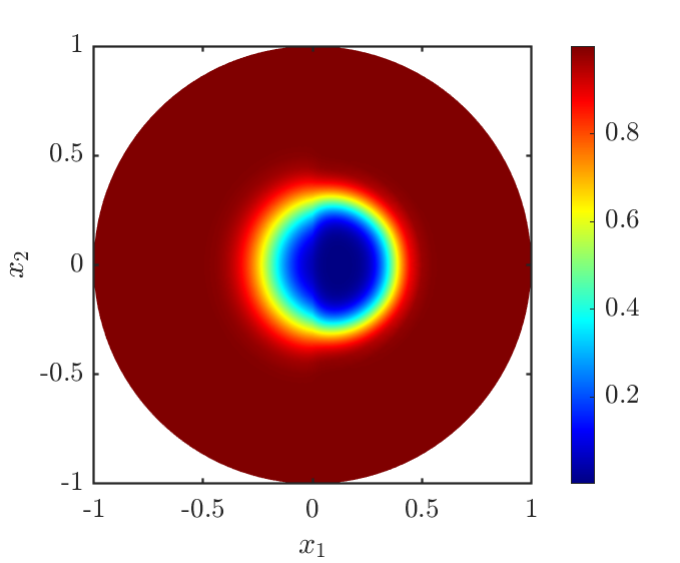}
        \end{minipage}%
    }}%
    \hspace{0.02\textwidth}%
    {\subfloat[{$t=0.03$}]{%
        \begin{minipage}[t]{0.31\textwidth}
        \centering
        \includegraphics[width=\linewidth]{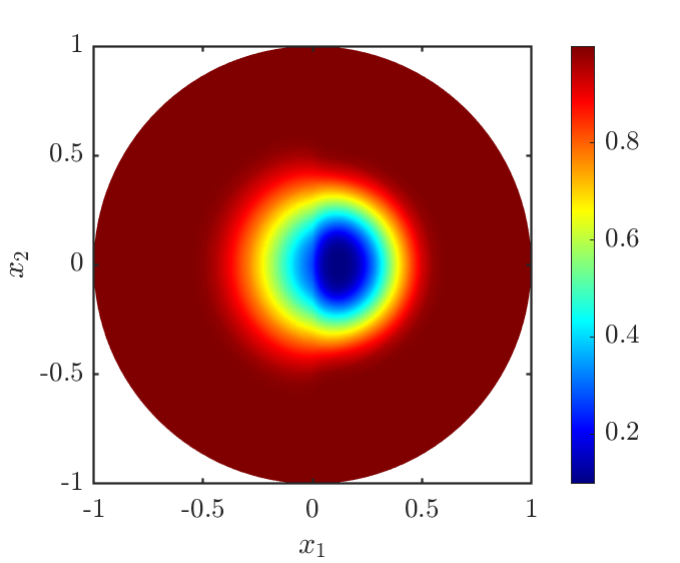}
        \end{minipage}%
    }}

    \vspace{0.25em}

    {\subfloat[{$t=0.0006$}]{%
        \begin{minipage}[t]{0.31\textwidth}
        \centering
        \includegraphics[width=\linewidth]{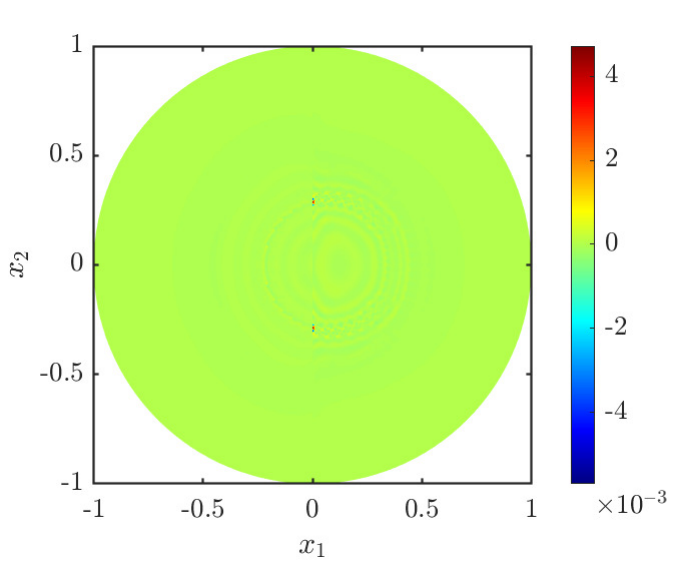}
        \end{minipage}%
    }}%
    \hspace{0.02\textwidth}%
    {\subfloat[{$t=0.01002$}]{%
        \begin{minipage}[t]{0.31\textwidth}
        \centering
        \includegraphics[width=\linewidth]{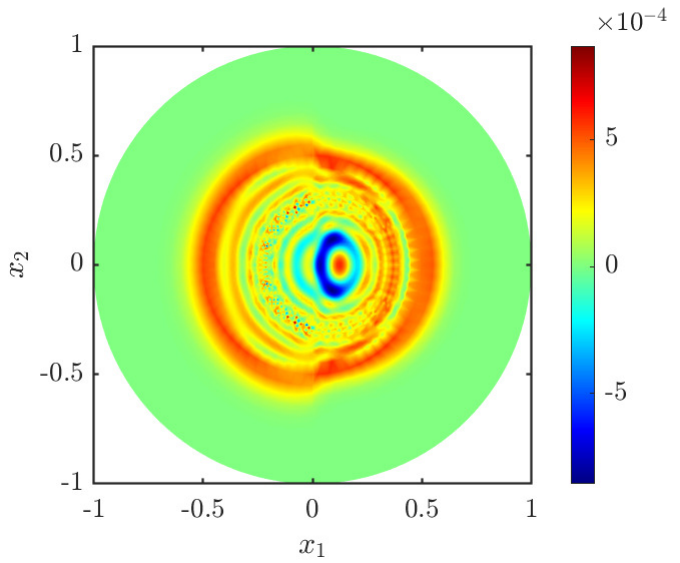}
        \end{minipage}%
    }}%
    \hspace{0.02\textwidth}%
    {\subfloat[{$t=0.03$}]{%
        \begin{minipage}[t]{0.31\textwidth}
        \centering
        \includegraphics[width=\linewidth]{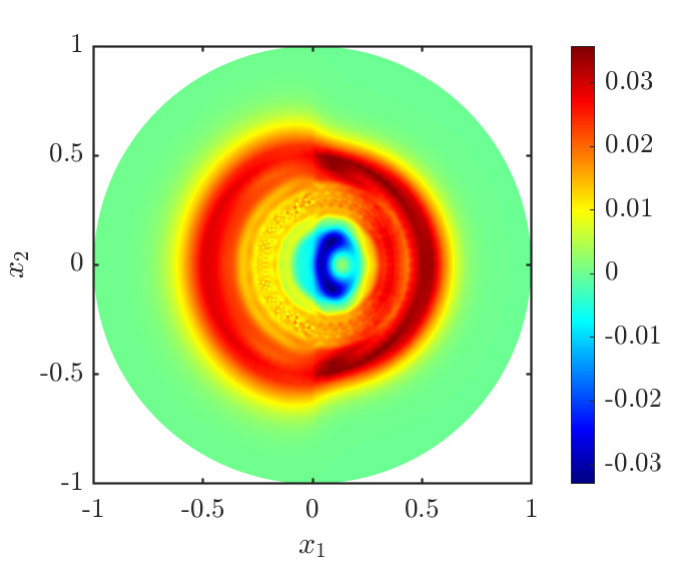}
        \end{minipage}%
    }}

    \caption{Heterogeneous Fisher--KPP problem; RQC-DMDc prediction of the density
    (top) and pointwise error relative to the full-order solution (bottom).}
    \label{fig:fisher_rqc_fields}
\end{figure}

Figures~\ref{fig:fisher_dmdc_fields} and~\ref{fig:fisher_rqc_fields} present the density predictions and pointwise errors. Both models reproduce the smoothing of the initially sharp density transition and the gradual filling of the initially empty region.
At \(t=0.01002\), the discrepancies are concentrated around the density transition and near the center of the disk. At \(t=0.03\), both models retain the main density profile, with localized pointwise errors of order \(10^{-2}\). The field comparisons show similar predictive accuracy for IA-DMDc and RQC-DMDc.

\begin{figure}[H]
    \centering
    \captionsetup[subfloat]{font=footnotesize,captionskip=2pt}

    {\subfloat[{Density \(u_1\).}]{%
        \begin{minipage}[t]{0.43\textwidth}
        \centering
        \includegraphics[width=\linewidth]{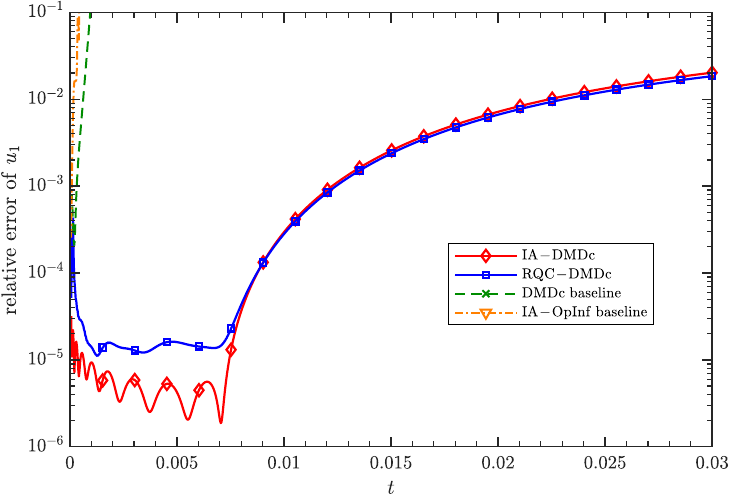}
        \end{minipage}%
    }}%
    \hspace{0.03\textwidth}%
    {\subfloat[{Density \(u_2\).}]{%
        \begin{minipage}[t]{0.43\textwidth}
        \centering
        \includegraphics[width=\linewidth]{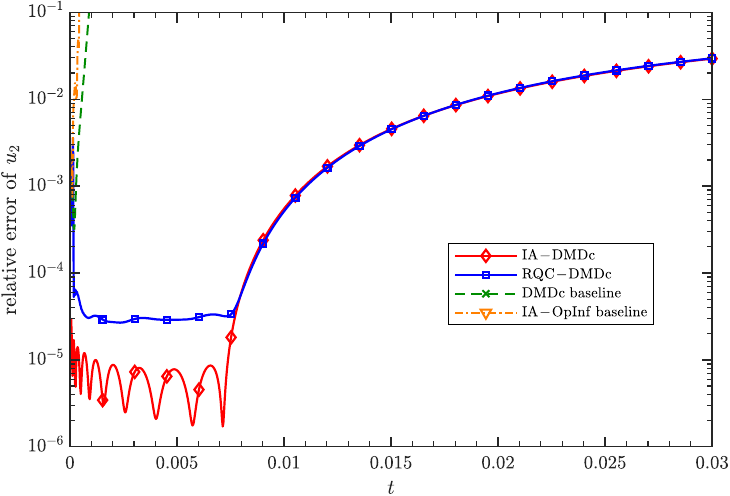}
        \end{minipage}%
    }}

    \caption{Heterogeneous Fisher--KPP problem; relative errors of the densities
    in the two subdomains. The training interval ends at 
    \(t=7.5\times10^{-3}\).}
    \label{fig:fisher_relative_errors}
\end{figure}

\begin{figure}[H]
    \centering
    \captionsetup[subfloat]{font=footnotesize,captionskip=2pt}

    {\subfloat[{Subdomain \(\Omega_1\).}\label{fig:fisher_free_energy_1}]{%
        \begin{minipage}[t]{0.43\textwidth}
        \centering
        \includegraphics[width=\linewidth]{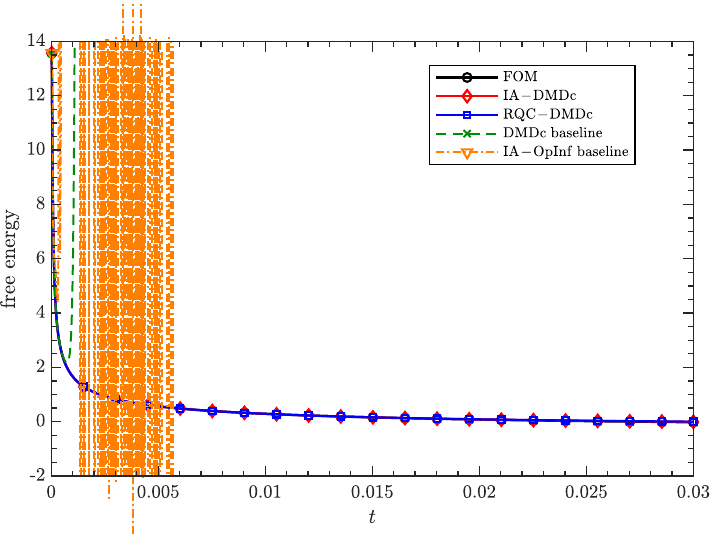}
        \end{minipage}%
    }}%
    \hspace{0.03\textwidth}%
    {\subfloat[{Subdomain \(\Omega_2\).}\label{fig:fisher_free_energy_2}]{%
        \begin{minipage}[t]{0.43\textwidth}
        \centering
        \includegraphics[width=\linewidth]{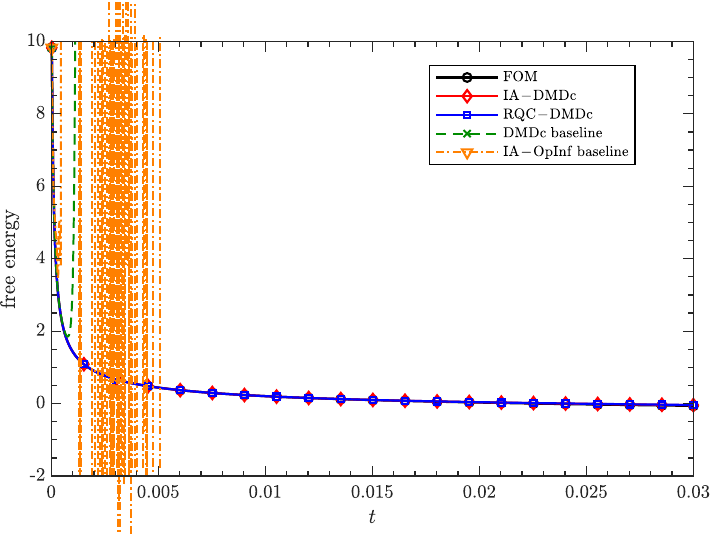}
        \end{minipage}%
    }}

    \caption{Heterogeneous Fisher--KPP problem; free energies in the two
    subdomains. The training interval ends near \(t=7.5\times10^{-3}\).}
    \label{fig:fisher_free_energies}
\end{figure}

Figure~\ref{fig:fisher_relative_errors} distinguishes the training accuracy from the subsequent predictive accuracy. IA-DMDc gives smaller relative errors during most of training. Beyond the training interval, the error histories become nearly coincident and increase gradually to order \(10^{-2}\) at the final time. RQC-DMDc gives a small reduction in the late-time error of \(u_1\), while the two predictions of \(u_2\) have essentially the same relative error. The quadratic enrichment therefore provides a limited accuracy gain in this example. To assess the integrated behavior of the predicted fields, we evaluate the subdomain free energies
$
E_i(t)
=
\int_{\Omega_i}
\left[
\frac{D_i}{2}|\nabla u_i|^2
+
\rho_i\left(\frac{u_i^3}{3}-\frac{u_i^2}{2}\right)
\right]\,\mathrm{d}{x},
i=1,2. 
$
The full-order histories in Figure~\ref{fig:fisher_free_energies} exhibit a rapid initial decrease followed by slower relaxation. Both IA-DMDc and RQC-DMDc closely reproduce these two stages in each subdomain, and their energy curves remain nearly coincident beyond the training interval. 
The DMDc baseline develops substantial errors early in the calculation, accompanied by an increase in the reported free energies. The IA-OpInf baseline also departs from the reference trajectory and exhibits large energy oscillations during the initial transient.

\subsection{Astorino--Grandmont problem}
\label{subsec:astorino_grandmont}

We consider a Navier--Stokes variant of the Astorino--Grandmont problem
introduced in~\cite{astorino2010convergence}. The fluid and solid domains are \(\Omega_{\mathrm f}=(0,1)\times(0,1)\), \(\Omega_{\mathrm s}=(0,1)\times(1,1.25)\), respectively, and share the interface \(\Sigma=(0,1)\times\{1\}\).
The structure satisfies the linear elastodynamics model
\begin{equation}
\label{eq:linear_elastic_solid}
\begin{aligned}
\rho_{\mathrm s}\partial_{tt}d_{\mathrm s}
-\nabla\cdot\sigma_{\mathrm s}(d_{\mathrm s})
&=f_{\mathrm s}
\qquad \text{in } \Omega_{\mathrm s},\\
\sigma_{\mathrm s}(d_{\mathrm s})
&=2\mu_{\mathrm s}\varepsilon_{\mathrm s}
+\lambda_{\mathrm s}(\nabla\cdot d_{\mathrm s})I,
\end{aligned}
\end{equation}
where \(\varepsilon_{\mathrm s}
=\tfrac12(\nabla d_{\mathrm s}+\nabla d_{\mathrm s}^{\top})\) is the
 strain tensor and \(f_{\mathrm s}\) is the body force. 
As in the original benchmark, the structural deformation is neglected in the geometric description, so \(\Omega_{\mathrm f}\) remains
fixed, although the kinematic and dynamic conditions are still imposed on
\(\Sigma\). Neumann conditions are prescribed on the left and right boundaries
of both subdomains; whereas the upper solid boundary and the lower fluid boundary are
subject to Dirichlet conditions. The material parameters are
\[
\rho_{\mathrm f}=1,\qquad \mu_{\mathrm f}=0.013,\qquad
\rho_{\mathrm s}=1.9,\qquad \mu_{\mathrm s}=3,\qquad \lambda_{\mathrm s}=4.5.
\]

We choose the initial data, boundary data, and source terms such that the
coupled Navier--Stokes--elasticity problem admits the manufactured solution
\begin{equation*}
\begin{aligned}
u_{\mathrm f}(x,y,t)
&=
\begin{bmatrix}
\cos(x+t)\sin(y+t)+\sin(x+t)\cos(y+t)\\
-\sin(x+t)\cos(y+t)-\cos(x+t)\sin(y+t)
\end{bmatrix},
\\[0.5em]
d_{\mathrm s}(x,y,t)
&=
\begin{bmatrix}
\sin(x+t)\sin(y+t)\\
\cos(x+t)\cos(y+t)
\end{bmatrix},
\\[0.5em]
p_{\mathrm f}(x,y,t)
&=2\mu_{\mathrm f}\left[
\sin(x+t)\sin(y+t)-\cos(x+t)\cos(y+t)
\right]
+2\mu_{\mathrm s}\cos(x+t)\sin(y+t).
\end{aligned}
\end{equation*}
The velocity is divergence free and satisfies
\((u_{\mathrm f}\cdot\nabla)u_{\mathrm f}=0\); hence, the convective term vanishes along this
trajectory even though it is retained in the full-order equations. Moreover,
\(\partial_t d_{\mathrm s}=u_{\mathrm f}\) satisfies the kinematic transmission condition, while the prescribed pressure and source terms enforce the dynamic condition.
The discretization contains \(N_{\mathrm f}=N_{\mathrm u}+N_{\mathrm p}=8403\) fluid degrees of freedom, with
\(N_{\mathrm u}=7442\) and \(N_{\mathrm p}=961\), together with \(N_{\mathrm s}=5002\) solid and
\(N_{\mathrm g}=122\) interface-control degrees of freedom. Snapshots are sampled every
\(\Delta t_{\mathrm s}=10^{-3}\,\mathrm{}\). The first 2500 snapshots,
\(t\in(0,2.5]\,\mathrm{s}\), are used for training, and the reduced models are
advanced to \(t=10\,\mathrm{s}\). The reduced ranks are
\(r_{\mathrm u}=5\), \(r_{\mathrm p}=4\), \(r_{\mathrm d}=3\), and \(r_{\mathrm g}=2\) for the fluid velocity,
fluid pressure, solid displacement, and interface control, respectively. 
The solid reduced model is identified by IA-DMDc; whereas for the fluid, we compare IA-DMDc, and RQC-DMDc using the same snapshots and reduced spaces. 
For RQC-DMDc, the quadratic weight is set \(\omega_{\mathrm H}=1\), which correspond to \(\rho_Q^{\mathrm{train}}=3.699\times10^{-4}\).
Since the convective term vanishes on the manufactured trajectory, this case tests whether the stability and quadratic corrections disrupt an already accurate linear reduced dynamics.

For this solution,
\(u_{\mathrm f}=\partial_t d_{\mathrm s}=[\sin(x+y+2t),-\sin(x+y+2t)]^{\top}\). The kinetic-energy
definition~\eqref{eq:kinetic_energy_definition} therefore gives
\begin{equation}
\label{eq:astorino_exact_kinetic_energies}
\begin{aligned}
E_{\mathrm f}^{\mathrm{kin}}(t)
&=\rho_{\mathrm f}\left[
\frac{1}{2}-\frac{1}{2}\sin^2(1)\cos(4t+2)
\right],
\\
E_{\mathrm s}^{\mathrm{kin}}(t)
&=\rho_{\mathrm s}\left[
\frac{1}{8}-\frac{1}{2}\sin(1)\sin\!\left(\frac{1}{4}\right)
\cos\!\left(4t+\frac{13}{4}\right)
\right].
\end{aligned}
\end{equation}

% ============================================================
% IA-DMDc--(H-)DMDc predictions and errors
% ============================================================

\begin{figure}[H]
    \centering
    \captionsetup[subfloat]{font=footnotesize,captionskip=1pt}

    % --------------------------------------------------------
    % Fluid velocity
    % --------------------------------------------------------
    {\subfloat[{$t=0.2\,\mathrm{s}$}\label{fig:astorino_dh_w_t2}]{%
        \begin{minipage}[t]{0.34\textwidth}
        \centering
        \includegraphics[width=\linewidth]
        {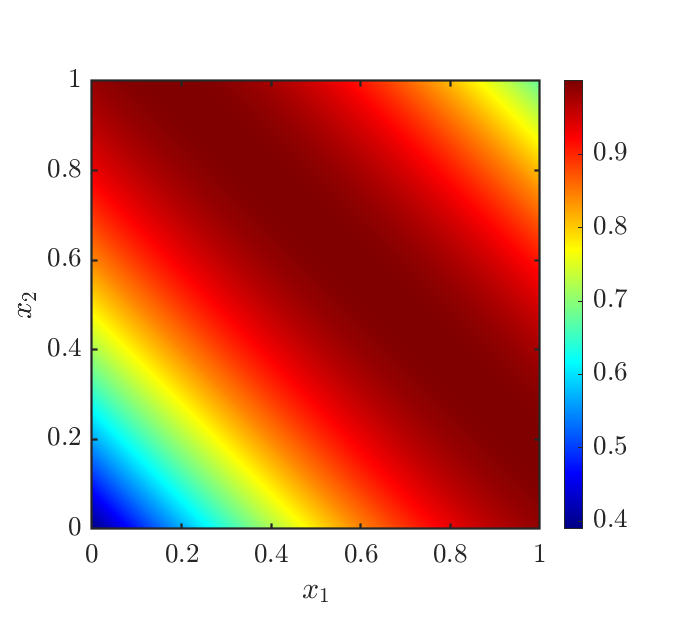}
        \end{minipage}%
    }}
    \hspace{0.04\textwidth}
    {\subfloat[{$t=10\,\mathrm{s}$}\label{fig:astorino_dh_w_t10}]{%
        \begin{minipage}[t]{0.34\textwidth}
        \centering
        \includegraphics[width=\linewidth]
        {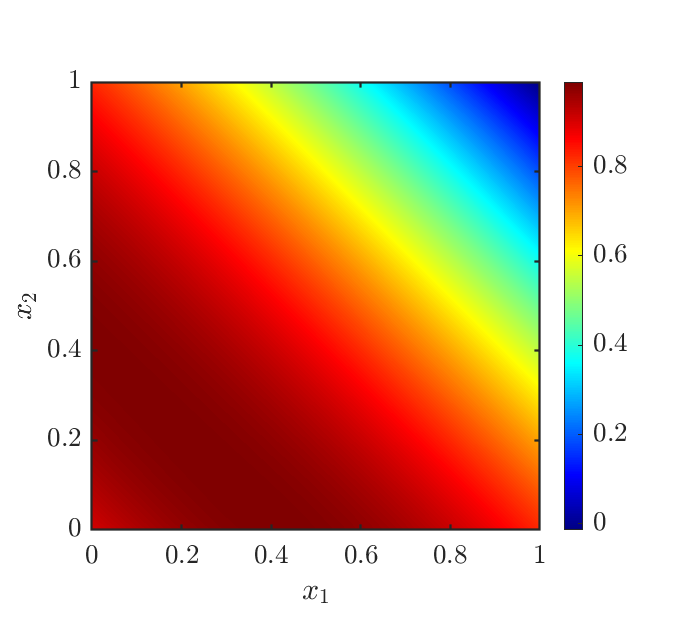}
        \end{minipage}%
    }}

    \vspace{-0.30em}

    % --------------------------------------------------------
    % Solid displacement
    % --------------------------------------------------------
    {\subfloat[{$t=0.2\,\mathrm{s}$}\label{fig:astorino_dh_dd_t2}]{%
        \begin{minipage}[t]{0.34\textwidth}
        \centering
        \includegraphics[width=\linewidth]
        {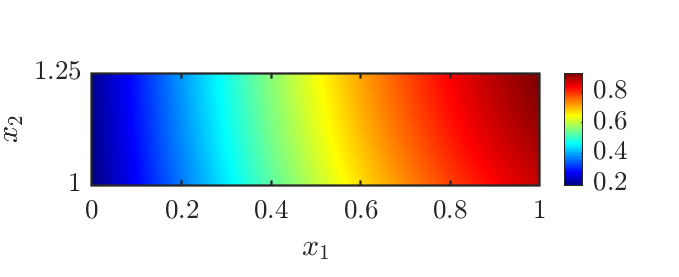}
        \end{minipage}%
    }}
    \hspace{0.04\textwidth}
    {\subfloat[{$t=10\,\mathrm{s}$}\label{fig:astorino_dh_dd_t10}]{%
        \begin{minipage}[t]{0.34\textwidth}
        \centering
        \includegraphics[width=\linewidth]
        {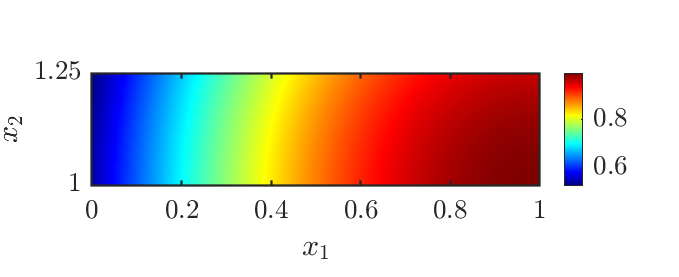}
        \end{minipage}%
    }}

    \vspace{-0.30em}

    % --------------------------------------------------------
    % Fluid-velocity error
    % --------------------------------------------------------
    {\subfloat[{$t=0.2\,\mathrm{s}$}\label{fig:astorino_dh_w_err_t2}]{%
        \begin{minipage}[t]{0.34\textwidth}
        \centering
        \includegraphics[width=\linewidth]
        {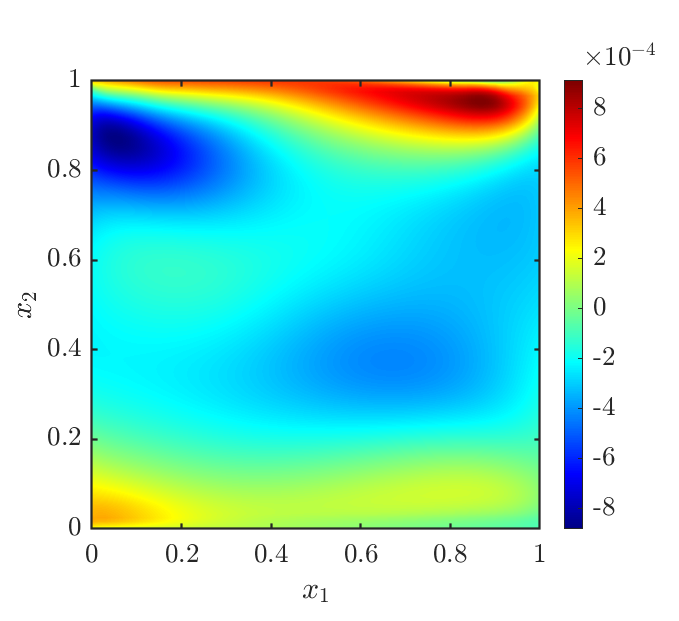}
        \end{minipage}%
    }}
    \hspace{0.04\textwidth}
    {\subfloat[{$t=10\,\mathrm{s}$}\label{fig:astorino_dh_w_err_t10}]{%
        \begin{minipage}[t]{0.34\textwidth}
        \centering
        \includegraphics[width=\linewidth]
        {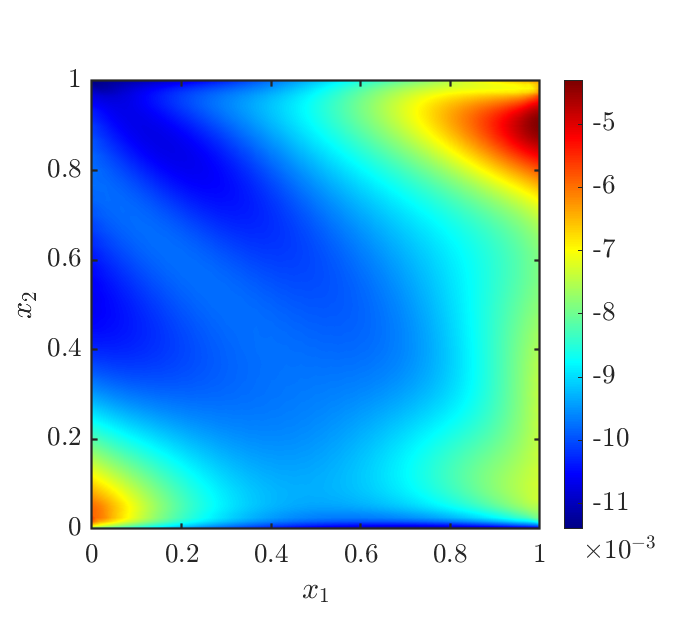}
        \end{minipage}%
    }}

    \vspace{-0.30em}

    % --------------------------------------------------------
    % Solid-displacement error
    % --------------------------------------------------------
    {\subfloat[{$t=0.2\,\mathrm{s}$}\label{fig:astorino_dh_dd_err_t2}]{%
        \begin{minipage}[t]{0.34\textwidth}
        \centering
        \includegraphics[width=\linewidth]
        {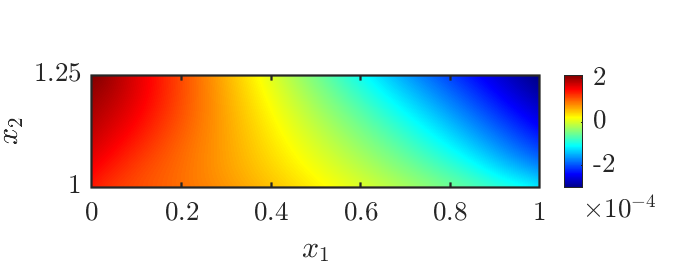}
        \end{minipage}%
    }}
    \hspace{0.04\textwidth}
    {\subfloat[{$t=10\,\mathrm{s}$}\label{fig:astorino_dh_dd_err_t10}]{%
        \begin{minipage}[t]{0.34\textwidth}
        \centering
        \includegraphics[width=\linewidth]
        {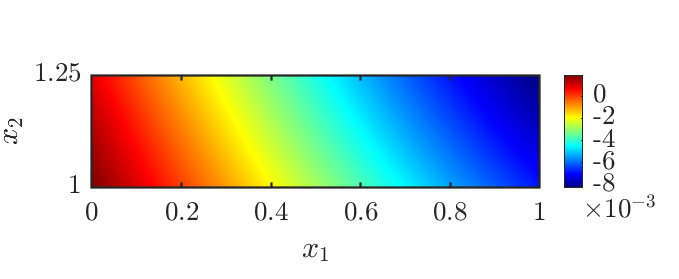}
        \end{minipage}%
    }}

    \caption{Astorino--Grandmont problem; IA-DMDc \(\mid\) IA-DMDc
    predictions and pointwise errors. From top to bottom, the rows show
    fluid velocity, solid displacement, fluid-velocity error, and
    solid-displacement error; the columns correspond to the two indicated
    times.}
    \label{fig:astorino_dh_predictions}
\end{figure}

% ============================================================
% Relative errors
% ============================================================

\begin{figure}[H]
    \centering

    {\subfloat[{Fluid velocity}\label{fig:astorino_w_re}]{%
        \begin{minipage}[t]{0.39\textwidth}
        \centering
        \includegraphics[width=\linewidth]
        {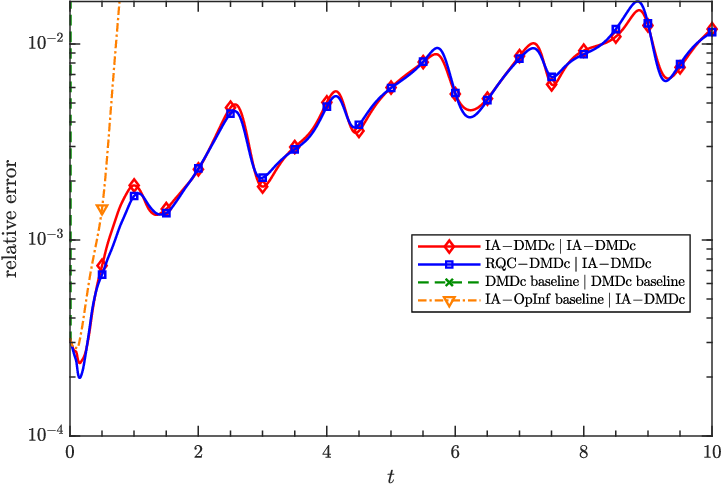}
        \end{minipage}%
    }}
    \hspace{0.04\textwidth}
    {\subfloat[{Solid displacement}\label{fig:astorino_dd_re}]{%
        \begin{minipage}[t]{0.39\textwidth}
        \centering
        \includegraphics[width=\linewidth]
        {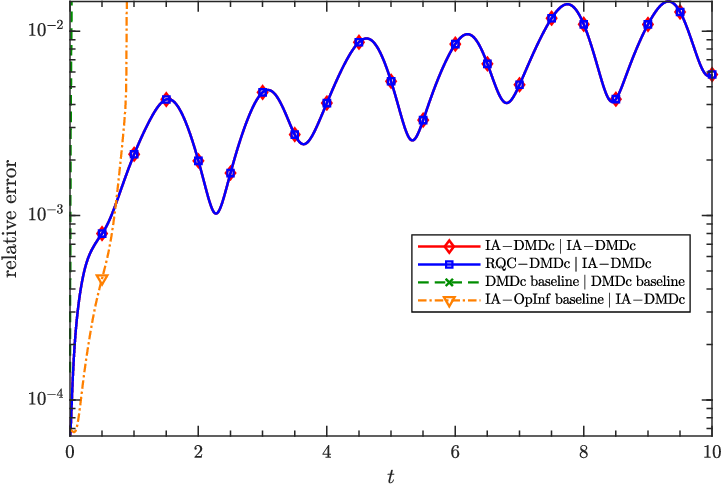}
        \end{minipage}%
    }}

    \caption{Astorino--Grandmont problem; relative errors of the fluid velocity
    and solid displacement. The training interval ends at
    $t=2.5\,\mathrm{s}$.}
    \label{fig:astorino_relative_errors}
\end{figure}

% ============================================================
% Kinetic energies
% ============================================================

\begin{figure}[H]
    \centering

    {\subfloat[{Fluid kinetic energy}\label{fig:astorino_w_kinetic_energy}]{%
        \begin{minipage}[t]{0.39\textwidth}
        \centering
        \includegraphics[width=\linewidth]
        {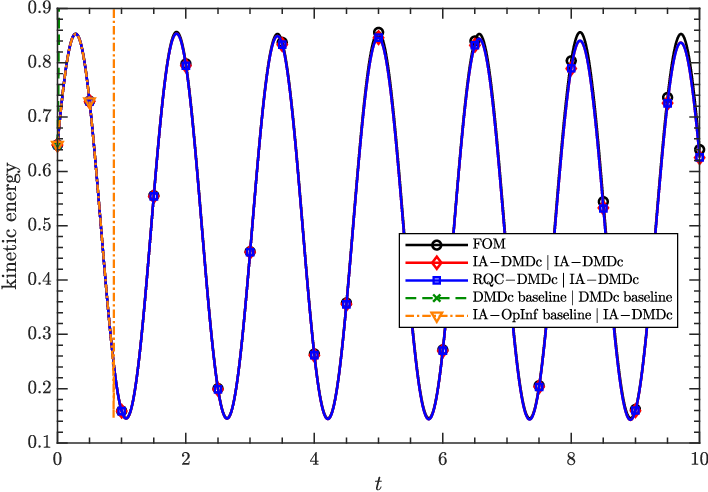}
        \end{minipage}%
    }}
    \hspace{0.04\textwidth}
    {\subfloat[{Solid kinetic energy}\label{fig:astorino_dd_kinetic_energy}]{%
        \begin{minipage}[t]{0.39\textwidth}
        \centering
        \includegraphics[width=\linewidth]
        {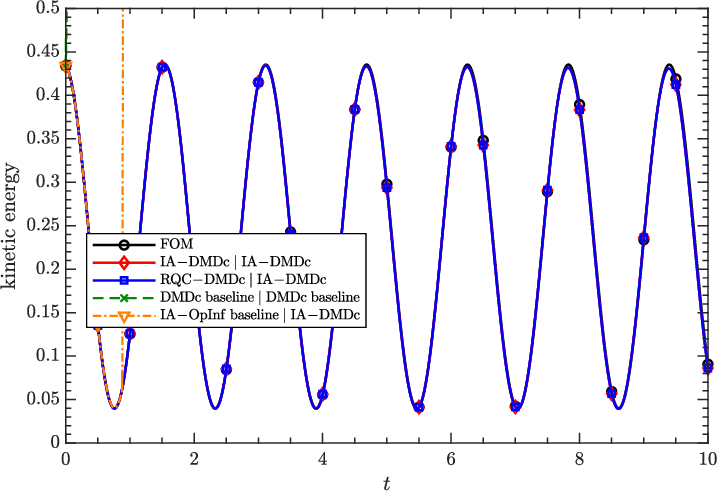}
        \end{minipage}%
    }}

    \caption{Astorino--Grandmont problem; fluid and solid kinetic energies.}
    \label{fig:astorino_kinetic_energies}
\end{figure}

Figure~\ref{fig:astorino_dh_predictions} presents the fields and pointwise errors obtained with IA-DMDc \(\mid\) IA-DMDc. The RQC-DMDc \(\mid\) IA-DMDc configuration produces similar fields and error distributions and is therefore omitted from this figure. Both configurations reproduce the fluid velocity and solid displacement at \(t=0.2\,\mathrm{s}\), with pointwise errors of order \(10^{-4}\). At \(t=10\,\mathrm{s}\), the spatial phase and amplitude remain accurately represented, with pointwise errors of order \(10^{-3}\)--\(10^{-2}\).
The relative-error histories in Figure~\ref{fig:astorino_relative_errors} show comparable accuracy for the two interface-aware configurations throughout the prediction interval, with neither model uniformly more accurate for both fields. The fluid and solid kinetic energies in
Figure~\ref{fig:astorino_kinetic_energies} also closely follow the analytical histories in~\eqref{eq:astorino_exact_kinetic_energies}. 
This agreement demonstrates accurate reproduction of the energy oscillations beyond the training interval. The DMDc baseline and IA-OpInf baseline do not provide a comparable stable extrapolation. 
The similar performance of IA-DMDc and RQC-DMDc is consistent with the vanishing convective contribution along the manufactured trajectory. In this setting, the quadratic enrichment provides little additional accuracy, while its regularized identification retains the predictive
quality of the linear interface-aware model.

\subsection{Freely oscillating cylinder}
\label{subsec:free_oscillating_cylinder}

We study incompressible flow past a circular cylinder constrained to move
 in the vertical direction. At \(\mathrm{Re}=100\), the
separated wake develops a von K\'arm\'an vortex street. Synchronization of the
vortex-shedding and structural frequencies can then produce sustained
vortex-induced vibration in the lock-in regime~\cite{williamson2004vortex}.
The computational domain is \([-3D,9D]\times[-3D,3D]\), and the cylinder is
initially centered at the origin. We set
\[
D = 1, \qquad
\rho_{\mathrm f} = 1, \qquad
\mu_{\mathrm f} = 0.008, \qquad
U_{\mathrm f} = 0.8, \qquad
\mathrm{Re}=\frac{\rho_{\mathrm f} U_{\mathrm f}D}{\mu_{\mathrm f}}=100.
\]
We prescribe uniform inflow on the left boundary, zero pressure on the right boundary, and homogeneous Neumann conditions on the horizontal boundaries. The moving fluid mesh is obtained by harmonic extension. The cylinder is constrained to move only in the vertical direction and is
governed by the mass--spring system
\[
m_{\mathrm s} \ddot d_{\mathrm y}+k_{\mathrm s} d_{\mathrm y}=F_{\mathrm y},
\]
where \(d_{\mathrm y}\) is the cylinder displacement and \(F_{\mathrm y}\) is the vertical fluid
force. For the mass ratio \(m^*=m_{\mathrm s}/(\rho_{\mathrm f}D^2)\), we choose
$
m^* = 4, k_{\mathrm s} = 2.13,
$
which places the coupled system close to lock-in. We monitor the vertical displacement, drag and lift forces, and fluid kinetic energy. The discretization contains \(N_{\mathrm f}=N_{\mathrm u}+N_{\mathrm p}=7128\) fluid degrees of freedom, with \(N_{\mathrm u}=6314\) and \(N_{\mathrm p}=814\), one structural degree of freedom, and \(N_{\mathrm g}=120\) interface-control degrees of freedom. Snapshots are sampled every \(\Delta t_{\mathrm s}=0.075\,\mathrm{s}\).
We first run a fixed-cylinder calculation until vortex shedding is fully
developed and then release the cylinder (the release time is set to \(t=0\)). The first 2500 snapshots,
\(t\in(0,187.5]\,\mathrm{s}\), are used for training, and prediction is
continued to \(t=750\,\mathrm{s}\). The fluid and interface ranks are
\(r_{\mathrm u}=270\), \(r_{\mathrm p}=30\), and \(r_{\mathrm g}=20\). The cylinder has only one structural
degree of freedom, so its mass--spring equation is retained without reduction.

We compare IA-DMDc and RQC-DMDc for the fluid. For RQC-DMDc, we use \(\omega_{\mathrm H}=10^{-2}\), corresponding to \(\rho_Q^{\mathrm{train}}=9.097\times10^{-5}\). 
The selected reduced dimension leads to a large quadratic feature space relative to the available training data. This strongly underdetermined regression tests whether the normalization and block-selective regularization of RQC-DMDc can identify a useful nonlinear correction, a central issue in regularized operator inference \cite{mcquarrie2021regularised,sawant2023physics}.  The IA-OpInf baseline and 
DMDc baseline are included for comparison.
In the field plots, the gray and black circles mark the equilibrium and
instantaneous cylinder positions, respectively; the cross and plus symbols
mark their centers, and \(\Delta y_c\) denotes the vertical displacement.

% ============================================================
% IA-DMDc prediction and pointwise error
% ============================================================

\begin{figure}[H]
    \centering

    {\subfloat[{\(t=3\,\mathrm{s}\)}]{%
        \begin{minipage}[b]{0.30\textwidth}
        \centering
        \includegraphics[
        width=\linewidth,
        trim={25 5 15 5},
        clip
        ]{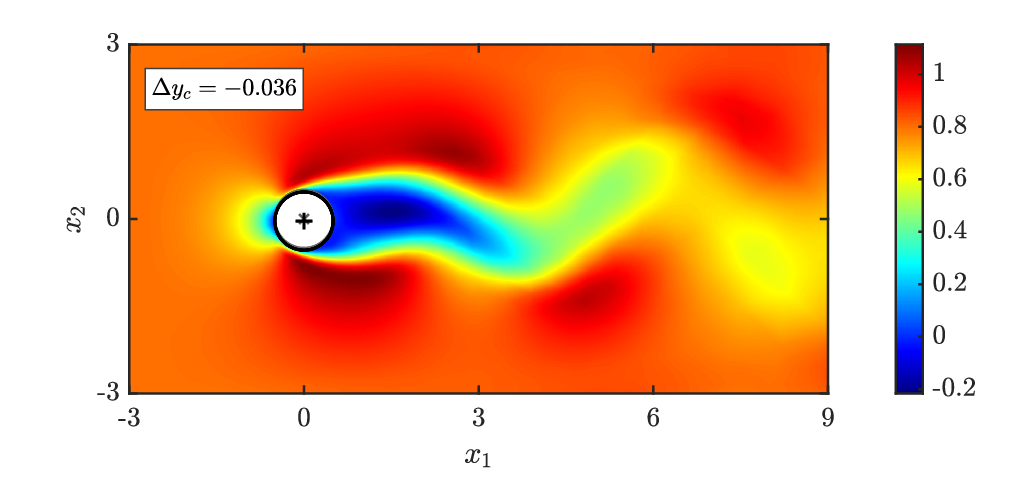}
        \end{minipage}%
    }}%
    \hspace{0.02\textwidth}%
    {\subfloat[{\(t=300\,\mathrm{s}\)}]{%
        \begin{minipage}[b]{0.30\textwidth}
        \centering
        \includegraphics[
        width=\linewidth,
        trim={25 5 15 5},
        clip
        ]{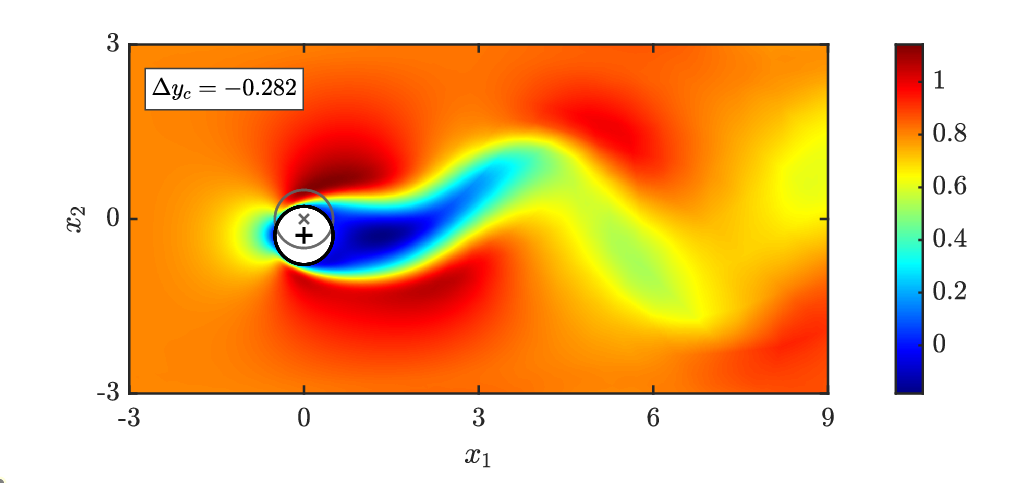}
        \end{minipage}%
    }}%
    \hspace{0.02\textwidth}%
    {\subfloat[{\(t=750\,\mathrm{s}\)}]{%
        \begin{minipage}[b]{0.30\textwidth}
        \centering
        \includegraphics[
        width=\linewidth,
        trim={25 5 15 5},
        clip
        ]{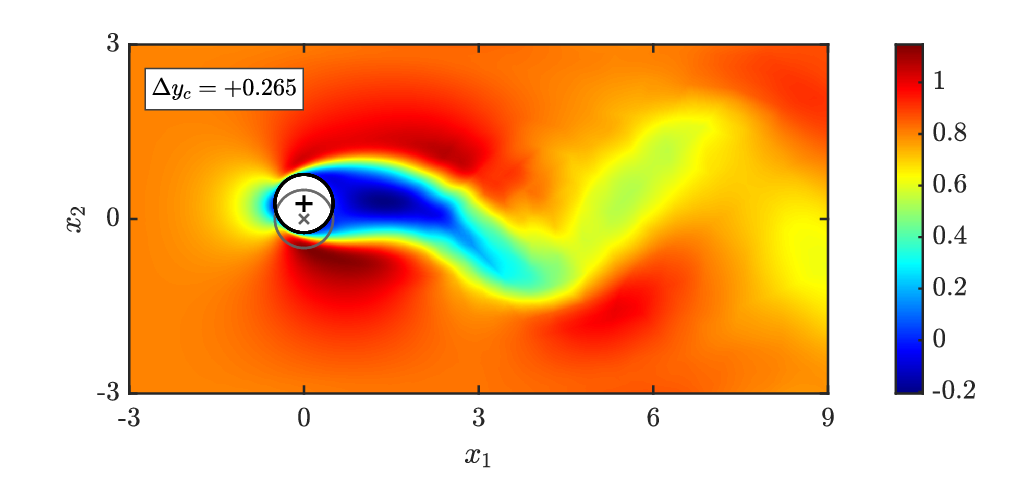}
        \end{minipage}%
    }}

    \vspace{-0.2em}

    {\subfloat[{\(t=3\,\mathrm{s}\)}]{%
        \begin{minipage}[b]{0.30\textwidth}
        \centering
        \includegraphics[
        width=\linewidth,
        trim={25 5 15 5},
        clip
        ]{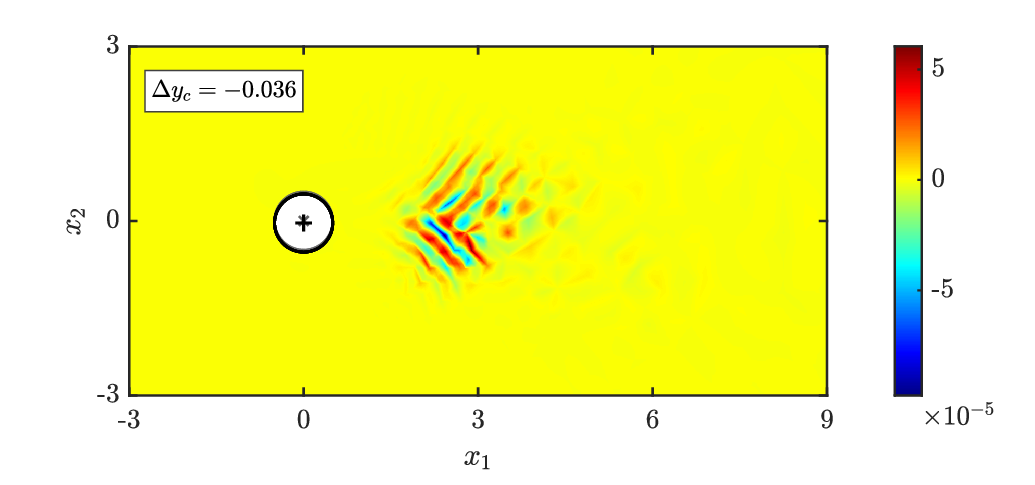}
        \end{minipage}%
    }}%
    \hspace{0.02\textwidth}%
    {\subfloat[{\(t=300\,\mathrm{s}\)}]{%
        \begin{minipage}[b]{0.30\textwidth}
        \centering
        \includegraphics[
        width=\linewidth,
        trim={25 5 15 5},
        clip
        ]{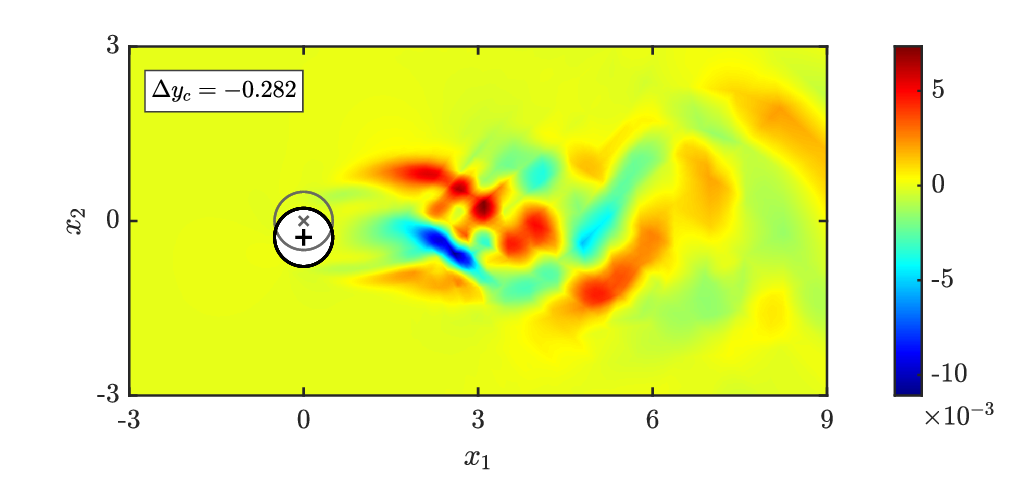}
        \end{minipage}%
    }}%
    \hspace{0.02\textwidth}%
    {\subfloat[{\(t=750\,\mathrm{s}\)}]{%
        \begin{minipage}[b]{0.30\textwidth}
        \centering
        \includegraphics[
        width=\linewidth,
        trim={25 5 15 5},
        clip
        ]{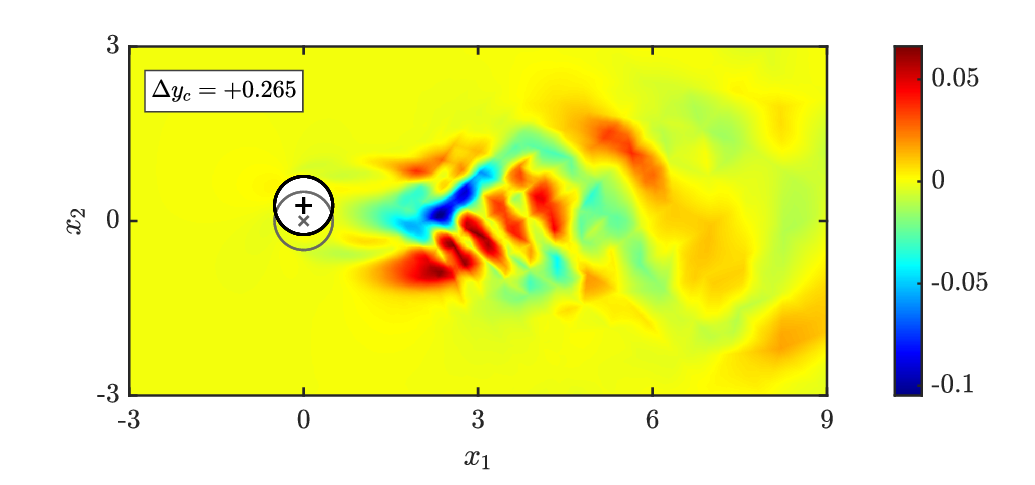}
        \end{minipage}%
    }}

    \caption{Freely oscillating cylinder; IA-DMDc prediction of the
    streamwise velocity and the corresponding pointwise error. The first row
    shows the IA-DMDc prediction, and the second row shows the error relative to
    the reference full-order solution.}
    \label{fig:oscillating_cylinder_dmdc_velocity_error}
\end{figure}

% ============================================================
% RQC-DMDc prediction and pointwise error
% ============================================================

\begin{figure}[H]
    \centering

    {\subfloat[{\(t=3\,\mathrm{s}\)}]{%
        \begin{minipage}[b]{0.30\textwidth}
        \centering
        \includegraphics[
        width=\linewidth,
        trim={25 5 15 10},
        clip
        ]{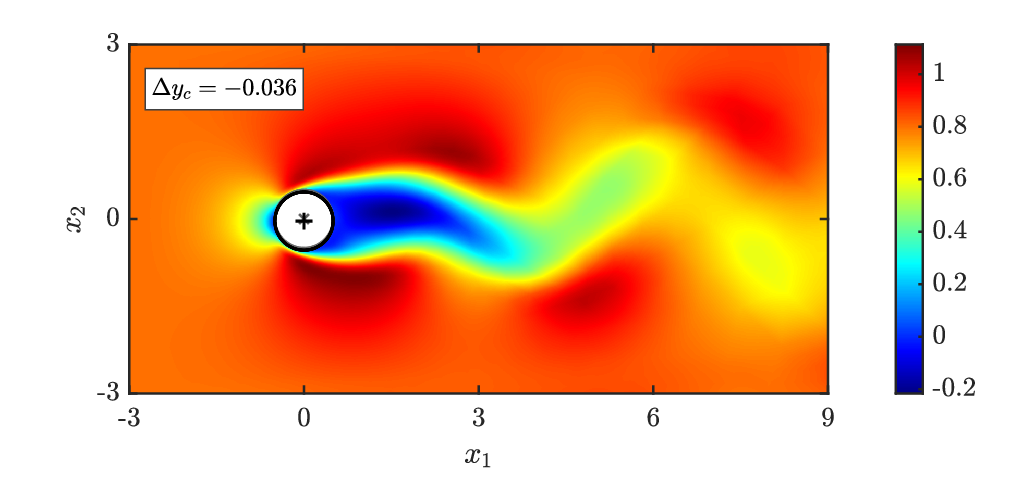}
        \end{minipage}%
    }}%
    \hspace{0.02\textwidth}%
    {\subfloat[{\(t=300\,\mathrm{s}\)}]{%
        \begin{minipage}[b]{0.30\textwidth}
        \centering
        \includegraphics[
        width=\linewidth,
        trim={25 5 15 10},
        clip
        ]{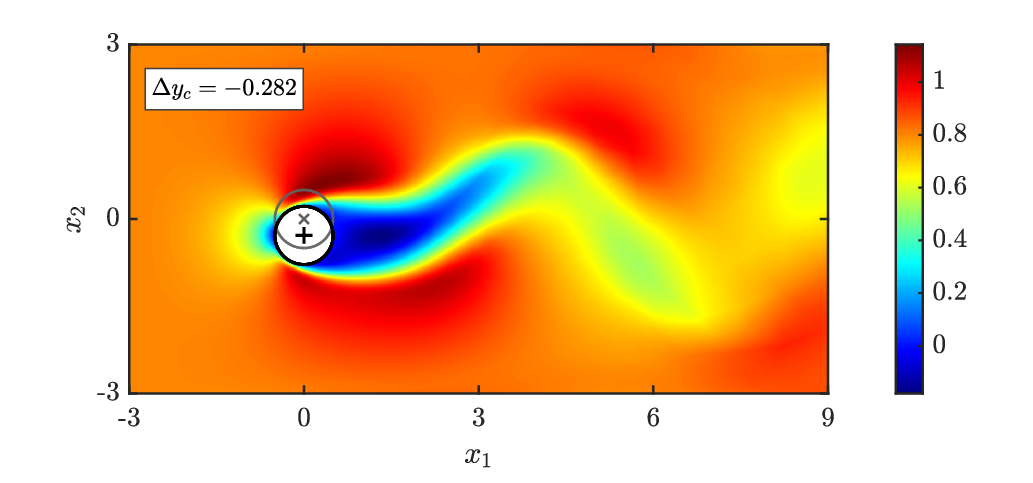}
        \end{minipage}%
    }}%
    \hspace{0.02\textwidth}%
    {\subfloat[{\(t=750\,\mathrm{s}\)}]{%
        \begin{minipage}[b]{0.30\textwidth}
        \centering
        \includegraphics[
        width=\linewidth,
        trim={25 5 15 10},
        clip
        ]{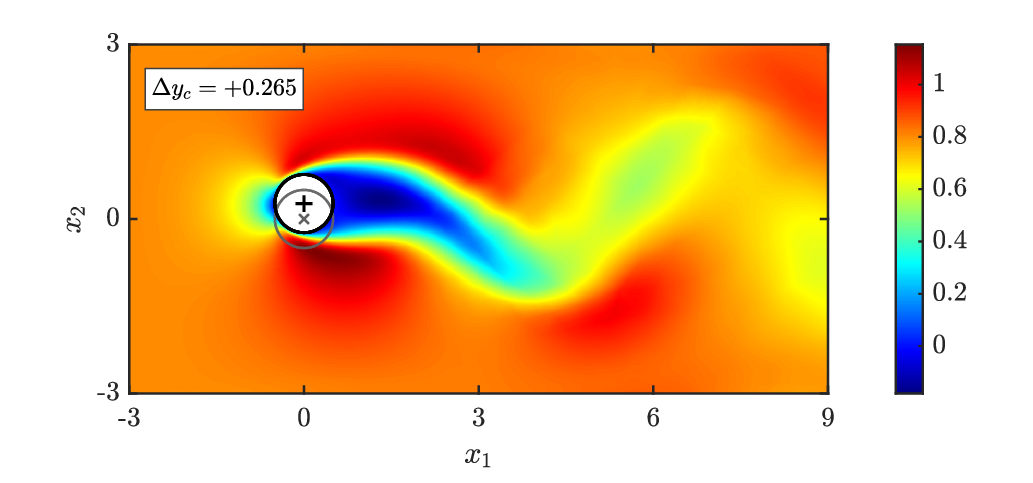}
        \end{minipage}%
    }}

    \vspace{-0.2em}

    {\subfloat[{\(t=3\,\mathrm{s}\)}]{%
        \begin{minipage}[b]{0.30\textwidth}
        \centering
        \includegraphics[
        width=\linewidth,
        trim={25 5 15 10},
        clip
        ]{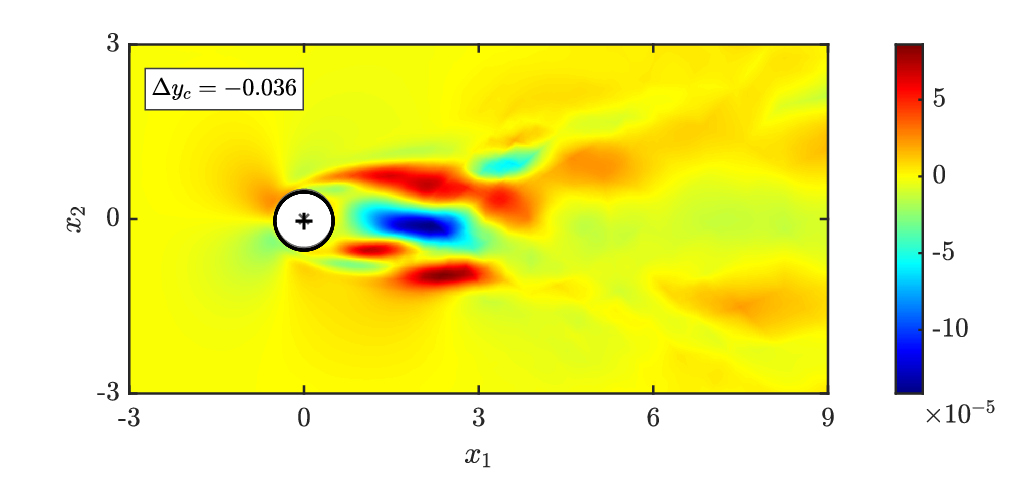}
        \end{minipage}%
    }}%
    \hspace{0.02\textwidth}%
    {\subfloat[{\(t=300\,\mathrm{s}\)}]{%
        \begin{minipage}[b]{0.30\textwidth}
        \centering
        \includegraphics[
        width=\linewidth,
        trim={25 5 15 10},
        clip
        ]{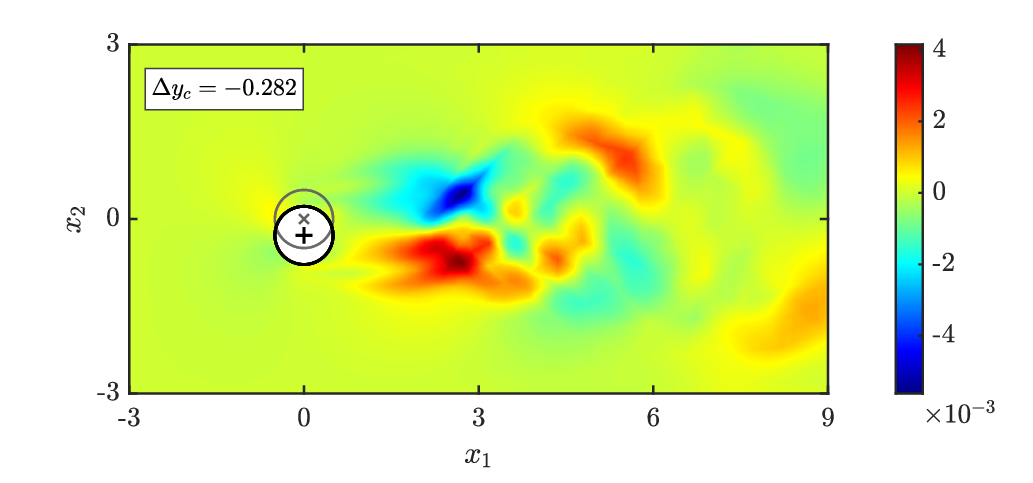}
        \end{minipage}%
    }}%
    \hspace{0.02\textwidth}%
    {\subfloat[{\(t=750\,\mathrm{s}\)}]{%
        \begin{minipage}[b]{0.30\textwidth}
        \centering
        \includegraphics[
        width=\linewidth,
        trim={25 5 15 10},
        clip
        ]{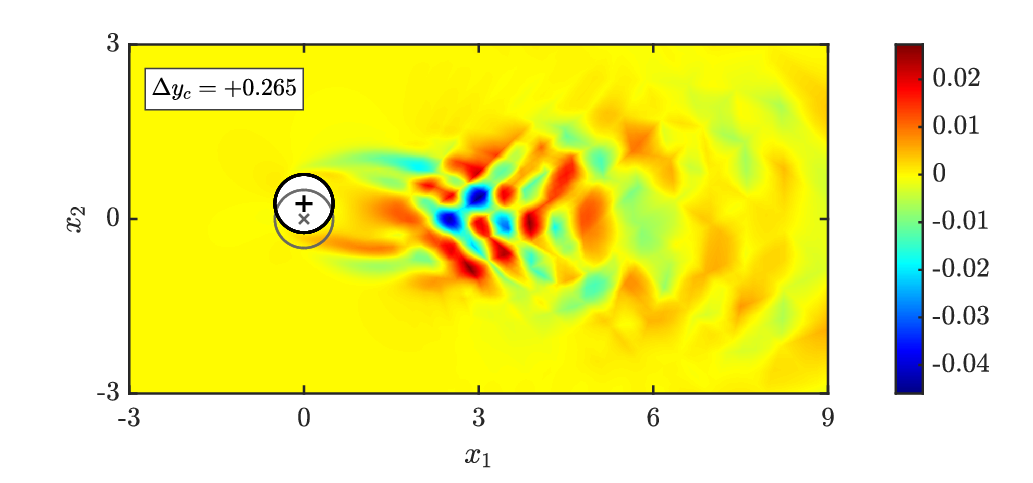}
        \end{minipage}%
    }}

    \caption{Freely oscillating cylinder; RQC-DMDc prediction of the
    streamwise velocity and the corresponding pointwise error. The first row
    shows the RQC-DMDc prediction, and the second row shows the error relative to
    the reference full-order solution.}
    \label{fig:oscillating_cylinder_rqc_dmdc_velocity_error}
\end{figure}

% ============================================================
% Relative error
% ============================================================

\begin{figure}[H]
    \centering
    \includegraphics[width=0.43\textwidth]
    {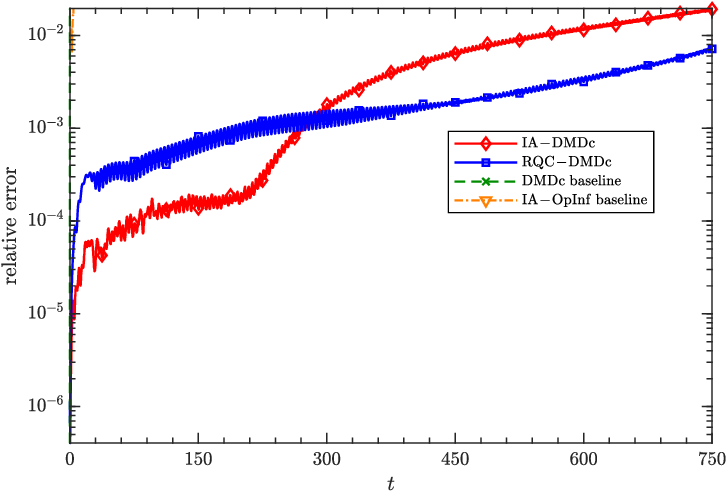}
    \caption{Freely oscillating cylinder; relative error of the
    streamwise-velocity prediction. The training interval ends at
    \(t=187.5\,\mathrm{s}\).}
    \label{fig:oscillating_cylinder_w_re}
\end{figure}

% ============================================================
% Cylinder vertical displacement
% ============================================================

\begin{figure}[H]
    \centering

    {\subfloat[{Full interval, \(t\in(0,750]\,\mathrm{s}\)}]{%
        \begin{minipage}[b]{0.39\textwidth}
        \centering
        \includegraphics[width=\linewidth]
        {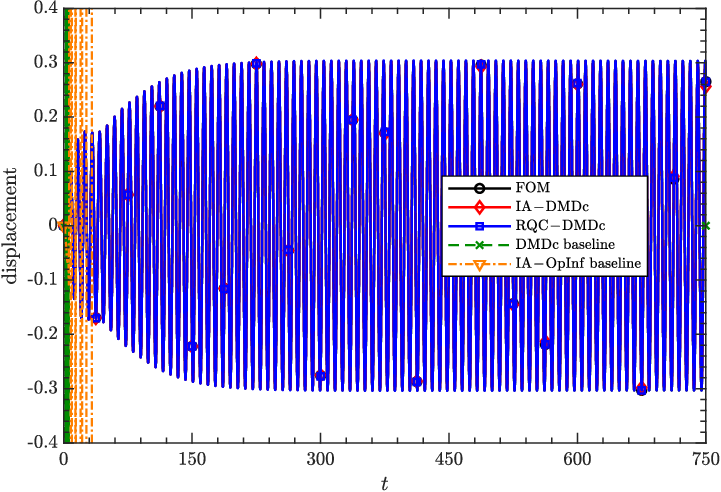}
        \end{minipage}%
    }}
    \hspace{0.04\textwidth}
    {\subfloat[{Final five seconds, \(t\in[745,750]\,\mathrm{s}\)}]{%
        \begin{minipage}[b]{0.39\textwidth}
        \centering
        \includegraphics[width=\linewidth]
        {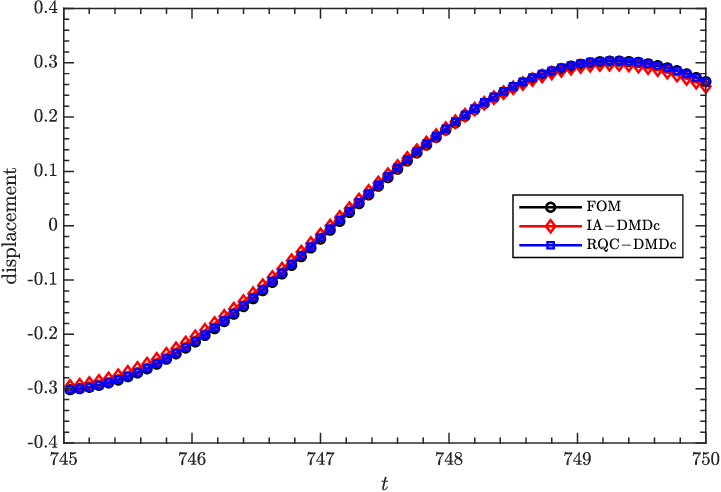}
        \end{minipage}%
    }}

    \caption{Freely oscillating cylinder; vertical (cross-flow) displacement
    over the full prediction interval and over the final five seconds.}
    \label{fig:oscillating_cylinder_displacement}
\end{figure}

% ============================================================
% Drag force
% ============================================================

\begin{figure}[H]
    \centering

    {\subfloat[{Full interval, $t\in(0,750]\,\mathrm{s}$}]{%
        \begin{minipage}[t]{0.39\textwidth}
        \centering
        \includegraphics[width=\linewidth]
        {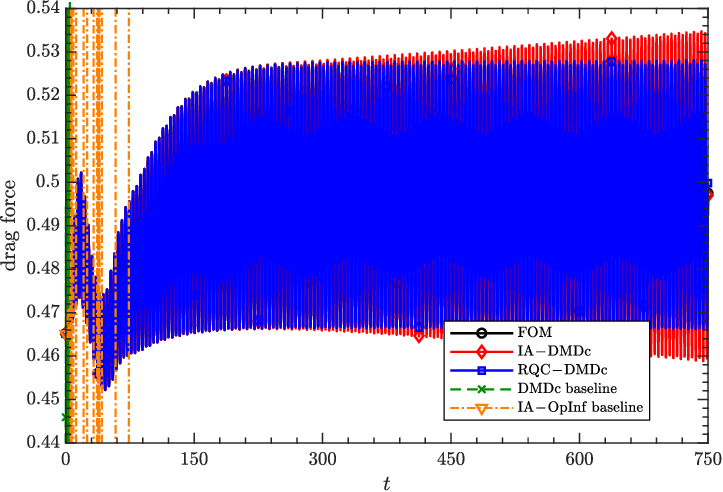}
        \end{minipage}%
    }}
    \hspace{0.04\textwidth}
    {\subfloat[{Final five seconds, $t\in[745,750]\,\mathrm{s}$}]{%
        \begin{minipage}[t]{0.39\textwidth}
        \centering
        \includegraphics[width=\linewidth]
        {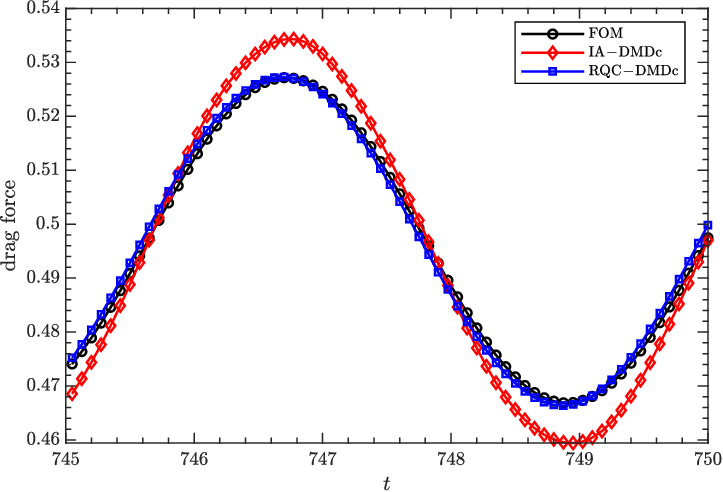}
        \end{minipage}%
    }}

    \caption{Freely oscillating cylinder; drag force exerted by the fluid over
    the full prediction interval and the final five seconds.}
    \label{fig:oscillating_cylinder_drag_force}
\end{figure}

% ============================================================
% Lift force
% ============================================================

\begin{figure}[H]
    \centering

    {\subfloat[{Full interval, $t\in(0,750]\,\mathrm{s}$}]{%
        \begin{minipage}[t]{0.39\textwidth}
        \centering
        \includegraphics[width=\linewidth]
        {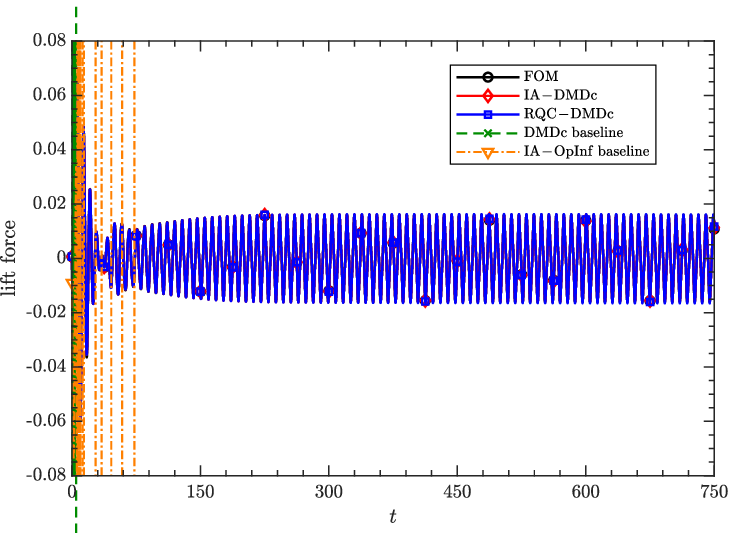}
        \end{minipage}%
    }}
    \hspace{0.04\textwidth}
    {\subfloat[{Final five seconds, $t\in[745,750]\,\mathrm{s}$}]{%
        \begin{minipage}[t]{0.39\textwidth}
        \centering
        \includegraphics[width=\linewidth]
        {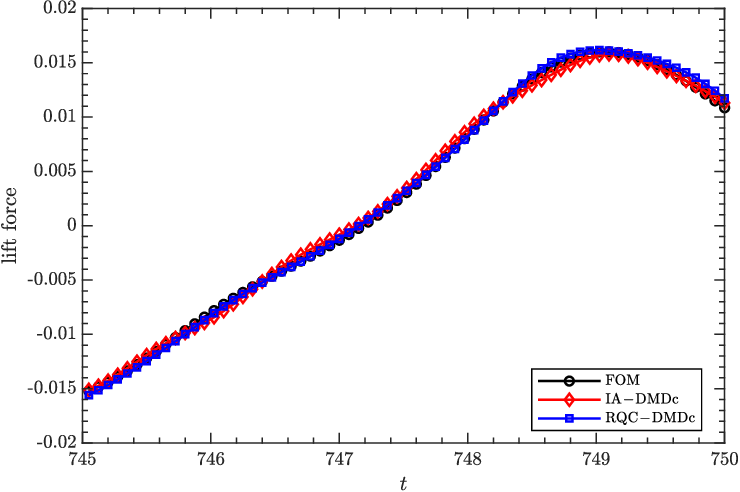}
        \end{minipage}%
    }}

    \caption{Freely oscillating cylinder; lift force exerted by the fluid over
    the full prediction interval and the final five seconds.}
    \label{fig:oscillating_cylinder_lift_force}
\end{figure}

% ============================================================
% Fluid kinetic energy
% ============================================================

\begin{figure}[H]
    \centering

    {\subfloat[{Full interval, $t\in(0,750]\,\mathrm{s}$}]{%
        \begin{minipage}[t]{0.39\textwidth}
        \centering
        \includegraphics[width=\linewidth]
        {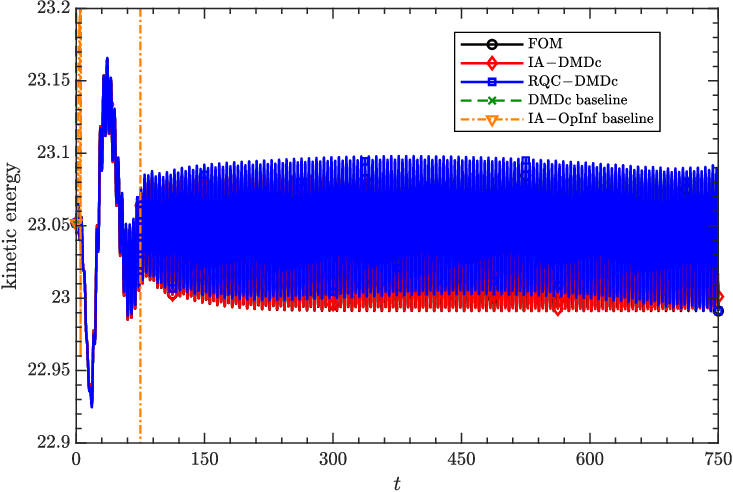}
        \end{minipage}%
    }}
    \hspace{0.04\textwidth}
    {\subfloat[{Final five seconds, $t\in[745,750]\,\mathrm{s}$}]{%
        \begin{minipage}[t]{0.39\textwidth}
        \centering
        \includegraphics[width=\linewidth]
        {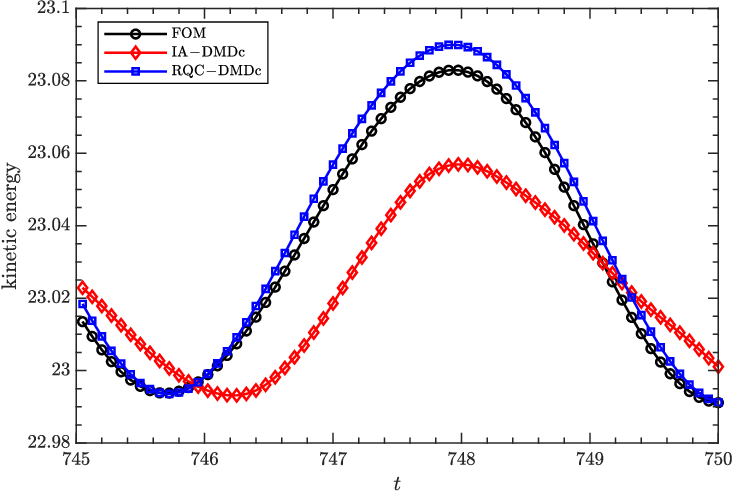}
        \end{minipage}%
    }}

    \caption{Freely oscillating cylinder; fluid kinetic energy over the full
    prediction interval and the final five seconds.}
    \label{fig:oscillating_cylinder_fluid_kinetic_energy}
\end{figure}

Figures~\ref{fig:oscillating_cylinder_dmdc_velocity_error}-\ref{fig:oscillating_cylinder_rqc_dmdc_velocity_error} compare the streamwise-velocity predictions and pointwise errors. Both models reproduce the wake during training.
At \(t=300\,\mathrm{s}\) and \(t=750\,\mathrm{s}\), RQC-DMDc more accurately predicts the position and strength of the alternating wake structures, whereas IA-DMDc develops
larger accumulated discrepancies.
Figure~\ref{fig:oscillating_cylinder_w_re} distinguishes training accuracy from long-time predictive accuracy. IA-DMDc has the smaller relative error during training and
the initial part of extrapolation. Its error subsequently grows more rapidly, and RQC-DMDc
becomes more accurate over the later prediction interval, retaining this advantage through \(t=750\,\mathrm{s}\).

Both models reproduce the sustained vertical cylinder oscillation in Figure~\ref{fig:oscillating_cylinder_displacement}. The drag and lift histories in
Figures~\ref{fig:oscillating_cylinder_drag_force} and~\ref{fig:oscillating_cylinder_lift_force} show that RQC-DMDc more closely follows the full-order force amplitudes and phase.
The fluid kinetic energy in Figure~\ref{fig:oscillating_cylinder_fluid_kinetic_energy} provides a consistent distinction: the RQC-DMDc curve follows the phase and amplitude of the reference energy more closely over the final oscillation cycles.
The DMDc and IA-OpInf baselines become unstable near the start of the calculation and do not produce usable long-time trajectories. RQC-DMDc maintains bounded and accurate predictions despite the strongly underdetermined quadratic regression. Its comparison with IA-OpInf shows that its block-selective regularization can extract a useful nonlinear correction without destabilizing recursive prediction.

\subsection{Turek problem}
\label{subsec:turek_benchmark}

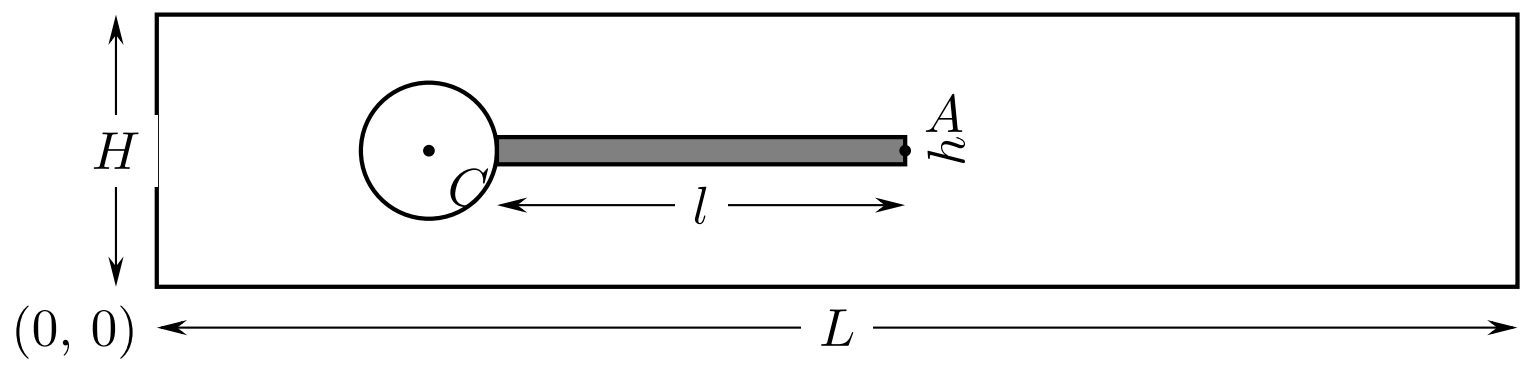
\begin{figure}[H]
\centering

\begin{tikzpicture}[scale=0.9, >={Stealth[length=2mm, width=1mm]}]  % Customize global arrow tip style
    % Outer rectangle
    \draw[thick] (0,0) rectangle (10,2);  % Full outer boundary

    % Labels for the dimensions of the outer rectangle
    \draw[<->] (-0.3,0) -- (-0.3,2) node[midway, fill=white] {$H$};
    \draw[<->] (0,-0.3) -- (10,-0.3) node[midway, fill=white] {$L$};

    % Inner circle representing the cylinder
    \draw[thick] (2,1) circle (0.5);  % Centered at (2,1) with a radius of 0.5
    
    % Mark the center of the circle and label it
    \filldraw [black] (2,1) circle (1pt) node[below right] {$C$};

    % Inner rectangle parameters
    \pgfmathsetmacro{\rectHeight}{0.2} % Height for the rectangle
    \pgfmathsetmacro{\rectWidth}{3} % Width of the rectangle
    \pgfmathsetmacro{\rectStartY}{1 - \rectHeight/2} % Vertical start position

    % Fill and draw the inner rectangle without the left side
    \fill[gray] (2.5, \rectStartY) rectangle (2.5 + \rectWidth, \rectStartY + \rectHeight);
    \draw[thick] (2.5, \rectStartY) rectangle (2.5 + \rectWidth, \rectStartY + \rectHeight);

    % Arrow indicating the height 'h' of the inner rectangle
    \draw[<->]  % Custom arrow tip for this line
        (2.5 + \rectWidth + 0.3, \rectStartY) -- 
        (2.5 + \rectWidth + 0.3, \rectStartY + \rectHeight) node[midway,fill=white, rotate=90] {$h$};

\filldraw [black] (2.5 + \rectWidth, \rectStartY+\rectHeight/2) circle (1pt) node[above right] {$A$};

    % Arrow indicating the length 'l' of the inner rectangle
    \draw[<->]  % Another style of custom arrow tip
        (2.5, \rectStartY - 0.3) -- 
        (2.5 + \rectWidth, \rectStartY - 0.3) node[midway, fill=white] {$l$};

    % Positioning the diagram in a nice place
    \node at (0,0) [below left] {(0, 0)};
\end{tikzpicture}
\caption{Turek problem; geometric configuration (not to scale). The principal dimensions are $H=0.41$, $L=2.5$, $r=0.05$, $l=0.35$, and $h=0.02$, with cylinder center $C=(0.2,0.2)$. The reference position $A(0)=(0.6,0.2)$ is used to monitor the beam-tip displacement.}
\label{fig:turek_geometry}
\end{figure}

We consider laminar flow past an elastic body using a linear-elastic variant of the FSI3 benchmark in~\cite{turek2006proposal}. 
%The interaction of the vortex street behind a rigid cylinder with an attached elastic beam provides a demanding test of both field accuracy and long-time robustness. 
We retain the benchmark geometry and physical parameters, use the incompressible
Navier--Stokes equations for the fluid, and model the beam with \eqref{eq:linear_elastic_solid}, setting \(c_0 = f_{\mathrm s}=0\). This structural law is the small-strain linearization of the
St.~Venant--Kirchhoff model used in the original benchmark.
Figure~\ref{fig:turek_geometry} shows the configuration. The channel has length \(L=2.5\) and height \(H=0.41\); the cylinder has radius \(r=0.05\) and center \(C=(0.2,0.2)\). The beam has length \(l=0.35\) and height \(h=0.02\).
Its upstream end is clamped to the cylinder, while the remaining three sides form the fluid--structure interface. We monitor the material point with reference position \(A(0)=(0.6,0.2)\) at the beam tip. No-slip conditions are imposed on the channel walls and the rigid cylinder. At the inlet, we prescribe
\begin{equation*}
    \boldsymbol{u}_{\mathrm{in}}(x_2)
    =
    \begin{bmatrix}
        \displaystyle
        4u_{\max}\frac{x_2(H-x_2)}{H^2}\\[1mm]
        0
    \end{bmatrix},
    \qquad
    u_{\max}=3,
    %\label{eq:turek_inlet_profile}
\end{equation*}
and a homogeneous Neumann condition is imposed at the outlet.

The physical parameters are
\begin{equation*}
    \rho_{\mathrm f} = 10^3,\qquad
    \rho_{\mathrm s} = 10^3,\qquad
    E = 5.6\times 10^6,\qquad
    \nu_{\mathrm s} = 0.4,\qquad
    \mu_{\mathrm f} = 1.
\end{equation*}
The Lam\'e coefficients in \eqref{eq:linear_elastic_solid} are
\(\mu_{\mathrm s}=E/[2(1+\nu_{\mathrm s})]\) and
\(\lambda_{\mathrm s}=E\nu_{\mathrm s}/[(1+\nu_{\mathrm s})(1-2\nu_{\mathrm s})]\).
The mean inlet velocity is \(\bar u=2\), and the cylinder diameter is
\(2r=0.1\), giving \(\mathrm{Re}=\rho_{\mathrm f}\bar u(2r)/\mu_{\mathrm f}=200\). We first hold
the beam fixed until a developed wake is obtained, then release it and
integrate the coupled system for \(T=12\,\mathrm{s}\). 
%The algebraic equations are scaled by \(\rho_{\mathrm f}\) to improve conditioning without changing the physical solution.
We use an unstructured triangular grid with local refinement near the
cylinder, beam, and wake. The discretized fluid state contains
$
    N_{\mathrm f}=N_{\mathrm u}+N_{\mathrm p}=10379,
    N_{\mathrm u}=9178,
    N_{\mathrm p}=1201
$
degrees of freedom; the structural and interface variables contain \(N_{\mathrm s}=294\)
and \(N_{\mathrm g}=162\) degrees of freedom, respectively. The comparatively small structural system is retained in full-order form, and only the fluid is reduced.
Snapshots are sampled every \(\Delta t_{\mathrm s}=0.01\,\mathrm{s}\). The first
\(N_{\mathrm{tr}}=300\) snapshots, \(t\in(0,3]\,\mathrm{s}\), are used for training, and
\(t\in(3,12]\,\mathrm{s}\) is reserved for out-of-sample prediction. The fluid
and interface ranks are \(r_{\mathrm u}=85\), \(r_{\mathrm p}=40\), and \(r_{\mathrm g}=25\). IA-DMDc and RQC-DMDc use the same data, reduced spaces, time discretization, and projected interface-input operator. 
RQC-DMDc uses \(\omega_{\mathrm H}=25\), giving \(\rho_Q^{\mathrm{train}}=9.249\times10^{-4}\).
The structural subsystem is retained in full-order form for both configurations. 
The quadratic regression is strongly underdetermined: the reduced quadratic feature space is much larger than the available training data. This case tests the predictive capability of the regularized quadratic model when identification is strongly underdetermined and the training interval covers only one quarter of the simulation horizon. 
The IA-OpInf baseline and DMDc baseline are included for comparison.

% ============================================================
% IA-DMDc prediction and pointwise error
% ============================================================

\begin{figure}[H]
    \centering

    {\subfloat[{$t=0.2\,\mathrm{s}$}]{%
        \begin{minipage}[b]{0.30\textwidth}
        \centering
        \includegraphics[
        width=\linewidth,
        trim={10 0 20 10},
        clip
        ]{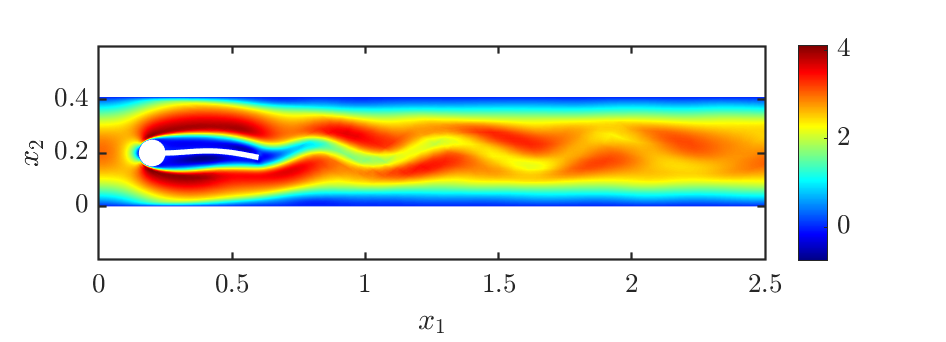}
        \end{minipage}%
    }}%
    \hspace{0.02\textwidth}%
    {\subfloat[{$t=3\,\mathrm{s}$}]{%
        \begin{minipage}[b]{0.30\textwidth}
        \centering
        \includegraphics[
        width=\linewidth,
        trim={10 0 20 10},
        clip
        ]{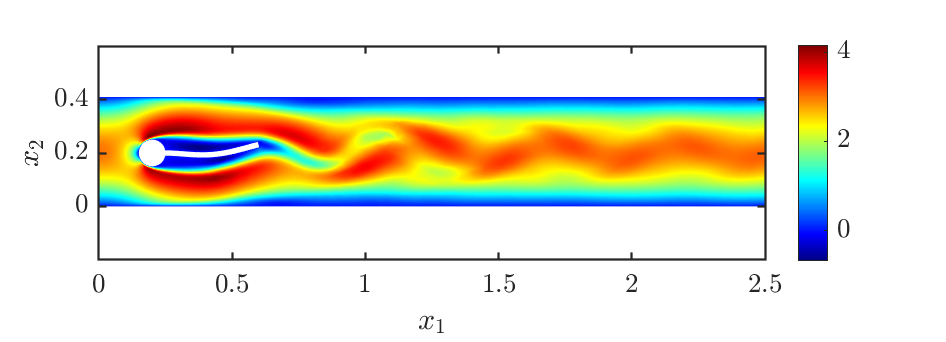}
        \end{minipage}%
    }}%
    \hspace{0.02\textwidth}%
    {\subfloat[{$t=12\,\mathrm{s}$}]{%
        \begin{minipage}[b]{0.30\textwidth}
        \centering
        \includegraphics[
        width=\linewidth,
        trim={10 0 20 10},
        clip
        ]{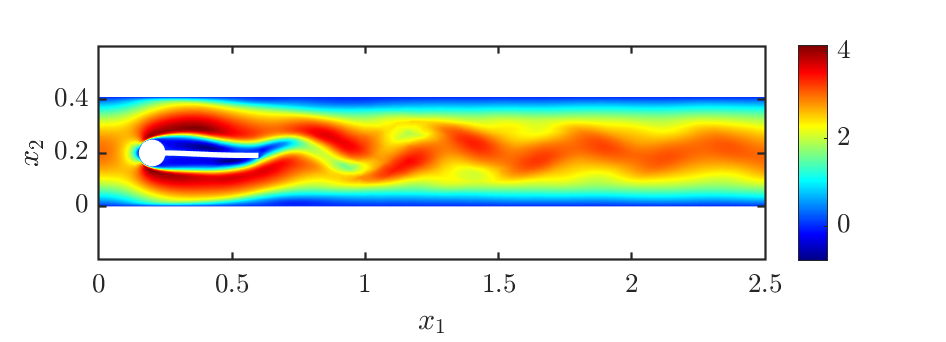}
        \end{minipage}%
    }}

    \vspace{-0.2em}

    {\subfloat[{$t=0.2\,\mathrm{s}$}]{%
        \begin{minipage}[b]{0.30\textwidth}
        \centering
        \includegraphics[
        width=\linewidth,
        trim={10 0 20 10},
        clip
        ]{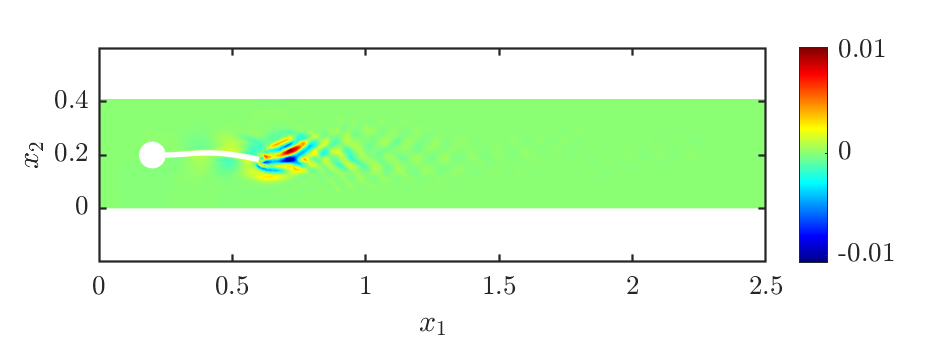}
        \end{minipage}%
    }}%
    \hspace{0.02\textwidth}%
    {\subfloat[{$t=3\,\mathrm{s}$}]{%
        \begin{minipage}[b]{0.30\textwidth}
        \centering
        \includegraphics[
        width=\linewidth,
        trim={10 0 20 10},
        clip
        ]{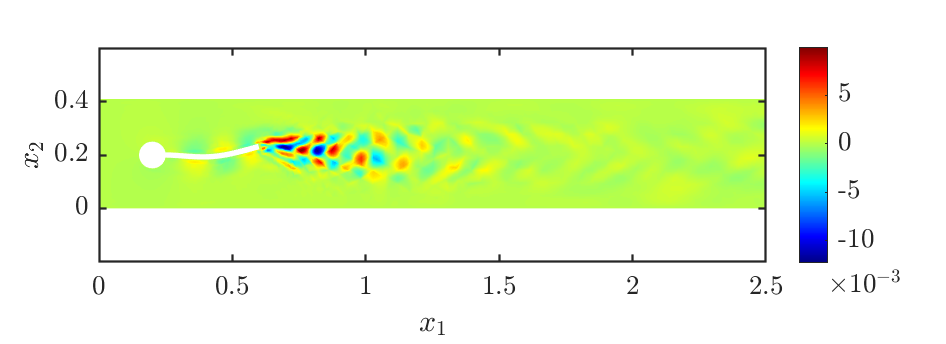}
        \end{minipage}%
    }}%
    \hspace{0.02\textwidth}%
    {\subfloat[{$t=12\,\mathrm{s}$}]{%
        \begin{minipage}[b]{0.30\textwidth}
        \centering
        \includegraphics[
        width=\linewidth,
        trim={10 0 20 10},
        clip
        ]{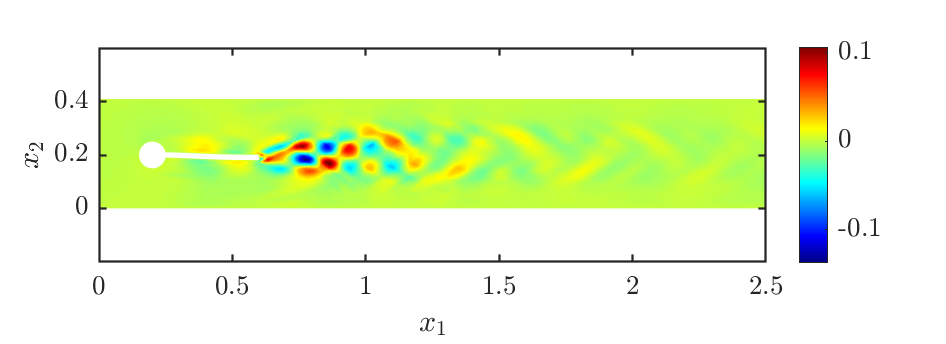}
        \end{minipage}%
    }}

    \caption{Turek problem; IA-DMDc prediction of the streamwise velocity
    (top) and pointwise error relative to the full-order solution (bottom).}
    \label{fig:turek_dmdc_field_error}
\end{figure}

% ============================================================
% RQC-DMDc prediction and pointwise error
% ============================================================

\begin{figure}[H]
    \centering

    {\subfloat[{$t=0.2\,\mathrm{s}$}]{%
        \begin{minipage}[b]{0.30\textwidth}
        \centering
        \includegraphics[
        width=\linewidth,
        trim={10 0 20 10},
        clip
        ]{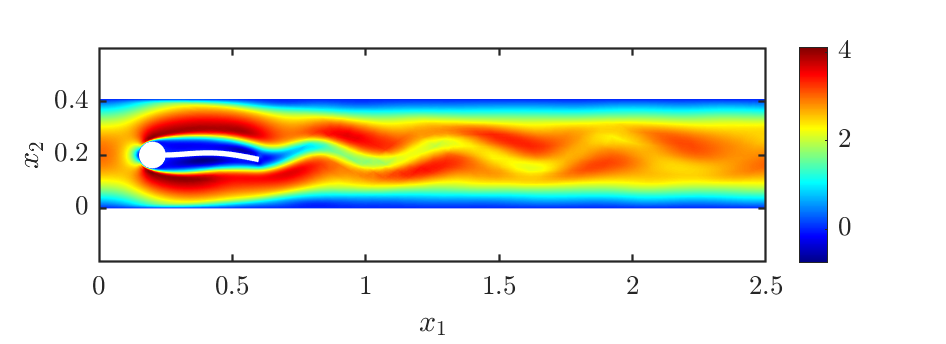}
        \end{minipage}%
    }}%
    \hspace{0.02\textwidth}%
    {\subfloat[{$t=3\,\mathrm{s}$}]{%
        \begin{minipage}[b]{0.30\textwidth}
        \centering
        \includegraphics[
        width=\linewidth,
        trim={10 0 20 10},
        clip
        ]{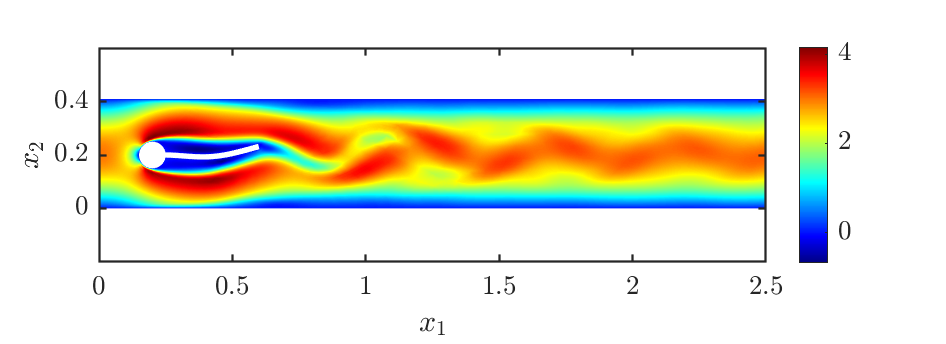}
        \end{minipage}%
    }}%
    \hspace{0.02\textwidth}%
    {\subfloat[{$t=12\,\mathrm{s}$}]{%
        \begin{minipage}[b]{0.30\textwidth}
        \centering
        \includegraphics[
        width=\linewidth,
        trim={10 0 20 10},
        clip
        ]{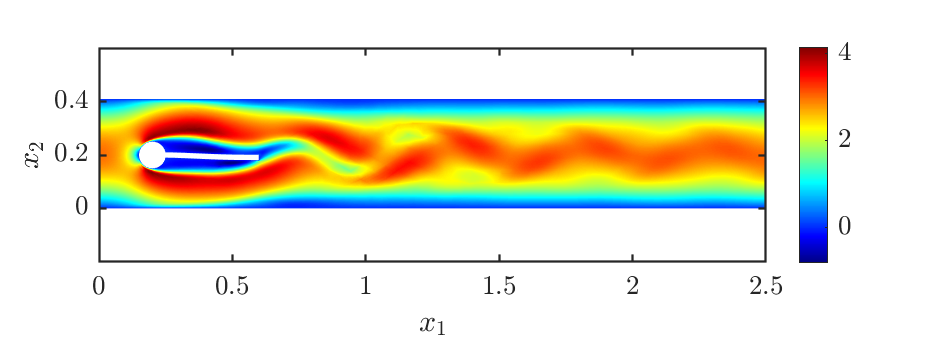}
        \end{minipage}%
    }}

    \vspace{-0.2em}

    {\subfloat[{$t=0.2\,\mathrm{s}$}]{%
        \begin{minipage}[b]{0.30\textwidth}
        \centering
        \includegraphics[
        width=\linewidth,
        trim={10 0 20 0},
        clip
        ]{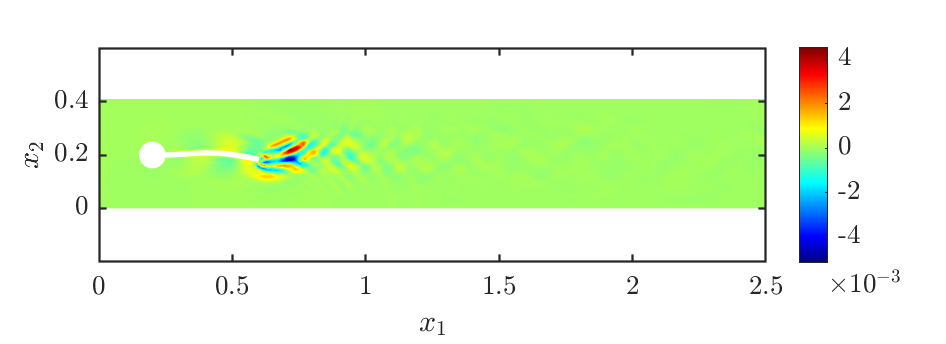}
        \end{minipage}%
    }}%
    \hspace{0.02\textwidth}%
    {\subfloat[{$t=3\,\mathrm{s}$}]{%
        \begin{minipage}[b]{0.30\textwidth}
        \centering
        \includegraphics[
        width=\linewidth,
        trim={10 0 20 0},
        clip
        ]{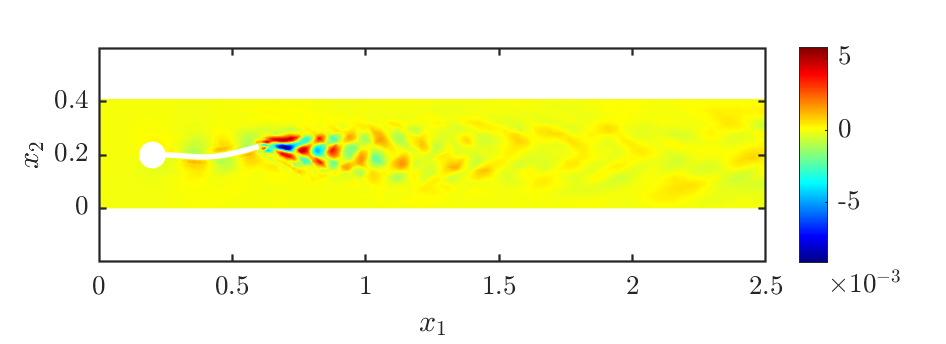}
        \end{minipage}%
    }}%
    \hspace{0.02\textwidth}%
    {\subfloat[{$t=12\,\mathrm{s}$}]{%
        \begin{minipage}[b]{0.30\textwidth}
        \centering
        \includegraphics[
        width=\linewidth,
        trim={10 0 20 0},
        clip
        ]{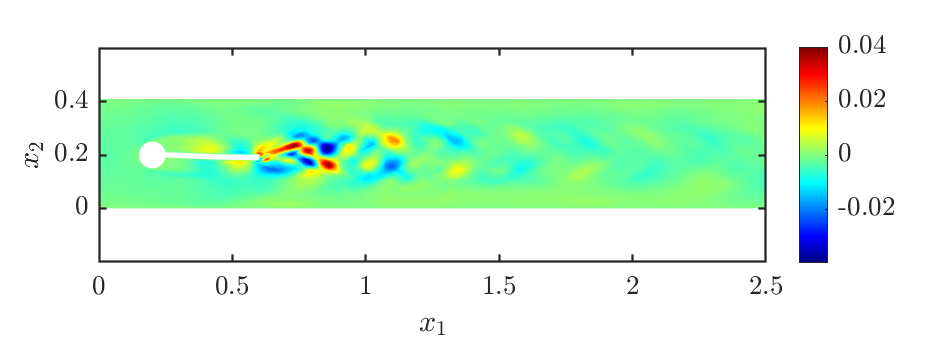}
        \end{minipage}%
    }}

    \caption{Turek problem; RQC-DMDc prediction of the streamwise velocity
    (top) and pointwise error relative to the full-order solution (bottom).}
    \label{fig:turek_rqc_dmdc_field_error}
\end{figure}

% ============================================================
% Relative error
% ============================================================

\begin{figure}[H]
    \centering

    \includegraphics[width=0.43\textwidth]
    {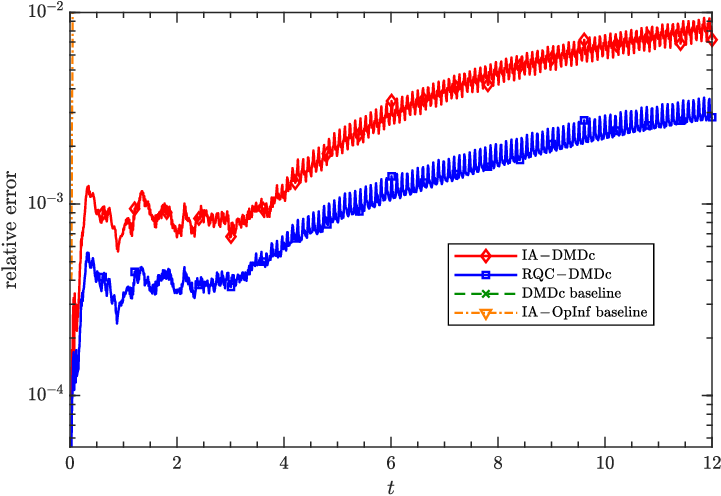}

    \caption{Turek problem; relative error of the streamwise velocity. The
    training interval ends at $t=3\,\mathrm{s}$.}

    \label{fig:turek_w_re}
\end{figure}

% ============================================================
% Beam-tip vertical displacement
% ============================================================

\begin{figure}[H]
    \centering

    {\subfloat[{Full interval, $t\in(0,12]\,\mathrm{s}$}]{%
        \begin{minipage}[t]{0.39\textwidth}
        \centering
        \includegraphics[width=\linewidth]
        {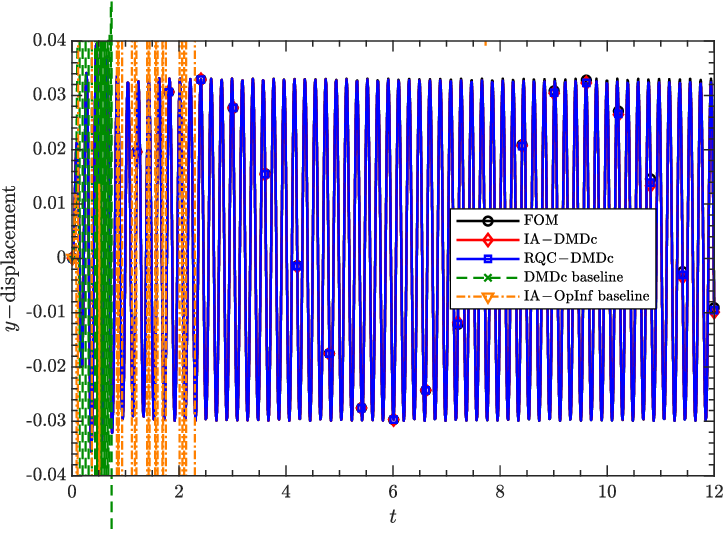}
        \end{minipage}%
    }}
    \hspace{0.04\textwidth}
    {\subfloat[{Final $0.5\,\mathrm{s}$, $t\in[11.5,12]\,\mathrm{s}$}]{%
        \begin{minipage}[t]{0.39\textwidth}
        \centering
        \includegraphics[width=\linewidth]
        {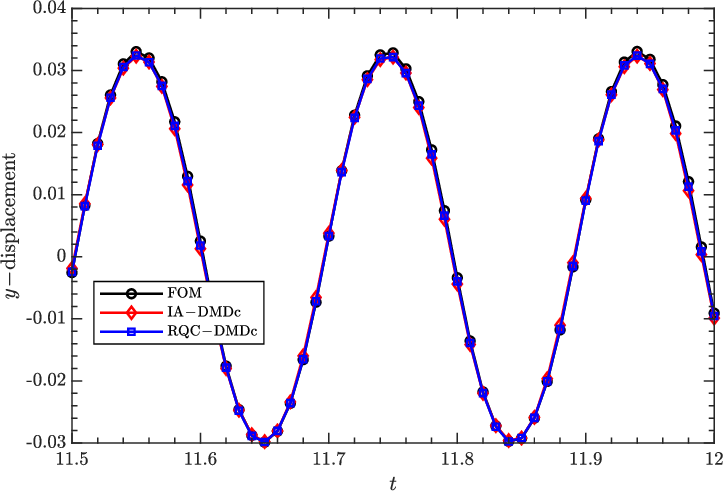}
        \end{minipage}%
    }}

    \caption{Turek problem; vertical displacement $d_{\mathrm A,\mathrm y}$ of the beam-tip
    monitoring point $A$ over the full interval and the final
    $0.5\,\mathrm{s}$.}

    \label{fig:turek_y_displacement}
\end{figure}

% ============================================================
% Drag force
% ============================================================

\begin{figure}[H]
    \centering

    {\subfloat[{Full interval, $t\in(0,12]\,\mathrm{s}$}]{%
        \begin{minipage}[t]{0.39\textwidth}
        \centering
        \includegraphics[width=\linewidth]
        {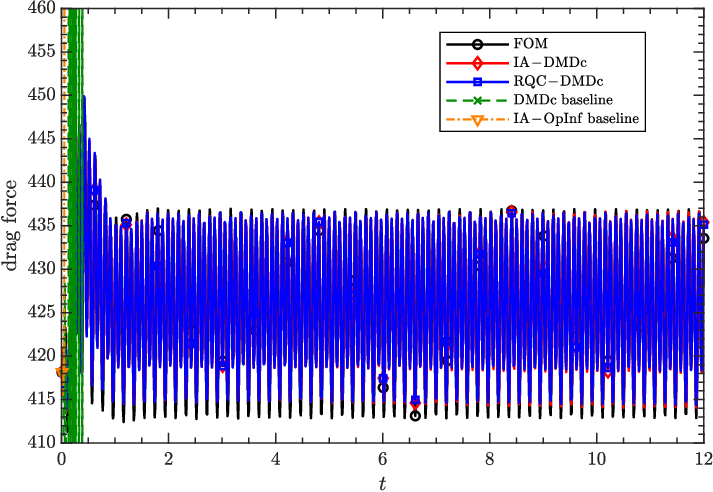}
        \end{minipage}%
    }}
    \hspace{0.04\textwidth}
    {\subfloat[{Final $0.5\,\mathrm{s}$, $t\in[11.5,12]\,\mathrm{s}$}]{%
        \begin{minipage}[t]{0.39\textwidth}
        \centering
        \includegraphics[width=\linewidth]
        {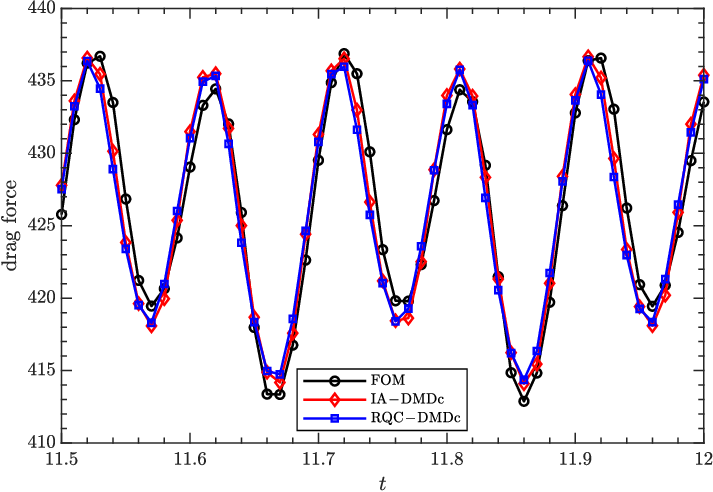}
        \end{minipage}%
    }}

    \caption{Turek problem; drag force exerted by the fluid on the
    cylinder--beam body over the full interval and the final
    $0.5\,\mathrm{s}$.}

    \label{fig:turek_drag_force}
\end{figure}

% ============================================================
% Lift force
% ============================================================

\begin{figure}[H]
    \centering

    {\subfloat[{Full interval, $t\in(0,12]\,\mathrm{s}$}]{%
        \begin{minipage}[t]{0.39\textwidth}
        \centering
        \includegraphics[width=\linewidth]
        {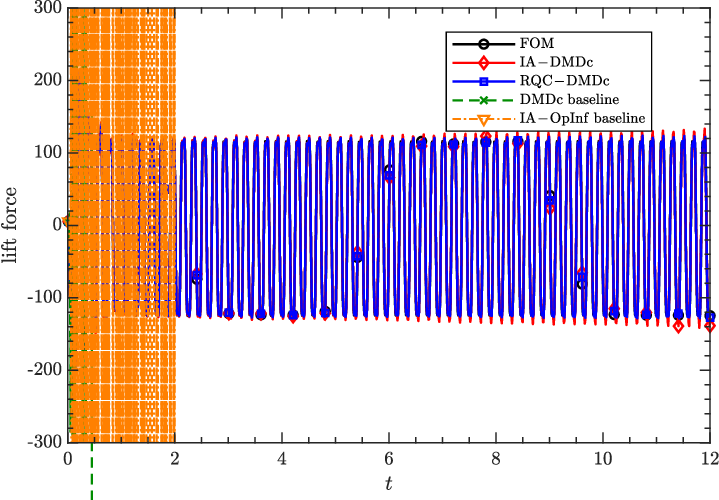}
        \end{minipage}%
    }}
    \hspace{0.04\textwidth}
    {\subfloat[{Final $0.5\,\mathrm{s}$, $t\in[11.5,12]\,\mathrm{s}$}]{%
        \begin{minipage}[t]{0.39\textwidth}
        \centering
        \includegraphics[width=\linewidth]
        {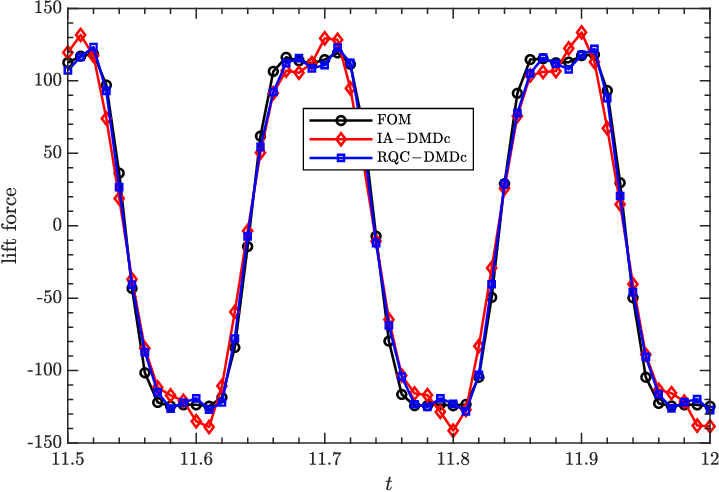}
        \end{minipage}%
    }}

    \caption{Turek problem; lift force exerted by the fluid on the
    cylinder--beam body over the full interval and the final
    $0.5\,\mathrm{s}$.}

    \label{fig:turek_lift_force}
\end{figure}

% ============================================================
% Fluid and solid kinetic energies
% ============================================================

\begin{figure}[H]
    \centering

    {\subfloat[{Fluid kinetic energy, $t\in(0,12]\,\mathrm{s}$}\label{fig:turek_fluid_kinetic_energy}]{%
        \begin{minipage}[t]{0.39\textwidth}
        \centering
        \includegraphics[width=\linewidth]
        {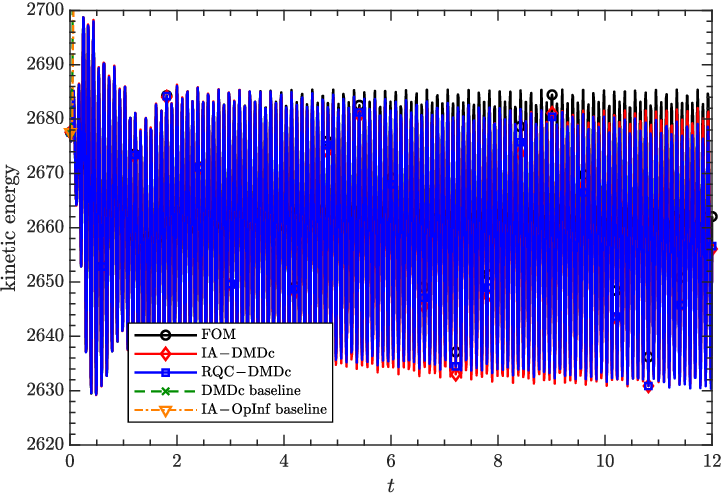}
        \end{minipage}%
    }}
    \hspace{0.04\textwidth}
    {\subfloat[{Fluid kinetic energy, $t\in[11.5,12]\,\mathrm{s}$}]{%
        \begin{minipage}[t]{0.39\textwidth}
        \centering
        \includegraphics[width=\linewidth]
        {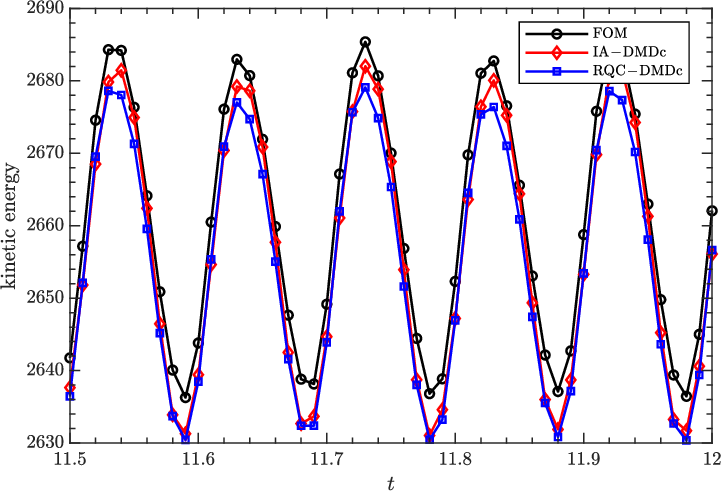}
        \end{minipage}%
    }}

    \vspace{0.2em}

    {\subfloat[{Solid kinetic energy, $t\in(0,12]\,\mathrm{s}$}\label{fig:turek_solid_kinetic_energy}]{%
        \begin{minipage}[t]{0.39\textwidth}
        \centering
        \includegraphics[width=\linewidth]
        {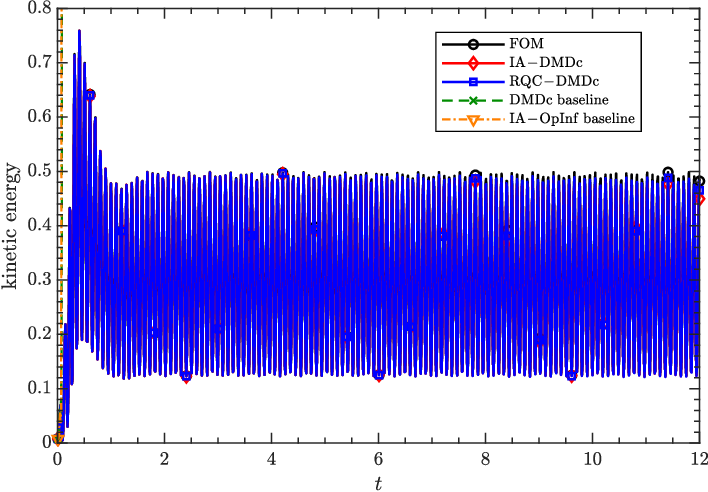}
        \end{minipage}%
    }}
    \hspace{0.04\textwidth}
    {\subfloat[{Solid kinetic energy, $t\in[11.5,12]\,\mathrm{s}$}]{%
        \begin{minipage}[t]{0.39\textwidth}
        \centering
        \includegraphics[width=\linewidth]
        {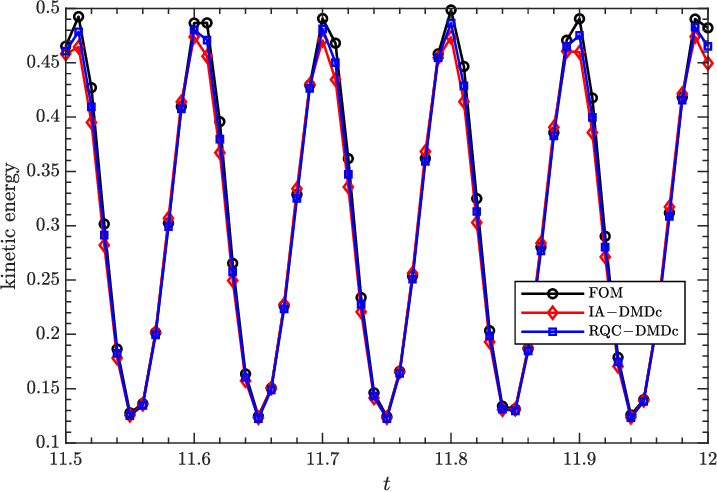}
        \end{minipage}%
    }}

    \caption{Turek problem; fluid and solid kinetic energies over the full
    interval (left) and the final $0.5\,\mathrm{s}$ (right).}

    \label{fig:turek_kinetic_energy}
\end{figure}

Figures~\ref{fig:turek_dmdc_field_error} and~\ref{fig:turek_rqc_dmdc_field_error} show that both
IA-DMDc and RQC-DMDc accurately reproduce the developing wake during training and retain the principal wake structures at \(t=12\,\mathrm{s}\).
The pointwise errors are concentrated near the beam tip and in the downstream wake, with smaller discrepancies for RQC-DMDc. Figure~\ref{fig:turek_w_re} provides a quantitative comparison:
the streamwise-velocity relative error of IA-DMDc remains below \(10^{-2}\) throughout the reported interval, demonstrating accurate prediction beyond the training window.
RQC-DMDc further reduces this error, although both models exhibit gradual error growth during extrapolation. 

Both models closely reproduce the dominant beam-tip oscillation in Figure~\ref{fig:turek_y_displacement}. The final \(0.5\,\mathrm{s}\) shows good agreement with the
full-order displacement in both amplitude and phase. The drag histories in Figure~\ref{fig:turek_drag_force} are also similar for the two reduced models and capture the mean force and the main oscillatory response, with small discrepancies near some extrema. 
The lift histories in Figure~\ref{fig:turek_lift_force} show a clearer distinction. IA-DMDc captures the dominant oscillation but exhibits modest overshoots near positive and negative extrema, while RQC-DMDc follows the reference lift more closely. The fluid and solid kinetic energies in Figure~\ref{fig:turek_kinetic_energy} confirm that both models retain the reference oscillation pattern over the prediction interval. In the final \(0.5\,\mathrm{s}\), both models slightly underestimate the fluid kinetic energy, with IA-DMDc reproducing the upper envelope more closely. For the solid kinetic energy, the two predictions remain close, and RQC-DMDc gives slightly better agreement near the peaks. 
The DMDc and IA-OpInf baselines become unstable early in the calculation and do not reproduce the out-of-sample coupled oscillations. 
Both interface-aware models remain bounded throughout the simulation and accurately capture the dominant fluid--structure response. This example demonstrates the predictive capability of linear IA-DMDc, while the regularized quadratic enrichment provides additional accuracy most clearly in the velocity field and lift history.

\section{Conclusions}
\label{sec:conclusions}

We developed an interface-aware reduced-order framework for optimization-based partitioned multiphysics systems. The prescribed interface-input operators are retained during local identification, and the common interface control is recomputed online by minimizing the physical interface mismatch.
The formulation accommodates ODE and DAE subproblems. For a DAE, tangent coordinates carry the differential evolution, while a complementary normal balance reconstructs the identifiable pressure algebraically.
We proposed RQC-DMDc to enrich the local model through jointly identified linear and compact quadratic operators. Feature normalization and block-dependent regularization control the enlarged regression without changing the prescribed interface action or identified DAE geometry. The resulting nonlinear reduced models are incorporated into the same implicit coupling formulation.

The numerical examples demonstrate the interface-aware formulation across scalar transmission and fluid--structure interaction problems. Across these examples, retaining the prescribed interface-input operator improves the robustness of the reduced prediction relative to the jointly identified DMDc baseline. 
IA-DMDc and RQC-DMDc give comparable field and energy predictions for the short-time Fisher--KPP problem and the manufactured Astorino--Grandmont solution. RQC-DMDc reduces velocity errors for the freely oscillating cylinder and the Turek problem, with additional improvements in selected force and energy histories.

Future work will consider parametric and multi-trajectory identification, data-driven or adaptive selection of stabilization and regularization parameters, structure-preserving nonlinear parameterizations, and three-dimensional applications with moving geometries.

\bibliography{references}

\end{document}